\documentclass[11pt, reqno]{amsart}
\usepackage{amsmath}
\usepackage{amssymb}
\usepackage{amsfonts}
\usepackage{dsfont}
\usepackage[usenames]{color}
\usepackage{version}

 \newtheorem{theorem}{Theorem}[section]
 \newtheorem{corollary}[theorem]{Corollary}
 \newtheorem{lemma}[theorem]{Lemma}
 \newtheorem{proposition}[theorem]{Proposition}
\newtheorem{question}[theorem]{Question}

\newtheorem{definition}[theorem]{Definition}
\newtheorem{example}[theorem]{Example}
\newtheorem{remark}[theorem]{Remark}
\newtheorem{fact*}{Fact}

\DeclareMathOperator{\spa}{span}

\newcommand\1{\boldsymbol {1}}
\newcommand\be{\begin{equation}}
\newcommand\ee{\end{equation}}
\newcommand\m{\mathcal{M}}

\newcommand\lam{\lambda}

\newcommand{\sad}{\mathcal{SA}_d}

\renewcommand{\L}{\mathcal{L}}
\newcommand{\M}{\mathcal{M}}
\newcommand{\N}{\mathcal{N}}
\newcommand{\T}{\mathbb{T}}

\newcommand{\D}{\mathbb{D}}
\newcommand{\C}{\mathbb{C}}

\newcommand{\nt}{\stackrel{\mathrm {nt}}{\to}}

\newcommand{\ip}[2]{\left\langle #1, #2 \right\rangle}

\newcommand{\vp}{\varphi}
\newcommand{\ph}{\varphi}

\newcommand\ga{\gamma}
\newcommand\de{\delta}

\newcommand\la{\lambda}

\newcommand\beq{\begin{equation}}

\newcommand\eeq{\end{equation}}
\newcommand{\inv}{^{-1}}
\newcommand\df{\stackrel{\rm def}{=}}

\newcommand\bbm{\begin{bmatrix}}
\newcommand\ebm{\end{bmatrix}}
\newcommand\bpm{\begin{pmatrix}}
\newcommand\epm{\end{pmatrix}}
\newcommand{\vectwo}[2]
{
   \begin{pmatrix} #1 \\ #2 \end{pmatrix}
}
\newcommand\half{\tfrac 12}

\newcommand\nin\noindent
\numberwithin{equation}{section}

\begin{document}

\title{Boundary behavior of functions in the Schur-Agler class of the polydisc}

\author{Jim Agler}
\address{Department of Mathematics, University of California at San Diego, CA \textup{92103}, USA}
\email{jagler@ucsd.edu}

\author{Connor Evans}
\address{School of Mathematics,  Statistics and Physics, Newcastle University, Newcastle upon Tyne NE\textup{1} \textup{7}RU, U.K.}
\email{CEvansMathematics@outlook.com}

\author{Zinaida A. Lykova}
\address{School of Mathematics,  Statistics and Physics, Newcastle University, Newcastle upon Tyne	NE\textup{1} \textup{7}RU, U.K.}
\email{Zinaida.Lykova@ncl.ac.uk}

\author{N. J. Young}
\address{School of Mathematics, Statistics and Physics, Newcastle University, Newcastle upon Tyne NE1 7RU, U.K.}
\email{Nicholas.Young@ncl.ac.uk}
%    \date is required; it is the date received by the editor.

\date{17th July, 2026}

\keywords{Schur class, Schur-Agler class, polydisc, Carath\'eodory condition, Hilbert space models, directional derivative}

\subjclass[2020]{32A30, 32S05, 47A56, 47A57 }
%32A30: Other generalizations of function theory of one complex variable (should also be assigned at least one classification number from Section 30)
%32S05: Local singularities
%30E20: Integration, integrals of Cauchy type, integral representations of analytic functions
%47B25: Symmetric and selfadjoint operators (unbounded)
%47A10: Spectrum, resolvent
%47A56: Functions whose values are linear operators (operator and matrix valued functions, etc., including analytic and meromorphic ones)
%47A57: Operator methods in interpolation, moment and extension problems
\thanks{Partially supported by National Science Foundation Grants
	DMS 1361720 and 1665260, Heilbronn Institute for Mathematical Research (Focused Research Grant) and the Engineering and Physical Sciences Research Council grants EP/N03242X/1 and EP/T517914/1. } 

\begin{abstract}
We describe a generalization of the notion of a Hilbert space model of a function $\varphi$ in the Schur-Agler class of the polydisc.  This generalization is well adapted to the investigation of boundary behavior of $\varphi$ at a mild singularity $\tau$ on the $d$-torus.  We prove the existence of a generalized model with an enhanced continuity property at such a singularity $\tau$. We use this result to prove the directional differentiability of a function $\varphi$ in the Schur-Agler class at a singular point on the  $d$-torus for which the Carath\'eodory condition holds and to calculate the corresponding directional derivative. The results of this paper extend to the polydisc $\mathbb{D}^d$ results of Agler, McCarthy, Tully-Doyle and Young which generalized to the bidisc the classical Julia-Wolff-Carath\'{e}odory theorem about analytic self-maps of $\mathbb{D}$. 
\end{abstract}

\maketitle
\tableofcontents
\section{Introduction}\label{intro}
In this paper we prove the existence of directional derivatives of certain analytic functions $\ph$ on the polydisc $\D^d$, where $d$ is a positive integer, at points of the boundary of $\D^d$ where a weak regularity property is satisfied. To do this we use a generalization of the notion of a Hilbert-space model of $\ph$.
Our results can be considered as higher-dimensional analogs of the  Julia-Carath\'{e}odory\footnote{In the Julia-Wolff-Carath\'eodory theorem, the contribution of Wolff was about attractive fixed points of analytic self-maps of the disc.  In this paper we are not concerned with self-maps of a domain nor with fixed points, so we shall refer to the classical result which we are generalizing as the ``Julia-Carath\'eodory theorem".}  theorem about analytic self-maps of $\D$. First we generalize theorems of Agler, McCarthy and Young \cite{AMYcara} on Carath\'eodory theorems for the bidisc via Hilbert space methods, and secondly we extend theorems proved by Agler, Tully-Doyle and Young in \cite{ATY2} on the boundary behavior of analytic functions on the bidisc.
\begin{definition}{\normalfont{\cite{Ag1}}}\label{modelpolydefinto}
Let $\varphi:\mathbb{D}^{d}\rightarrow\mathbb{C}$ be analytic. A {\em{model}} of $\varphi$ on $\mathbb{D}^{d}$ is a pair $(\mathcal{L},v)$, where $\mathcal{L}$ is a separable complex Hilbert space with an orthogonal decomposition $\mathcal{L}=\mathcal{L}_{1}\oplus\hdots\oplus\mathcal{L}_{d}$ and $v : \mathbb{D}^{d} \rightarrow \mathcal{L}$ is an analytic map such that, for all $\lambda, \mu \in \mathbb{D}^{d}$,
\begin{equation}\label{modeleqn}
1-\overline{\varphi(\mu)}\varphi(\lambda) = \Big\langle{\big(1_{\mathcal{L}}-\mu_{P}^{*}\lambda_{P}\big)v({\lambda}),v({\mu})} \Big\rangle_{\mathcal{L}},
\end{equation}
where $P$ is the $d$-tuple $P=(P_{1},\hdots,P_{d})$ on $\mathcal{L}$, for $j=1,\hdots,d$, $P_{j}:\mathcal{L}\rightarrow\mathcal{L}$ is the orthogonal projection onto $\mathcal{L}_{j}$,  and $\lambda_{P}$ denotes the operator $\lambda_{1}P_{1}+\hdots+\lambda_{d}P_{d}$.
\end{definition} 

The {\em Schur class} $\mathcal{S}(U)$ of a domain $U$ in $\mathbb{C}^d$ is defined to be the set of analytic functions $\ph$ on $U$ such that $|\ph(z)|\leq 1$ for all $z\in U$. The following result from \cite{Ag1} explains why models have been an effective tool for the study of the Schur class functions on  the bidisc.
\begin{theorem}\label{card2}
Let $\ph:\D^2 \to \C$ be analytic.  The following conditions are equivalent.
\begin{enumerate}
\item $\ph\in\mathcal{S}(\D^2)$;
\item $\ph$ has a model;
\item $\ph$ satisfies Ando's inequality, that is,
for all commuting pairs $(T_1,T_2)$ of contractions on a complex Hilbert space $\mathcal{H}$,
\[
\sup_{r<1}\|\ph(rT_1,rT_2)\|_{\mathcal{B}(\mathcal{H})} \leq \sup_{\lambda \in \D^2} | \ph(\lambda) |.
\]
\end{enumerate}
\end{theorem}
In this paper, for any complex Hilbert space $\mathcal{H}$, we denote by $\mathcal{B(\mathcal{H})}$ the Banach algebra of bounded linear operators $T: \mathcal{H}\rightarrow \mathcal{H}$, with composition as multiplication and the operator norm. 

Agler also proved in \cite{Ag1} that, for any positive integer $d$ and for any analytic function $\ph$ on $\D^d$, $\ph$ has a model if and only if $\ph$ satisfies  the $d$-dimensional analog of von Neumann's inequality, that is, conditions (2) and (3) of Theorem \ref{card2}  (with the obvious modifications for functions of $d$ variables)  are equivalent for all $\ph\in\mathcal{S}(\D^d)$ when $d \geq 2$.  It is well known that the $d$-dimensional analog of von Neumann's inequality fails for some functions from $ \mathcal{S}(\D^d)$, when $d \geq 3$, see \cite{Varopouloscite,crabbdavie}.
That is, conditions (1) and (3) of Theorem \ref{card2} are not equivalent for all $\ph\in\mathcal{S}(\D^d)$ when $d > 2$. 
Therefore,  not all functions $\ph\in\mathcal{S}(\D^d)$ have models when $d \geq 3$.
It is natural to consider the class of functions for which the  analog of von Neumann's inequality does hold.
\begin{definition}{\normalfont{\cite{Ag1}}}\label{schuraglerclassdefhere}
The {\em{Schur-Agler class}}, denoted by $\mathcal{SA}_{d}$, is the set of analytic functions $\varphi:\mathbb{D}^{d}\rightarrow \mathbb{C}$ such that, for all commuting $d$-tuples $(T_{1},\hdots, T_{d})$ of contractions acting on any complex Hilbert space $\mathcal{H}$,
$$\|\varphi\|_{(d,\infty)} \df \sup_{r<1} \| \varphi(rT_{1},\hdots, rT_{d}) \|_{\mathcal{B}(\mathcal{H})} \leq 1.$$
\end{definition}
\noindent It was proved in \cite{Ag1} that the analytic functions on $\D^d$ that have models are precisely the functions in 
 $\mathcal{SA}_d$. 

Since functions in the Schur-Agler class of the polydisc have models, 
it is natural to use models to elucidate the function-theoretic properties of such functions and, in particular, to study the boundary behavior of such functions near a mild singularity.   By a {\em singularity} of a function $\ph\in\mathcal{S}(\D^d)$ we mean a point $\tau$ in the topological boundary $\partial \D^d$ of the polydisc such that $\ph$ does not extend analytically to any neighborhood of $\tau$.  We say that  a function $\ph\in\mathcal{S}(\D^d)$ {\em satisfies the Carath\'eodory condition at a point $\tau$ in $\partial \D^d$ } if
\be \label{Carapoint}
\liminf_{\lam \to \tau, \lam \in \D^d} \frac{1-|\ph(\lam)|}{1-\|\lam\|_\infty} < \infty.
\ee
When this relation holds we also say that {\em $\tau$ is a carapoint for} $\ph$.
\begin{definition}{\normalfont{({Non-tangential approach on $\mathbb{D}^{d}$}})}\label{nontangetniallimitonDd-intro}
Let $S\subseteq \mathbb{D}^{d}$ and $\tau \in \partial\mathbb{D}^{d}$ be such that $\tau \in \overline{S}$. If there exists a constant $c>0$ such that, for all $\lambda \in S$,
$$ \|\lambda - \tau \|_{\infty}\leq c ( 1 - \|\lambda\|_{\infty}),$$
then we say $S$ approaches $\tau$ non-tangentially, denoted by $S {\overset{nt}\rightarrow}\tau$.
\end{definition}

The classical Julia-Carath\'{e}odory Theorem \cite{caratheodory}, which Carath\'eodory proved in 1929, states:
\begin{theorem}\label{caratheodedorytheorem}
    Let $\tau\in\mathbb{T}$ and let $\varphi\in\mathcal{S}(\mathbb{D})$. The following conditions are equivalent.
    \begin{enumerate}
        \item $\tau$ is a carapoint for $\varphi$;
          \item the quantity
        $$\sup_{\lambda\in S}\dfrac{1-\lvert \varphi(\lambda)\rvert}{1-\lvert\lambda\rvert}<\infty$$
        whenever $S\subseteq \mathbb{D}$ and $S{\overset{nt}\rightarrow}\tau$;
        \item $\varphi$ is directionally differentiable at $\tau$, that is, 
          there exists a complex number $\omega$ such that,
     whenever $\delta \in \C$ and
        $\tau+t\delta\in\mathbb{D}^{d}$ for sufficiently small positive $t$, 
$$\lim_{t\rightarrow 0^{+}} \dfrac{\varphi(\tau+t\delta)-\omega}{t}~~{\text{exists}}.$$
Moreover,  $|\omega|=1$;
        \item $\varphi$ is non-tangentially differentiable at $\tau$, that is, there exist $\omega\in\mathbb{C}$  and $\eta\in\mathbb{C}$ such that
$$\lim_{\substack{\lambda\rightarrow \tau \\ \lambda\in S}}\dfrac{\varphi(\lambda)-\omega-\eta (\lambda-\tau)}{\lvert\lambda-\tau\rvert}=0$$
whenever $S\subset \mathbb{D}$ and $S{\overset{nt}\rightarrow}\tau$. Moreover, $|\omega|=1$;
        \item $\varphi$ is non-tangentially continuously differentiable at $\tau$, that is, condition 
        {\rm (4)} holds and, in addition, 
        $$\lim_{\substack{\lambda\rightarrow \tau \\ \lambda\in S}}\varphi'(\lambda)=\eta$$
whenever $S\subseteq \mathbb{D}$ and $S{\overset{nt}\rightarrow}\tau$.
        \end{enumerate}
\end{theorem}
The gist of Theorem \ref{caratheodedorytheorem} is
that if a function satisfies the Carath\'{e}odory condition at a point $\tau \in \T$ then it satisfies some apparently stronger regularity conditions there, in particular, it has an “angular derivative” at $\tau$ (meaning that condition (4) holds). 

Carath\'eodory's result has been extended to functions of several variables by numerous authors, notably by K. W\l odarczyk \cite{woody}, W. Rudin \cite{rudin}, F. Jafari \cite{jaf93},  M. Abate \cite{Abate} and by Agler, McCarthy and Young \cite{AMYcara}. A version of the Julia-Carath\'eodory theorem for functions in the Schur-Agler class of the polydisc was proved by J.A. Ball and V. Bolotnikov in \cite[Lemma 2.1]{BallandBolot1}.

In particular, in \cite{AMYcara} Theorem \ref{caratheodedorytheorem}  was extended to analytic functions on the bidisc $\mathbb{D}^{2}$. For $\ph \in \mathcal{S}(\mathbb{D}^{2})$ and for 
$\tau \in \partial \D^2$, the authors  proved that if $\ph$ satisfies the Carath\'eodory condition at $\tau$ then  the directional derivatives $D_{\delta}\ph(\tau)$ exist for all appropriate directions $\de \in\C^2$. 
They also give an example of a rational inner function $\psi \in \mathcal{S}(\mathbb{D}^{2})$, 
$\psi:\D^2 \to \D$, given by the formula
\[
\psi(\lam)= \frac{2\lam_1\lam_2 -\lam_1-\lam_2}{2-\lam_1-\lam_2},
\]
which satisfies Carath\'{e}odory’s condition at the point $\chi:= (1,1)$ but
is not non-tangentially differentiable at  $\chi$. In \cite{ATY2} Agler, Tully-Doyle and  Young  introduced the notion of a ``generalized model", which enabled them  to present an alternative, more algebraic, proof of the existence of directional derivatives at $\tau$.

In this paper we analyse the boundary behavior of functions in the Schur-Agler class $\mathcal{SA}_{d}$ of the polydisc $\mathbb{D}^{d}$ of any dimension $d\geq 2$. We exploit Hilbert space methods which were introduced in \cite{ATY2} for analytic functions on the bidisc $\mathbb{D}^{2}$. 
We shall show below that, for a function $\ph\in\mathcal{SA}_d$, the Carath\'eodory condition at a point $\tau\in\T^d$ implies that $\ph$ has a non-tangential limit at $\tau$, and that $\varphi$ is directionally differentiable at $\tau$.  

  An essential component of our approach will be the notion of a generalized model, which we now define. 
\begin{definition}\label{inneroperator}
Let $\mathcal{H}$ be a Hilbert space and let 
$I : \mathbb{D}^{d} \mapsto \mathcal{B}(\mathcal{H})$
be an operator-valued analytic function. We say that $I$ is {\em{inner}} if, for almost all $\lambda\in\mathbb{T}^{d}$ with respect to Lebesgue measure, the radial limit 
$$\lim_{r\rightarrow 1^{-}}I(r\lambda)$$ 
exists with respect to the strong operator topology and is an isometric operator. 
\end{definition}
The following definition, an extension of Definition \ref{modelpolydefinto} in which $\lam_P$ is replaced by a general analytic operator-valued function $I$, will be central to our efforts.
\begin{definition}\label{genmodel}
Let $\vp:\D^d\to\C$ be analytic. The triple $(\mathcal M, u, I)$ is a {\em generalized model of $\vp$} if 
\begin{enumerate}
 \item $\M$ is a  separable complex Hilbert space, 
 \item $u: \D^d \to \M$ is analytic on $\mathbb{D}^{d}$, 
 \item $I:\mathbb{D}^{d}\rightarrow\mathcal{B}(\mathcal{M})$ is analytic on $\mathbb{D}^{d}$ and $\|I(\lambda)\|_{\mathcal{B(M)}}<1$ for all $\lambda\in\mathbb{D}^{d}$, and 
\item the equation
\begin{equation}\label{modeleqn-I}
1-\overline{\varphi(\mu)}\varphi(\lambda) = \Big\langle{\big(1_{\mathcal{M}}-I(\mu)^{*}I(\lambda)\big)u(\lambda),u(\mu)} \Big\rangle_{\mathcal{M}}
\end{equation}
holds for all $\lambda, \mu \in \D^d$.
\end{enumerate}
A generalized model $(\mathcal M, u, I)$ is {\em inner} if $I(\cdot)$ is an inner $\mathcal{B(\M)}$\text{-valued function on}~$\mathbb{D}^{d}$.
\end{definition}
 This simple relaxation of the definition of model enables us to concentrate information about singular behavior of $\ph$ in the inner function $I(\cdot)$ rather than the model function $u$.  This fact proved efficacious in \cite{ATY2} for the description of directional derivatives at a carapoint $\tau$ for $\ph \in\mathcal{S}(\D^2)$ since  one obtains a generalized model in which the function $u:\D^d \to \m$ extends continuously to the singular point $\tau$.
True, the cost is an inner operator-valued function $I(\cdot)$ which is, in principle, more complicated than the function $\lambda_P$. However, it transpires that one can construct a generalized model in which the inner operator-valued function $I(\cdot)$ is given by an explicit (and fairly simple) formula, see Theorem \ref{mainExistence-intro}.                 
 In this paper we carry out a similar programme to that of \cite{ATY2} in dimensions $d \geq 2$.
 We formalise the appropriate continuity property of $u$ in the following definition.
\begin{definition}{\normalfont{\cite[Definition 3.4]{ATY2}({\bf{$C$-point of the model}})}}
Let $(\mathcal{M},u,I)$ be a generalized model of a function $\varphi\in\mathcal{SA}_{d}$. A point $\tau\in \partial\mathbb{D}^{d}$ is a {\em{$C$-point of the model}} $(\mathcal{M},u,I)$ if, for every subset $S$ of $\mathbb{D}^{d}$ that approaches $\tau$ non-tangentially, $u$ extends continuously to $S\cup\{\tau\}$, with respect to the norm topology on $\mathcal{M}$. 
\end{definition}
In contrast, if $\tau$ is a carapoint of $\ph$ and $(\M,u)$ is a model of $\ph$ in the sense of Definition \ref{modelpolydefinto}, then $\tau$ is a ``$B$-point" of $(\M,u)$ in the following sense.
\begin{definition}{\normalfont{\cite[Definition 3.4]{ATY2}({\bf{$B$-point of the model}})}} Let $(\mathcal{L},u)$ be a  model of a function $\varphi\in\mathcal{SA}_{d}$. A point $\tau\in \partial\mathbb{D}^{d}$ is a {\em{$B$-point of the model}} if $u$ is bounded on every subset of $\mathbb{D}^{d}$ that approaches $\tau$ non-tangentially. 
\end{definition}
We shall prove in Proposition \ref{prop5.2cara} that, if $\tau$ is a carapoint of $\ph\in\mathcal{SA}_{d}$ and $(\L,u)$ is a model of $\ph$, then $\tau $ is a $B$-point of the model, but $\tau$ may fail to be a $C$-point of the model.

We shall interpret $\mathbb{C}^{d}$ as a Banach algebra with co-ordinate-wise operations and identity $\mathds{1}=(1,1,\hdots,1)$ so that, for example, if $\lambda=(\lambda_{1},\hdots,\lambda_{d})$, then
$$\bigg(\dfrac{\mathds{1}}{\mathds{1}-\lambda}\bigg) = \bigg(\dfrac{1}{1-\lambda_{1}},\hdots, \dfrac{1}{1-\lambda_{d}}\bigg).$$
Following Definition \ref{Carapoint}, we define carapoints of {\em operator-valued} functions on $\D^d$ in the following way.
\begin{definition} 
Let $I$ be a contractive operator-valued analytic function on $\D^d$. We say that $\tau\in\T^d$ is a {\em carapoint} for $I$ if
\[
\liminf_{\la\to\tau}\frac{1-\|I(\la)\|}{1-\|\la\|_\infty} < \infty.
\]
\end{definition}

Consider a function $\varphi\in\mathcal{SA}_{d}$ having a carapoint at $\tau\in \T^d$.
The following definition formalises the notion of a generalized model of $\varphi$ for which
$\tau$ is a $C$-point of the model.
\begin{definition}\label{desingularized} 
Corresponding to a function $\varphi\in\mathcal{SA}_{d}$ and a carapoint $\tau\in\mathbb{T}^{d}$ of $\varphi$, we define a {\em{desingularized model of $\varphi$ relative to $\tau$}} to be a generalized model $(\mathcal{M},u,I)$ of $\varphi$ such that the following conditions hold:
 \begin{enumerate}
\item $\tau$ is a $C$-point for $(\mathcal{M},u,I)$;
\item $\tau$ is a carapoint for $I$ and $\lim_{\lambda{\overset{nt}\rightarrow}\tau}I(\lambda)=1_{\mathcal{M}}$;
\item $I:\mathbb{D}^{d}\rightarrow \mathcal{B(M)}$ is an inner $\mathcal{B(M)}$-valued function.
\end{enumerate}
\end{definition}
The main tool we shall use in this paper is the following theorem (see Theorem \ref{generalizedmodeltheorem}). 
\begin{theorem}\label{mainExistence-intro}
Let $\tau \in\mathbb{T}^{d}$ be a carapoint for $\varphi \in\mathcal{SA}_{d}$. There exists a 
desingularized model $(\mathcal{M},u,I)$ of $\varphi$ relative to $\tau$. 
Furthermore, we can express $I$ in the form\footnote{Here the notation $\left(\frac{\mathds{1}}{\mathds{1}-\bar\tau\lam}\right)_Y\inv$ denotes $\left( \frac{1}{1-\bar\tau_1\lam_1}Y_1+\dots+\frac{1}{1-\bar\tau_d\lam_d}Y_d \right)\inv$.}
\begin{equation}\label{Ioperatorgeneralizedmodel-intro}
I(\lambda) = 1_{\mathcal{M}} - \bigg( \dfrac{\mathds{1}}{\mathds{1}-\overline{\tau}\lambda}\bigg)_{Y}^{-1}, \  \text{for} \ \lambda \in \mathbb{D}^{d},
\end{equation}
for some $d$-tuple of positive contractions $(Y_{1},\hdots,Y_{d})$ on $\mathcal{M}$ such that $\sum_{j=1}^{d}Y_{j}=1_{\mathcal{M}}$.  At every point $\lambda\in\mathbb{T}^{d}$ such that $\lambda_{j}\neq\tau_{j}$ for $j=1,\hdots,d$, $I(\cdot)$ extends analytically to $\lambda$ and $\|I(\lambda)\|_{\mathcal{B(M)}}=1$.
\end{theorem}

In Remark \ref{case-d=2} we compare the known results for $\ph\in\mathcal{S}(\D^2)$ from \cite{ATY2} with our Theorem \ref{mainExistence-intro}. In \cite[Theorem 3.6]{ATY2} the authors considered a function 
$\ph\in\mathcal{S}(\D^2)$ with a carapoint $\tau \in \T^2$  and used a model $(\mathcal{L},v)$ of $\varphi$ and a realization $(\alpha,\beta, \gamma, D)$ of $(\mathcal{L},v)$ to construct an inner generalized model $(\m,u,I_{ATY})$. We show that 
 $I=I_{ATY}$ and the two generalized models coincide.

We now address directional derivatives of $\varphi \in\mathcal{SA}_{d}$ at a carapoint $\tau \in\mathbb{T}^{d}$. The fact that $\tau$ is a carapoint for $\varphi \in\mathcal{SA}_{d}$ implies that the non-tangential limit at $\tau$,
$$\lim_{\lambda{\overset{nt}\rightarrow}\tau}\varphi(\lambda)=\omega,$$
exists, and its value $\omega$ satisfies $|\omega|=1$  (see Proposition \ref{phi-nt-omega}). 
When this non-tangential limit exists
we shall often use the notation $\varphi(\tau) = \lim_{\lambda{\overset{nt}\rightarrow}\tau}\varphi(\lambda)$.
For $\tau\in\mathbb{T}^{d}$, one can see that $\delta$ from $\mathbb{C}^{d}$ points into $\mathbb{D}^{d}$ at $\tau$ if and only if $-\delta \in\mathbb{C}^{d}_{+}(\tau)$, where
$$\mathbb{C}_{+}^{d}(\tau) = \big\{( z_{1},\hdots, z_{d})\in\mathbb{C}^{d} : \mathrm{Re}( z_{j}\overline{\tau}_{j})>0~~\text{for}~~j=1,\hdots,d\big\}.$$
\begin{definition}
Let $\varphi\in\mathcal{SA}_{d}$ and let $\tau\in\mathbb{T}^{d}$ be a carapoint for $\varphi$. We say that
$\varphi$ is {\em directionally differentiable} at $\tau$ if, whenever $-\delta \in\mathbb{C}^{d}_{+}(\tau)$,
$$\lim_{t\rightarrow 0^{+}} \dfrac{\varphi(\tau+t\delta)-\ph(\tau)}{t}~~{\text{exists}}.$$
When the limit exists we denote it by 
$D_{\delta}\varphi(\tau).$
\end{definition}
In Section \ref{directional_der} the following function-theoretic result is proved by the use of  desingularized models of $\varphi\in\mathcal{SA}_{d}$ (see Theorem \ref{directionalderivtheorem}). It is a generalization of the implication (1)$\Rightarrow$(3) in the classical Julia-Carath\'eodory Theorem. 
\begin{theorem}
If $\tau\in\mathbb{T}^{d}$ is a carapoint for $\varphi\in\mathcal{SA}_{d}$, then $\varphi$ is directionally differentiable at $\tau$.
\end{theorem}

Yet another dividend of desingularized models is an explicit formula for the directional derivative $D_{\delta}\varphi(\tau)$ of $\varphi$ at a carapoint $\tau$, as described in the following theorem.
\begin{theorem}\label{dir-derivative-tau-intro}
    Let $\tau\in\mathbb{T}^{d}$ be a carapoint for $\varphi\in\mathcal{SA}_{d}$. 
There exists an analytic function $h:\mathbb{C}_{+}^{d}(\tau)\rightarrow\mathbb{C}$ such that
    \begin{enumerate}
        \item $\mathrm{Re}(-h(z))>0$ for all $ z\in\mathbb{C}_{+}^{d}(\tau)$;
        \item \begin{equation}\label{hfunction-intro}
            h(\tau) = -\liminf_{\lambda{\rightarrow}\tau} \dfrac{1-\lvert \varphi(\lambda)\rvert}{1-\|\lambda \|_{\infty}};
        \end{equation}
        \item for all $ \delta\in\mathbb{C}_{+}^{d}(\tau)$, 
        \begin{equation}\label{directionalderivativeD-intro}
        D_{- \delta}\varphi(\tau) = \varphi(\tau)h( \delta).
\end{equation}
    \end{enumerate}
Moreover, we can express the function $h:\mathbb{C}_{+}^{d}(\tau)\rightarrow\mathbb{C}$ by
$$ h( z) = -\bigg\langle \bigg(\dfrac{\mathds{1}}{\overline{\tau} z}\bigg)_{Y}^{-1}u(\tau),u(\tau)\bigg\rangle_{\mathcal{M}} \ ~~\text{for all}~ z\in\mathbb{C}_{+}^{d}(\tau),$$  
for some desingularized model $(\mathcal{M},u,I)$ of $\varphi$ with $I$ of the form
\begin{equation}\label{Ioperatorgeneralizedmodel-intro-2}
I(\lambda) = 1_{\mathcal{M}} - \bigg( \dfrac{\mathds{1}}{\mathds{1}-\overline{\tau}\lambda}\bigg)_{Y}^{-1} \
\text{for} \ \lambda \in \mathbb{D}^{d},
\end{equation}
for some $d$-tuple of positive contractions $(Y_{1},\hdots,Y_{d})$ on $\mathcal{M}$ such that $\sum_{j=1}^{d}Y_{j}=1_{\mathcal{M}}$. 
\end{theorem}
 
For Hilbert spaces $\mathcal{H}$ and $\mathcal{K}$ we shall use notation that allows the identification of elements of the algebraic tensor product $\mathcal{H} \otimes\mathcal{K}$ with operators from $\mathcal{K}$ to $\mathcal{H}$. For $\gamma \in \mathcal{H}$ and $ \beta \in \mathcal{K}$, we define $\gamma \otimes \beta \in \mathcal{B}(\mathcal{K}, \mathcal{H})$ by the formula
\beq \label{u-otimes-v}
(\gamma \otimes \beta) (x) = \ip{x}{ \beta}_{\mathcal{K}} \gamma \;\; \text{for all} \; x \in \mathcal{K}.
\eeq

A well-known lurking isometry argument proceeds from a model $(\M,u)$ of a function $\ph\in \sad$ to a realization of $\ph$ \cite{Ag1}.  The identical argument applied to a {\em generalized} model $(\M, u, I)$ produces a {\em generalized} notion of realization, see Theorem \ref{genRealiz} in Section \ref{caratheodory_section}.
\begin{theorem} \label{genRealiz-intro}
 If $(\L, u, I)$ is a generalized model of $\vp \in \sad$ then there exist a Hilbert space $\M$ containing $\L$, a scalar $a \in \C$, vectors $\beta, \gamma \in \M$ and a linear operator $D:\M \to \M$ such that the operator
\beq \label{defL-intro}
L =  \begin{bmatrix}
  a & 1 \otimes \beta \\
  \gamma \otimes 1 & D
 \end{bmatrix}
\eeq
is unitary on $\C \oplus \M$ and, for all $\lambda \in \D^d$, 
\beq\label{Leq-intro}
 L \vectwo{1}{I(\lambda)u(\lambda)} = \vectwo{\vp(\lambda)}{u(\lambda)},
\eeq
and consequently, for all $\la\in\D^d$,
\beq\label{formphi-intro}
\ph(\la) = a+ \ip{I(\la)(1_{\M}-DI(\la))^{-1}\ga}{\beta}_{\M}.
\eeq
\end{theorem}
The ordered 4-tuple $(a, \beta, \gamma, D)$, as  in equation \eqref{defL-intro}, will be called a {\em realization} of the (generalized) model $(\L, u, I)$ of $\ph$ if $L$ is a contraction and equation \eqref{Leq-intro} holds.  It will be called a {\em unitary realization} if in addition $L$ is unitary on $\C\oplus \L$. 
 
The paper is organized as follows. 
In Section \ref{caratheodory_section} we discuss the Carath\'eodory condition and study properties of functions $\varphi \in \mathcal{SA}_{d}$ at $B$-points and at $C$-points. 
In Sections \ref{boundary-results} and \ref{B-point-realization}
we prove   numerous results concerning the boundary behavior of a function $\varphi \in \mathcal{SA}_{d}$. These results allow us to move seamlessly between results formulated in terms of functions, to their equivalents restated in terms of models. These re-statements make such functions more amenable to analysis. 
In these first three sections we follow and use methods of the paper ``A Carath\'{e}odory Theorem for the bidisc via Hilbert space methods'' by Agler, McCarthy and Young, \cite{AMYcara}. 
For the purpose of the present paper we need to generalize several definitions and theorems from \cite{AMYcara} from functions in $\mathcal{S}(\D^2)$ to functions in $\mathcal{SA}_{d}$, $d\geq 2$. The new definitions and statements are very similar to those in \cite{AMYcara}, but to establish them requires some new calculations.

In Sections \ref{interaction}, \ref{analytic_I} and \ref{Ioflambda} we develop the theory of the 
 operator-valued analytic function $I(\lambda)$ on $\mathbb{D}^{d}$, which appears in equation \eqref{Ioperatorgeneralizedmodel-intro} and which is essential for our main result (Theorem \ref{mainExistence-intro}). 
In particular we investigate conditions for the invertibility and analyticity of the operator-valued functions on $\mathbb{D}^{d}$ given by
$(\mathds{1}-\lambda)_{T}~\text{and}~\bigg(\dfrac{\mathds{1}}{\mathds{1}-\lambda}\bigg)_{T}$
for a $d$-tuple of operators $T=(T_{1},\hdots,T_{d})$ on $\mathcal{M}$. 
In Section \ref{main_result} we  prove the main result of our paper Theorem \ref{mainExistence-intro}.  Starting with a model $(\mathcal{L},v)$  of $\varphi\in\mathcal{SA}_{d}$ with a $B$-point $\tau \in \mathbb{T}^d$, we  arrive at a model $(\mathcal{M},u,I)$ of a more general type. It is called the desingularization of the model $(\mathcal{L},v)$ at $\tau$, see Definition \ref{desingularized}, and satisfies  additional conditions, in particular, $\tau$ is a $C$-point for $(\mathcal{M},u,I)$ and $\lim_{\lambda{\overset{nt}\rightarrow}\tau}I(\lambda)=1_{\mathcal{M}}$.
In Section \ref{directional_der} the model $(\mathcal{M},u,I)$ is utilised to prove the directional differentiability of a function $\varphi\in\mathcal{SA}_{d}$ at a carapoint $\tau$ for suitable directions into the polydisc, see Theorem \ref{dir-derivative-tau-intro}.
In Section \ref{example-in-D3} we give an example of a function $\ph_3\in\mathcal{SA}_{3}$
with a carapoint in $\T^3$ and in Section \ref{model3} we
present an explicit $9$-dimensional model $(\mathcal{L},v)$ of $\ph_3$.

\section{Special types of singularities for Schur-Agler class functions}\label{caratheodory_section}
In this section we follow the paper ``A Carath\'{e}odory Theorem for the bidisc via Hilbert space methods'' by Agler, McCarthy and Young, \cite{AMYcara}. The authors of \cite{AMYcara} proved a Julia-Carath\'{e}odory Theorem for the Schur class functions on $\mathbb{D}^{2}$.  To make this paper self-contained, in this section, we extend their results to the class $\mathcal{SA}_{d}$. We must first make suitable definitions of carapoints, $B$-points, $C$-points and angular gradients which permit the formulation of analogues of the main results of \cite{AMYcara} and allow the proofs to remain valid, with minimal changes. To investigate the boundary behaviour of a function $\varphi\in\mathcal{SA}_{d}$ by means of models, we require the definitions of carapoints, $B$-points and $C$-points of a model for the given  function $\varphi$.

We use the notation $\| \cdot \|_{\infty}$ to denote the $l^{\infty}$ norm on $\mathbb{C}^{d}$; that is, for all $\lambda\in\mathbb{C}^{d}$, 
$$\|\lambda \|_{\infty} = \max\{\lvert\lambda_{1}\rvert,\hdots, \lvert\lambda_{d}\rvert\}.$$
Although we mostly concern ourselves with points $\tau\in\mathbb{T}^{d}$, we include for completeness companion results for
$$\partial\mathbb{D}^{d}:= \big\{z\in\mathbb{C}^{d}:\lvert z_{j} \rvert \leq 1 ~\text{for}~j=1,\hdots, d~\text{and}~\lvert z_{j}\rvert=1~\text{for at least one}~j=1,\hdots, d \big\}.$$
\begin{definition}{\normalfont\label{JuliaQuotientOnDd}({\bf{Julia quotient for $\mathbb{D}^{d}$}})}
Let $\varphi \in \mathcal{S}(\mathbb{D}^{d})$. The Julia quotient of $\varphi$, denoted $J_{\varphi}$, is the function on $\mathbb{D}^{d}$ given by
$$J_{\varphi}({\lambda})=\dfrac{1-\lvert{\varphi(\lambda)}\rvert}{1-\|{\lambda}\|_{\infty}}.$$
\end{definition}
Note here and throughout, we use the superscript notation $\lambda^{(n)}$ to denote an element of a sequence $(\lambda^{(n)})_{n=1}^{\infty}$. If $(\lambda^{(n)})_{n=1}^{\infty}$ is a sequence in $\mathbb{D}^{d}$, then we use the subscript notation $\lambda^{(n)}_{j}$ to denote the $j$\textsuperscript{th} component of the point $\lambda^{(n)}$.

\begin{definition}{\normalfont{({\bf{Non-tangential approach on $\mathbb{D}^{d}$}}})}\label{nontangetniallimitonDd}
Let $S\subseteq \mathbb{D}^{d}$ and $\tau \in \partial\mathbb{D}^{d}$ be such that $\tau \in \overline{S}$. If there exists a constant $c>0$ such that, for all $\lambda \in S$,
$$ \|\lambda - \tau \|_{\infty}\leq c ( 1 - \|\lambda\|_{\infty}),$$
then we say that $S$ approaches $\tau$ non-tangentially, denoted $S {\overset{nt}\rightarrow}\tau$. Similarly, if $(\lambda^{(n)})_{n\in\mathbb{N}}$ is a sequence in $\mathbb{D}^{d}$ satisfying $\lambda^{(n)} \rightarrow \tau$ and, for some $c>0$,
\begin{equation}\label{nontangentialdefpoly}
 \| \lambda^{(n)} - \tau \|_{\infty} \leq c ( 1 - \|\lambda^{(n)}\|_{\infty} ),
 \end{equation}
for all $n\geq 1$, then we say that $\lambda^{(n)}$ has non-tangential limit $\tau$. We denote this relation by $\lambda^{(n)} {\overset{nt}\rightarrow} \tau$.
\end{definition}
The notion of carapoints essentially remains the same, but with $\lambda\in\mathbb{D}^{d}$ instead of $\lambda\in\mathbb{D}$ and now the Julia quotient of Definition \ref{JuliaQuotientOnDd}.
\begin{definition}{\normalfont({\bf{Carapoints for $\varphi$ on $\mathbb{D}^{d}$}})}\label{generalizedcara}
Let $\varphi\in\mathcal{S}(\mathbb{D}^{d})$ and let $\tau\in\partial\mathbb{D}^{d}$. We say that $\tau$ is a {\em{carapoint}} for $\varphi$ if
$$\liminf_{\substack{\lambda\rightarrow\tau \\ \lambda\in\mathbb{D}^{d}}} J_{\varphi}(\lambda)<\infty,$$
or equivalently, there exists a sequence $\{\lambda^{(n)}\}$ in $\mathbb{D}^{d}$ such that $\lambda^{(n)}\rightarrow \tau$ and 
$$\sup_{n} J_{\varphi}(\lambda^{(n)})<\infty.$$
\end{definition}
The following definition is the natural analogue for a function on the polydisc of a non-tangentially differentiable function on the unit disc. 
\begin{definition}{\normalfont({\bf{Angular gradient}})}\label{angularderivative}
Let $\varphi\in\mathcal{S}(\mathbb{D}^{d})$ and let $\tau\in\mathbb{T}^{d}$ be a carapoint of $\varphi$. We say that $\varphi$ has an {\em{angular gradient}} at $\tau\in\mathbb{C}^{d}$ if there exist $\omega\in\mathbb{C}$ of unit modulus and $\eta\in\mathbb{C}^{d}$ such that
$$\lim_{\substack{\lambda\rightarrow \tau \\ \lambda\in S}}\dfrac{\varphi(\lambda)-\omega-\eta\cdot (\lambda-\tau)}{\|\lambda-\tau\|_{\infty}}=0$$
whenever $S\subset \mathbb{D}^{d}$ and $S{\overset{nt}\rightarrow}\tau$. Here
$$\eta \cdot (\lambda-\tau)=\eta_{1}(\lambda_{1}-\tau_{1})+\hdots+\eta_{d}(\lambda_{d}-\tau_{d})$$
and the {\em angular gradient} of $\varphi$ at $\tau$ is the vector $\begin{bmatrix}
    \eta_{1}\\
    \vdots\\
    \eta_{d}
\end{bmatrix}.$
\end{definition}
 Let $\varphi\in\mathcal{SA}_{d}$, $\tau\in\partial\mathbb{D}^{d}$, and let $\{\lambda^{(n)}\}$ be a sequence in $\mathbb{D}^{d}$ such that $\lambda^{(n)}{\overset{nt}\rightarrow}\tau$. We will prove in  
Proposition \ref{prop5.2cara} that $\tau$ is a carapoint of $\varphi$ if and only if,  for every model $(\mathcal{M},u)$ of $\varphi$,  $\tau$ is a $B$-point of the model $(\mathcal{M},u)$ and  $ \lim_{n \to \infty} |\varphi(\lambda^{(n)})| =1$.  We follow  the authors of \cite{AMYcara} and call a carapoint of $\varphi$  a $B$-point of $\varphi$.

\begin{definition}{\normalfont{({\bf{$B$-point for $\varphi$}})}}\label{B-point-phi}
Let $\varphi\in\mathcal{S}(\mathbb{D}^{d})$ and let $\tau\in \partial\mathbb{D}^{d}$. We say that $\tau$ is a {\em{$B$-point}} for $\varphi$ if there exists a sequence $\{\lambda^{(n)}\}$ in $\mathbb{D}^{d}$ converging to $\tau$ such that 
$$\sup_{n}\dfrac{1-\lvert{\varphi(\lambda^{(n)})}\rvert}{1-\|\lambda^{(n)}\|_{\infty}}<\infty.$$
\end{definition}
\begin{example} {\normalfont{\cite[Page 602]{AMYcara} ({\bf $B$-point for a function $\varphi\in\mathcal{S}(\mathbb{D}^{2})$})}}\label{bpointforfuncexample} \rm
 Consider the rational inner function 
       $$\varphi(\lambda)=\dfrac{\frac{1}{2}\lambda_{1}+\frac{1}{2}\lambda_{2}-\lambda_{1}\lambda_{2}}{1-\frac{1}{2}\lambda_{1}-\frac{1}{2}\lambda_{2}},~\lambda\in\mathbb{D}^{2}.$$
It was shown in \cite[Section 6]{AMYcara}  that $(1,1)$ is a $B$-point for 
$\varphi\in\mathcal{S}(\mathbb{D}^{2})$, $\varphi$ has non-tangential limit $1$ at $(1,1)$ and does not extend continuously to $\mathbb{D}^{2}\cup\{(1,1)\}$. Despite these facts, the directional derivative $D_{\delta}\varphi((1,1))$ exists for every direction $\delta$
pointing into the bidisc at $(1,1)$.
    \end{example}

\begin{definition}{\normalfont{({\bf{$C$-point of $\varphi$}})}}
Let $\varphi\in\mathcal{S}(\mathbb{D}^{d})$, let $\tau\in\partial \mathbb{D}^{d}$ be a carapoint of $\varphi$ and suppose $\varphi$ has a nontangential limit at $\tau$. When $\varphi$ has an angular gradient at $\tau$, we say that $\tau$ is a $C$-point of $\varphi$.
\end{definition}
\begin{example}\normalfont{\cite[Page 588]{AMYcara} (\bf{$C$-point for a function $\varphi\in\mathcal{S}(\mathbb{D}^{2})$})} \rm
Consider the rational function 
        $$\varphi(\lambda)=\dfrac{-4\lambda_{1}\lambda_{2}^{2}+\lambda_{2}^{2}+3\lambda_{1}\lambda_{2}-\lambda_{1}+\lambda_{2}}{\lambda_{2}^{2}-\lambda_{1}\lambda_{2}-\lambda_{1}-3\lambda_{2}+4},~\lambda\in\mathbb{D}^{2}.$$
It was shown in \cite{AMYcara} that $\varphi\in\mathcal{S}(\mathbb{D}^{2})$ and $\varphi$ has a $C$-point at $(1,1)\in\mathbb{T}^{2}$.
    \end{example}

\begin{definition}{\normalfont{\cite[Definition 3.4]{ATY2} ({\bf{$C$-point of the model}})}}
Let $(\mathcal{M},u)$ be a model of a function $\varphi\in \mathcal{SA}_{d}$.
 A point $\tau\in\partial\mathbb{D}^{d}$ is a {\em{$C$-point of the model}} if, for every subset $S$ of $\mathbb{D}^{d}$ that approaches $\tau$ nontangentially, $u$ extends continuously to $S\cup\{\tau\}$, with respect to the norm topology on $\mathcal{M}$. 
\end{definition}

For models on the disc, if $\tau\in\mathbb{T}$ is a $B$-point for $\varphi\in\mathcal{S}(\mathbb{D})$, then for all 
sequences $(\lambda^{(n)}) $ in $\mathbb{D}$ such that 
$\lambda^{(n)}{\overset{nt}\rightarrow}\tau$, $u(\lambda^{(n)})$ tends weakly to a  unique limit in the model space, \cite[Proposition 5.30]{AMYOperatorbook}. For the bidisc and upwards, this is no longer the case, as the set of limits of sequences $u(\lambda^{(n)})$ as $\lambda^{(n)}{\overset{nt}\rightarrow}\tau$ is potentially larger than a single vector, \cite[Page 147]{AMYOperatorbook}. Because of this technicality, we shall require the following notion, similar to that found in \cite[Section 4]{AMYcara}.
\begin{definition}{\normalfont{({\bf{Cluster Set}})}}\label{Clustersetattau}
We define the {\em{cluster set}} of a model $(\mathcal{M},u)$ for a function $\varphi \in \mathcal{SA}_{d}$ at a $B$-point $\tau$ to be the set of weak limits in $\mathcal{M}$ of the weakly convergent sequences $\{u(\lambda^{(n)})\}$ as $\{\lambda^{(n)}\}$ ranges over all sequences in $\mathbb{D}^{d}$ that tend to $\tau$ and satisfy
\begin{equation}\label{boundedconditionforBpointonDd}
\sup_{n} J_{\varphi}(\lambda^{(n)})<\infty.
\end{equation}
We denote the cluster set at $\tau$ by $Y_{\tau}$.
\end{definition}
It follows from Proposition \ref{prop4.2generalized}(1) that the cluster set of a model of a function $\varphi \in \mathcal{SA}_{d}$ at a $B$-point is nonempty. 
\begin{theorem}{\normalfont{\cite[Theorem 9.103]{AMYOperatorbook}} ({\bf{Realization formula for the Schur-Agler class}})}\label{realizationtheoremschuragler}
    A function $\varphi$ belongs to $\mathcal{SA}_{d}$ if and only if there exists  an orthogonally decomposed Hilbert space 
    $\mathcal{M}=\mathcal{M}_{1}\oplus\hdots\oplus\mathcal{M}_{d}$,  a scalar $a\in\mathbb{C}$, vectors $\beta,\gamma\in\mathcal{M}$ and a contraction $D$ on $\mathcal{M}$ such that
    $$V = \begin{bmatrix}
        a & 1\otimes \beta \\
        \gamma\otimes 1 & D 
    \end{bmatrix} : \mathbb{C}\oplus \mathcal{M}\rightarrow \mathbb{C}\oplus \mathcal{M}$$
 is a unitary operator and
\begin{equation}\label{aglerrealizations}
\varphi(\lambda)=a+\langle \lambda_{P}(1_{\mathcal{M}}-D\lambda_{P})^{-1}\gamma,\beta\rangle_{\mathcal{M}}
\end{equation}  
for all $\lambda\in\mathbb{D}^{d}$, where $P$ is the $d$-tuple $P=(P_{1},\hdots,P_{d})$ on $\mathcal{M}$,  $P_{j}:\mathcal{M}\rightarrow\mathcal{M}$ is the orthogonal projection onto $\mathcal{M}_{j}$ for $j=1,\hdots,d$ and $\lambda_{P}$ denotes the operator $\lambda_{1}P_{1}+\hdots+\lambda_{d}P_{d}$ on $\mathcal{M}$. Moreover, for all $\lambda\in\mathbb{D}^{d}$,
\begin{align}
\gamma & =( 1_{\mathcal{M}}-D \lambda_{P} ) u(\lambda) \label{gamma-D-u}\\
\varphi(\lambda) &=a+\langle \lambda_{P}u(\lambda),\beta\rangle_{\mathcal{M}}.
\end{align}
\end{theorem}
\begin{remark} \label{mu=la} \rm Let $(\mathcal{M},u)$ be a model of $\varphi \in \mathcal{SA}_{d}$, see Definition \ref{modelpolydefinto}.
By equation (\ref{modeleqn}), for $\mu=\lambda\in\mathbb{D}^{d}$,
\begin{align}
1-\lvert\varphi(\lambda)\rvert^{2} &= \Big\langle{\big(1_{\mathcal{M}}-\lambda_{P}^{*}\lambda_{P}\big)u({\lambda}),u({\lambda})} \Big\rangle_{\mathcal{M}}\nonumber\\
&= \Big\langle{\big(1_{\mathcal{M}}-\sum_{i=1}^{d}\overline{\lambda_{i}}P_{i}\cdot\sum_{j=1}^{d}\lambda_{j}P_{j}\big)u({\lambda}),u({\lambda})} \Big\rangle_{\mathcal{M}}.\label{whenlambdaequalsmuinpolymodel}
\end{align}
Recall that, by Definition \ref{modelpolydefinto}, each $P_{j}:\mathcal{M}\rightarrow\mathcal{M}$ is the orthogonal projection onto $\mathcal{M}_{j}$ for $j=1,\hdots,d$. We have $P_{i}^{*}=P_{i}$ and $P_{i}P_{j}=0$ for every $i\neq j$ and $P_{j}^{2}=P_{j}$ for all $i,j=1,\hdots,d$. We shall denote $P_{j}u(\lambda)$ by $u_{j}(\lambda)$ for all $j=1,\hdots,d$ and $\lambda\in\mathbb{D}^{d}$. Hence
$$\sum_{i=1}^{d}\overline{\lambda_{i}}P^{*}_{i}\cdot\sum_{j=1}^{d}\lambda_{j}P_{j} = \sum_{j=1}^{d}\lvert \lambda_{j}\rvert^{2}P_{j},$$
and thus, for all $\lambda\in\mathbb{D}^{d}$,
\begin{align}
1-\lvert\varphi(\lambda)\rvert^{2} &= \Big\langle{\big(1_{\mathcal{M}}-\sum_{j=1}^{d}\lvert \lambda_{j}\rvert^{2}P_{j}\big)u({\lambda}),u({\lambda})} \Big\rangle_{\mathcal{M}}\nonumber\\
    &=\sum_{j=1}^{d}\Big\langle (1-\lvert \lambda_{j}\rvert^{2})u_{j}({\lambda}),u_{j}({\lambda}) \Big\rangle_{\mathcal{M}_{j}}\nonumber\\
    &=\sum_{j=1}^{d}(1-\lvert \lambda_{j}\rvert^{2})\|u_{j}(\lambda)\|^{2}_{\mathcal{M}_{j}}.\label{whenlambdametmu}
\end{align}
We shall use equation (\ref{whenlambdametmu}) extensively throughout this paper.
\end{remark}

A well-known lurking isometry argument proceeds from a model $(\M,u)$ of a function $\ph\in \sad$ to a realization of $\ph$ \cite{Ag1}. The identical argument applied to a {\em generalized} model $(\M, u, I)$ produces a {\em generalized} notion of realization.
\begin{theorem} \label{genRealiz}
 If $(\L, u, I)$ is a generalized model of $\vp \in \sad$ then there exist a Hilbert space $\M$ containing $\L$, a scalar $a \in \C$, vectors $\beta, \gamma \in \M$ and a linear operator $D:\M \to \M$ such that the operator
\beq \label{defL}
L =  \begin{bmatrix}
  a & 1 \otimes \beta \\
  \gamma \otimes 1 & D
 \end{bmatrix}
\eeq
is unitary on $\C \oplus \M$ and, for all $\lambda \in \D^d$, 
\beq\label{Leq}
 L \vectwo{1}{I(\lambda)u(\lambda)} = \vectwo{\vp(\lambda)}{u(\lambda)},
\eeq
and consequently, for all $\la\in\D^d$,
\beq\label{formphi}
\ph(\la) = a+ \ip{I(\la)(1_{\M}-DI(\la))^{-1}\ga}{\beta}_{\M}.
\eeq
\end{theorem}
\begin{proof}
By equation \eqref{modeleqn-I}, for all $\la,\mu\in\D^d$,
\[
1+\ip{I(\la)u(\la)}{I(\mu)u(\mu)}_{\L} = \overline{\ph(\mu)}\ph(\la) + \ip{u(\la)}{u(\mu)}_{\L}.
\]
We may interpret this equation as an equality between the gramians of two families of vectors in $\C\oplus\L$.
Accordingly we may define an isometric operator
\[
L_0: \spa \left\{ \vectwo{1}{I(\lambda)u(\lambda)} : \la\in\D^d\right\}\to \spa\left\{ \vectwo{\vp(\lambda)}{u(\lambda)}: \la\in\D^d \right\}
\]
by the equation 
\beq
 L_0 \vectwo{1}{I(\lambda)u(\lambda)} = \vectwo{\vp(\lambda)}{u(\lambda)}.
\eeq
 If necessary we may enlarge $\C\oplus\L$ to a space $\C\oplus\M$ in which the domain and range of $L_0$ have equal codimension, and then we may extend $L_0$ to a unitary operator $L$ on $\C\oplus\M$.
Let us represent $L$ as a 2 x 2 block matrix, 
\[
L =  \begin{bmatrix}
  a & 1 \otimes \beta \\
  \gamma \otimes 1 & D
 \end{bmatrix}: \C \oplus \M \to \C \oplus \M,
\]
then
\beq \label{main-connections}
\vectwo{\vp(\lambda)}{u(\lambda)} = L \vectwo{1}{I(\lambda)u(\lambda)} = \begin{bmatrix}
  a & 1 \otimes \beta \\
  \gamma \otimes 1 & D
 \end{bmatrix}\vectwo{1}{I(\lambda)u(\lambda)} \ \text{for all} \ \la\in\D^d.
\eeq
For each $\lambda\in\mathbb{D}^{d}$,  we define $I(\lambda)$ to be equal to $0$ on $\M \ominus \L$, so that $I(\lam)$ becomes an element of $\mathcal{B}(\M)$.
In particular, by the equality \eqref{main-connections}, for all $\lambda\in\mathbb{D}^{d}$,
\begin{align}
    \varphi(\lambda)&=a+\langle I(\lambda) u(\lambda),\beta\rangle_{\mathcal{M}},\\
    \gamma &=(1_{\mathcal{M}}-D I(\lambda)) u(\lambda). 
\end{align}
As $L$ is unitary, $D$ is a contraction and, by assumption (3) of Definition \ref{genmodel}, $\|I(\lambda)\|_{\mathcal{B(\M)}}<1$ for all $\lambda\in\mathbb{D}^{d}$, and so $(1_{\M} -DI(\la))$ is invertible for all $\lambda\in\mathbb{D}^{d}$.
Consequently, we may solve the last equation for $u(\lambda)$ to obtain the equation 
\[
u(\lambda)= (1_{\M} -DI(\la))^{-1}\gamma \ \text{for all} \ \la\in\D^d.
\]
Therefore, for all $\la\in\D^d$,
\beq\label{formphi-2}
\ph(\la) = a+ \ip{I(\la)(1_{\M}-DI(\la))^{-1}\ga}{\beta}_{\M}.
\eeq
\end{proof}

\section{Boundary results for Schur-Agler functions}\label{boundary-results}
\begin{proposition}\label{prop4.2generalized}
Let $\tau\in\partial\mathbb{D}^{d}$ be a $B$-point for $\varphi\in\mathcal{SA}_{d}$, let $(\mathcal{M},u)$ be a model of $\varphi$ and let $Y_{\tau}$ be the cluster set of the  model $(\mathcal{M},u)$ at $\tau$.
\begin{enumerate}
\item If
   $$\liminf_{\substack {\lambda{\overset{nt}\rightarrow}\tau \\ \lambda \in \mathbb{D}^{d}}} \dfrac{1-\lvert \varphi(\lambda) \rvert}{1-\|\lambda\|_{\infty}}=\alpha,$$
then there exists an $x\in Y_{\tau}$ such that $\| x\|^{2}_{\mathcal{M}}\leq \alpha$.
\item If $x\in Y_{\tau}$ and $\lvert\tau_{j}\rvert< 1$ for any $j=1,\hdots,d$, then $x_{j}=0$.
\item There exists $\omega\in\mathbb{T}$ such that, for all $x\in Y_{\tau}$ and $\lambda\in\mathbb{D}^{d}$,
    \begin{equation}\label{omegaequationintheoremstatement}
   1-\overline{\omega}\varphi(\lambda)=\sum_{\lvert \tau_{j} \rvert = 1}\Big(1-\overline{\tau}_{j}\lambda_{j}\Big)\langle u_{j}(\lambda),x_{j}\rangle_{\mathcal{M}_{j}}.
    \end{equation}
\end{enumerate} 
\end{proposition}
\begin{proof} (1).~ Pick a sequence $\{\lambda^{(n)}\}$ in $\mathbb{D}^{d}$ converging to $\tau$ such that 
\begin{equation}\label{alphaequaltothisequation}
\lim_{n\rightarrow \infty} \dfrac{1-\lvert \varphi(\lambda^{(n)})\rvert}{1-\|\lambda^{(n)}\|_{\infty}}=\alpha.
\end{equation}
Since $\tau$ is a $B$-point of $\varphi$, $\alpha<\infty$. Passing to a subsequence if necessary, we can suppose that $\varphi(\lambda^{(n)})$ converges to a point $\omega\in\overline{\mathbb{D}}$. Equation (\ref{alphaequaltothisequation}) implies that $\omega\in\mathbb{T}$. Since $1-\lvert \lambda^{(n)}_{j} \rvert^{2} \geq 1-\|\lambda^{(n)}\|^{2}_{\infty}$, for all $j=1,2, \dots$,
\begin{align*}
\dfrac{1-\lvert \varphi(\lambda^{(n)})\rvert^{2}}{1-\|\lambda^{(n)}\|_{\infty}^{2}}&=\dfrac{1}{1-\|\lambda^{(n)}\|_{\infty}^{2}}\Big\langle (1_{\mathcal{M}}-\lambda_{P}^{*}\lambda_{P})u(\lambda^{(n)}),u(\lambda^{(n)})\Big\rangle_{\mathcal{M}}\nonumber\\
&=\dfrac{1}{1-\|\lambda^{(n)}\|_{\infty}^{2}}\sum_{j=1}^{d}\big(1-\lvert \lambda^{(n)}_{j}\rvert^{2}\big) \| u_{j}(\lambda^{(n)})\|_{\mathcal{M}}^{2}\\
&\geq \sum_{j=1}^{d}\| u_{j}(\lambda^{(n)})\|_{\mathcal{M}}^{2}
= \|u(\lambda^{(n)})\|_{\mathcal{M}}^{2}.\nonumber
\end{align*}
Thus, $\|u(\lambda^{(n)})\|_{\mathcal{M}}$ is bounded. By the compactness and metrizability of closed balls in the separable complex Hilbert space $\mathcal{M}$ in the weak topology, we can take a further subsequence to arrange that $u(\lambda^{(n)})$ converges weakly to some $x\in\mathcal{M}$, where $\|x\|_{\mathcal{M}}^{2}\leq \alpha$. Clearly $x\in Y_{\tau}$.\\

(2) and (3). ~Consider $x\in Y_{\tau}$. Pick a sequence $\{\lambda^{(n)}\}$ in $\mathbb{D}^{d}$ converging to $\tau$ such that $u(\lambda^{(n)})\rightarrow x$ weakly and relation (\ref{boundedconditionforBpointonDd}) holds. Passing to a subsequence if necessary, we can guarantee that $\varphi(\lambda^{(n)})\rightarrow \omega$ for some $\omega\in\overline{\mathbb{D}}$. Since $x\in Y_{\tau}$, the relation (\ref{boundedconditionforBpointonDd}) holds, and we have $\lvert\omega\rvert=1$. For the model $(\mathcal{M},u)$ of $\varphi$, let $\mu=\lambda^{(n)}$ as $n\rightarrow \infty$ in equation (\ref{modeleqn}) to obtain
\begin{equation}\label{omegaisnowinhere}
1-\overline{\omega}\varphi(\lambda)=\sum_{j=1}^{d}(1-\overline{\tau}_{j}\lambda_{j})\langle u_{j}(\lambda),x_{j}\rangle_{\mathcal{M}_{j}}.
\end{equation}
If $\tau\notin\mathbb{T}^{d}$, say $\lvert \tau_{j}\rvert <1$ for some $j=1,\hdots,d$, then on setting $\lambda=\lambda^{(n)}$ in equation (\ref{omegaisnowinhere}) and letting $n\rightarrow\infty$ we find that
$$1-\lvert\omega\rvert^{2}=0=\sum_{j=1}^{d}(1-\lvert\tau_{j}\rvert^{2})\| x_{j}\|^{2}_{\mathcal{M}_{j}}.$$
Therefore $x_{j}=0$ whenever $\tau_{j}\in\mathbb{D}$. This proves statement (2) and equation (\ref{omegaequationintheoremstatement}) holds for all $\lambda\in\mathbb{D}^{d}$. Suppose also that $y\in Y_{\tau}$ is such that $u(\mu_{n})\rightarrow y$ where $\mu_{n}\rightarrow \tau$ and $\varphi(\mu_{n})\rightarrow\zeta$. Put $\lambda=\mu_{n}$ in equation (\ref{omegaequationintheoremstatement}) and let $n\rightarrow \infty$ to obtain $1-\overline{\omega}\zeta=0$, that is, $\omega=\zeta$. Thus $\omega$ is uniquely determined and is independent of the choice of $x\in Y_{\tau}$, and so (3) holds.
\end{proof}

The following theorem 
will be essential to prove the equivalence of various definitions of a $B$-point $\tau$ for $\varphi\in\mathcal{SA}_{d}$.
\begin{theorem}\label{theorem4.9generalized}
Let $\varphi\in\mathcal{SA}_{d}$, $\tau\in\partial\mathbb{D}^{d}$. If there is a sequence $\{\lambda^{(n)}\}$ in $\mathbb{D}^{d}$ converging to $\tau$ such that
$$\lim_{n\rightarrow \infty}\dfrac{1-\lvert \varphi(\lambda^{(n)})\rvert}{1-\|\lambda^{(n)}\|_{\infty}}=\alpha<\infty,$$
then there exists $\omega\in\mathbb{T}$ such that, for all $\lambda\in\mathbb{D}^{d}$, 
\begin{equation}\label{alphamaxequtionforgeneralized}
\dfrac{\lvert \varphi(\lambda)-\omega\rvert^{2}}{1-\lvert \varphi(\lambda) \rvert^{2}} \leq \alpha \max_{\lvert \tau_{j} \rvert=1} \dfrac{\lvert \lambda_{j}-\tau_{j}\rvert^{2}}{1-\lvert \lambda_{j} \rvert^{2}}.
\end{equation}
\end{theorem}
\begin{proof}
We will prove the inequality for $\tau\in\mathbb{T}^{d}$, as by the considerations of Proposition \ref{prop4.2generalized}(2), the conclusion then follows for $\tau\in\partial\mathbb{D}^{d}$.

Pick a model $(\mathcal{M},u)$ of $\varphi$. By Proposition \ref{prop4.2generalized}, there exist an $\omega\in\mathbb{T}$ and an $x\in Y_{\tau}$ such that 
\begin{equation}\label{normxlessalpha}
\|x\|^{2}_{\mathcal{M}}\leq \alpha,
\end{equation}
and equation (\ref{omegaequationintheoremstatement}) holds for all $\lambda\in\mathbb{D}^{d}$. If $\alpha=0$ then $x=0$ and so, by equation (\ref{omegaequationintheoremstatement}), 
$$1-\overline{\omega}\varphi(\lambda)=0 \; \text{which implies that } \; \varphi(\lambda)=\omega.$$
Therefore $\varphi(\lambda)$ is the constant function $\omega$ and so inequality (\ref{alphamaxequtionforgeneralized}) holds when $\alpha=0$.
Consider a fixed $\lambda\in\mathbb{D}^{d}$. From equation (\ref{omegaequationintheoremstatement})
 we have 
\begin{align*}
\lvert \varphi(\lambda)-\omega\rvert^{2}&=\lvert 1-\overline{\omega}\varphi(\lambda)\rvert^{2}\\
&= \Bigg\lvert \sum_{j=1}^{d} \big\langle (1-\overline{\tau}_{j}\lambda_{j})u_{j}(\lambda),x_{j}\big\rangle_{\mathcal{M}_{j}}\Bigg\rvert^{2} \\
&=\Bigg\lvert \Bigg\langle\begin{bmatrix}
	(1-\overline{\tau}_{1}\lambda_{1})u_{1}(\lambda) \\
    \vdots \\
    (1-\overline{\tau}_{d}\lambda_{d})u_{d}(\lambda) \\
\end{bmatrix},\begin{bmatrix}
x_{1}\\
\vdots\\
x_{d}
\end{bmatrix}\Bigg\rangle_{\mathcal{M}}\Bigg\rvert^{2}\\
&\leq\Bigg\|\begin{bmatrix}
	(1-\overline{\tau}_{1}\lambda_{1})u_{1}(\lambda) \\
    \vdots \\
    (1-\overline{\tau}_{d}\lambda_{d})u_{d}(\lambda) \\
\end{bmatrix}\Bigg\|^{2}_{\mathcal{M}}\Bigg\|\begin{bmatrix}
x_{1}\\
\vdots\\
x_{d}
\end{bmatrix}\Bigg\|^{2}_{\mathcal{M}}\\
 &=\bigg(\sum_{j=1}^{d}\|(1-\overline{\tau}_{j}\lambda_{j})u_{j}(\lambda)\|_{\mathcal{M}_{j}}^{2}\bigg)\bigg(\sum_{j=1}^{d}\|x_{j}\|_{\mathcal{M}_{j}}^{2}\bigg).
 \end{align*}
Thus
 \begin{align}\label{phi-omega}
 \lvert \varphi(\lambda)-\omega\rvert^{2}& \leq \|x\|_{\mathcal{M}}^{2}\sum_{j=1}^{d}\lvert 1-\overline{\tau}_{j}\lambda_{j}\rvert^{2}\|u_{j}(\lambda)\|_{\mathcal{M}_{j}}^{2}\nonumber\\
&= \|x\|_{\mathcal{M}}^{2}\sum_{j=1}^{d}\lvert \tau_{j}-\lambda_{j}\rvert^{2}\dfrac{1-\lvert \lambda_{j}\rvert^{2}}{1-\lvert \lambda_{j}\rvert^{2}}\|u_{j}(\lambda)\|_{\mathcal{M}_{j}}^{2}\nonumber\\
&\leq \|x\|_{\mathcal{M}}^{2}\bigg(\max_{j=1,\hdots,d}\dfrac{\lvert \tau_{j}-\lambda_{j}\rvert^{2}}{1-\lvert \lambda_{j}\rvert^{2}}\bigg)\sum_{j=1}^{d}\Big({1-\lvert \lambda_{j}\rvert^{2}}\Big)\|u_{j}(\lambda)\|_{\mathcal{M}_{j}}^{2}\nonumber\\
&=\|x\|_{\mathcal{M}}^{2}\bigg(\max_{j=1,\hdots,d}\dfrac{\lvert \tau_{j}-\lambda_{j}\rvert^{2}}{1-\lvert \lambda_{j}\rvert^{2}}\bigg)(1-\lvert \varphi(\lambda) \rvert^{2}),
\end{align}
where the last line in the above equation is derived from equation (\ref{whenlambdametmu}). Hence by dividing through by $1-\lvert\varphi(\lambda)\rvert^{2}$ and using the inequality (\ref{normxlessalpha}) in the last equation above, we prove inequality (\ref{alphamaxequtionforgeneralized}).\end{proof}

\begin{proposition}\label{prop5.4}
    Let $\varphi\in\mathcal{SA}_{d}$ and let $(\mathcal{M},u)$ be a model of $\varphi$. If $\tau\in\partial\mathbb{D}^{d}$ is a $B$-point for $\varphi$, then $u(\lambda)$ is bounded on any set in $\mathbb{D}^{d}$ that approaches $\tau$ nontangentially.
Moreover, if $S \subset \mathbb{D}^{d}$ approaches $\tau$ nontangentially,   $c>0$ is a constant such that 
$$\|\lambda-\tau\|_{\infty} \leq c(1-\|\lambda\|_{\infty})$$
for all $\lambda\in S$, and $\alpha$ is defined by 
    $$\liminf_{\substack{\lambda\rightarrow \tau \\ \lambda\in\mathbb{D}^{d}}} \dfrac{1-\lvert \varphi(\lambda) \rvert}{1-\|\lambda\|_{\infty}}=\alpha,$$
then
$$\|u(\lambda)\|\leq 2c\sqrt{\alpha}~~\text{for all}~\lambda\in S.$$
\end{proposition}

\begin{proof}
For simplicity, we consider the case when $\tau\in\mathbb{T}^{a}\times \mathbb{D}^{b}$ for some positive integers $a,b$ such that $a+b=d$. The proof extends to all other permutations of the co-ordinates of $\tau\in\partial\mathbb{D}^{d}$. 
By Proposition \ref{prop4.2generalized}, there is an $x\in Y_{\tau}$ such that $\|x\|_{\mathcal{M}}^{2}\leq \alpha$ and
there exists $\omega\in\mathbb{T}$ such that  the equation
\begin{equation}\label{equationforxjequalszerowhentaunotoncircle}
1-\overline{\omega}\varphi(\lambda)=\sum_{\lvert \tau_{j} \rvert = 1}\Big(1-\overline{\tau}_{j}\lambda_{j}\Big)\langle u_{j}(\lambda),x_{j}\rangle_{\mathcal{M}_{j}}
\end{equation}
holds for all $\lambda\in\mathbb{D}^{d}$.  Note in the equation above, by Proposition \ref{prop4.2generalized}, for $\tau\in\partial\mathbb{D}^{d}$ such that $\lvert \tau_{j}\rvert<1$ for some $j=1,\hdots,d$, we have $x_{j}=0$.
Fix a set $S$ that approaches $\tau$ nontangentially with the constant $c>0$. For $\lambda\in S$, 
$$(1-\|\lambda\|_{\infty}^{2})\|u(\lambda)\|_{\mathcal{M}}^{2}=\sum_{j=1}^{d}(1-\|\lambda\|_{\infty}^{2})\|u_{j}(\lambda)\|_{\mathcal{M}_{j}}^{2}
\leq \sum_{j=1}^{d}(1-\lvert\lambda_{j}\rvert^{2})\|u_{j}(\lambda)\|_{\mathcal{M}_{j}}^{2}.$$
By equation (\ref{whenlambdametmu}),
$$(1-\|\lambda\|_{\infty}^{2})\|u(\lambda)\|_{\mathcal{M}}^{2} \leq 1-\lvert \varphi(\lambda)\rvert^{2}
= 1-\overline{\omega}\varphi(\lambda)+\overline{(\omega-\varphi(\lambda))}\varphi(\lambda).$$
By equation (\ref{equationforxjequalszerowhentaunotoncircle}),
\begin{align*}
(1-\|\lambda\|_{\infty}^{2})\|u(\lambda)\|_{\mathcal{M}}^{2}&\leq 2\lvert 1-\overline{\omega}\varphi(\lambda)\rvert \\
&=2\bigg\lvert\sum_{ \lvert\tau_{j}\rvert=1} \Big( 1-\overline{\tau}_{j}\lambda_{j}\Big)\langle u_{j}(\lambda),x_{j}\rangle_{\mathcal{M}_{j}}\bigg\rvert\\
&=2\Bigg\lvert\Bigg\langle\begin{bmatrix}
	(1-\overline{\tau_{1}}\lambda_{1})u_{1}(\lambda) \\
    \vdots \\
    (1-\overline{\tau_{a}}\lambda_{a})u_{a}(\lambda)  \\
\end{bmatrix},\begin{bmatrix}
	x_{1} \\
    \vdots \\
    x_{a}  \\
\end{bmatrix}\Bigg\rangle_{\mathcal{M}_{1}\oplus\hdots\oplus\mathcal{M}_{a}}\Bigg\rvert\\
&\leq 2\Bigg\|\begin{bmatrix}
	(1-\overline{\tau_{1}}\lambda_{1})u_{1}(\lambda) \\
    \vdots \\
    (1-\overline{\tau_{a}}\lambda_{a})u_{a}(\lambda)  \\
\end{bmatrix}\Bigg\|_{_{\mathcal{M}_{1}\oplus\hdots\oplus\mathcal{M}_{a}}}\Bigg\|\begin{bmatrix}
	x_{1} \\
    \vdots \\
    x_{a} \\
\end{bmatrix}\Bigg\|_{_{\mathcal{M}_{1}\oplus\hdots\oplus\mathcal{M}_{a}}}\\
&\leq 2\|\tau-\lambda\|_{\infty}\|u(\lambda)\|_{\mathcal{M}}\|x\|_{\mathcal{M}}.
\end{align*}
By equation (\ref{nontangentialdefpoly}),
$$(1-\|\lambda\|_{\infty}^{2})\|u(\lambda)\|_{\mathcal{M}}^{2}\leq 2c(1-\|\lambda\|_{\infty})\|u(\lambda)\|_{\mathcal{M}}\sqrt{\alpha}.$$
For $\lambda\in S$, $\|\lambda\|_{\infty} < 1$ and  $ 0< 1-\|\lambda\|_{\infty} < 1-\|\lambda\|_{\infty}^{2}$, and so the above calculations show that
$$(1-\|\lambda\|_{\infty}^{2})\|u(\lambda)\|_{\mathcal{M}}^{2}\leq 2c(1-\|\lambda\|_{\infty}^{2})\|u(\lambda)\|_{\mathcal{M}}\sqrt{\alpha}.$$
Hence for $\lambda\in S$, $\|u(\lambda)\|\leq 2c\sqrt{\alpha}$.\end{proof}

\begin{proposition}\label{phi-nt-omega}
    If $\varphi\in\mathcal{SA}_{d}$ and $\tau\in\partial\mathbb{D}^{d}$ is a $B$-point for $\varphi$ then there exists $\omega\in\mathbb{T}$ such that 
 $$ \lim_{\lambda{\overset{nt}\rightarrow} \tau} \varphi(\lambda)= \omega.$$  
\end{proposition}
\begin{proof}
Let $(\mathcal{M},u)$ be a model of $\varphi$, and let $\alpha$ be defined by 
    $$\liminf_{\substack{\lambda\rightarrow \tau \\ \lambda\in\mathbb{D}^{d}}} \dfrac{1-\lvert \varphi(\lambda) \rvert}{1-\|\lambda\|_{\infty}}=\alpha.$$
By Proposition \ref{prop4.2generalized}, there is an $x\in Y_{\tau}$ such that $\|x\|_{\mathcal{M}}^{2}\leq \alpha$ and
there exists $\omega\in\mathbb{T}$ such that the equation
\begin{equation}\label{omega-tau}
1-\overline{\omega}\varphi(\lambda)=\sum_{\lvert \tau_{j} \rvert = 1}\Big(1-\overline{\tau}_{j}\lambda_{j}\Big)\langle u_{j}(\lambda),x_{j}\rangle_{\mathcal{M}_{j}}
\end{equation}
holds for all $\lambda\in\mathbb{D}^{d}$.  Note that in the equation above, by Proposition \ref{prop4.2generalized}, for $\tau\in\partial\mathbb{D}^{d}$ such that $\lvert \tau_{j}\rvert<1$ for some $j=1,\hdots,d$, we have $x_{j}=0$.
Fix a set $S$ that approaches $\tau$ nontangentially with the constant $c>0$. By equation \eqref{omega-tau},
for $\lambda\in S$, 
\begin{align}
|\varphi(\lambda) - \omega|& =| 1-\overline{\omega}\varphi(\lambda)| \leq \sum_{\lvert \tau_{j} \rvert = 1}
\lvert \Big(1-\overline{\tau}_{j}\lambda_{j}\Big)\langle u_{j}(\lambda),x_{j}\rangle_{\mathcal{M}_{j}}\rvert\\
  &\leq \sum_{\lvert \tau_{j} \rvert = 1} \lvert 1-\overline{\tau}_{j}\lambda_{j} \rvert \|u_{j}(\lambda)\|_{\mathcal{M}_{j}}
 \|x_{j}\|_{\mathcal{M}_{j}}.
\end{align}
By Proposition \ref{prop5.4},
$$\|u(\lambda)\|\leq 2c\sqrt{\alpha}~~\text{for all}~\lambda\in S.$$
Therefore, for all $\lambda\in S$,
\begin{equation}
|\varphi(\lambda) - \omega| \leq 2c \alpha \sum_{\lvert \tau_{j} \rvert = 1} \lvert 1-\overline{\tau}_{j}\lambda_{j} \rvert.
\end{equation}
Thus, 
$$ \lim_{\lambda{\overset{nt}\rightarrow\tau, \lambda \in S} } \varphi(\lambda)= \omega.$$  
\end{proof}

We will need the language of horocycles and horospheres to formulate the following two corollaries to Theorem \ref{theorem4.9generalized}.
\begin{definition}{\normalfont{\cite{AMYcara}}}\label{horocycle}
A {\em{horocycle}} in $\mathbb{D}$ is a set of the form $E(\tau,R)$ for some $\tau\in \overline{\mathbb{D}}$ and $R>0$, where 
$$E(\tau,R)=\bigg\{\lambda\in\mathbb{D}:\dfrac{\lvert \lambda-\tau\rvert^{2}}{1-\lvert\lambda\rvert^{2}}<R\bigg\}
\; \text{ if } \; \tau \in \mathbb{T}, \;\; \text{and} \; \;  E(\tau,R)=\mathbb{D} \; \;
\text{if} \;\; \tau \in \mathbb{D}.$$
For $\tau\in \partial\mathbb{D}^{d}$ and $R>0$, the $d$-dimensional {\em{horosphere}} is defined to be the set 
$$E(\tau,R)=E(\tau_{1},R)\times \hdots \times E(\tau_{d},R).$$
\end{definition}
Let $D(z,r)$ denote the Euclidean disc in $\mathbb{C}$ with centre $z$ and radius $r>0$. Then
\begin{equation}\label{horocycleequaltoadisc}
E(\tau,R)=D\bigg(\dfrac{\tau}{R+1},\dfrac{R}{R+1}\bigg).
\end{equation}
Note that $E(\tau_{j},R)$ is the circular disc internally tangent to $\mathbb{T}$ at $\tau_{j}$ having radius $R/(R+1)$, \cite[Section 4]{AMYcara}.

\begin{corollary}\label{horocyclecorollary}
Let $\varphi\in\mathcal{SA}_{d}$, $\tau\in\partial\mathbb{D}^{d}$. If there is a sequence $\{\lambda^{(n)}\}$ in $\mathbb{D}^{d}$ converging to $\tau$ such that
$$\lim_{n\rightarrow \infty}\dfrac{1-\lvert \varphi(\lambda^{(n)})\rvert}{1-\|\lambda^{(n)}\|_{\infty}}=\alpha<\infty,$$
then there exists $\omega\in\mathbb{T}$ such that, for any $R>0$, 
$$\varphi(E(\tau,R))\subset E(\omega,\alpha R).$$
\end{corollary}
\begin{proof}
    Let $\lambda\in E(\tau_{1},R)\times\hdots\times E(\tau_{d},R)$, that is, 
$$\dfrac{\lvert \lambda_{j}-\tau_{j}\rvert^{2}}{1-\lvert\lambda_{j}\rvert^{2}}<R~~\text{for all}~j=1,\hdots,d.$$
By Theorem \ref{theorem4.9generalized}, 
$$\dfrac{\lvert \varphi(\lambda)-\omega\rvert^{2}}{1-\lvert \varphi(\lambda)\rvert^{2}}\leq \alpha R.$$
Therefore $\varphi(\lambda)\in E(\omega,\alpha R)$.\end{proof}

The following statement allows us to test whether $\tau$ is a $B$-point using only the values of $\varphi$ along the radius through $\tau$. 
\begin{corollary}\label{coro4.14cara}
Let $\varphi\in\mathcal{SA}_{d}$ and $\tau\in\partial\mathbb{D}^{d}$. The following conditions are equivalent:
\begin{enumerate}
    \item $\tau$ is a $B$-point for $\varphi$;
    \item $\dfrac{1-\lvert \varphi(\lambda) \rvert}{1-\|\lambda\|_{\infty}}$ is bounded on the radius $\{r\tau : 0<r<1\}$;
    \item $\dfrac{1-\lvert \varphi(\lambda) \rvert}{1-\|\lambda\|_{\infty}}$ is bounded on every set in $\mathbb{D}^{d}$ that approaches $\tau$ non-tangentially.
\end{enumerate}
Moreover, if $\tau$ is a $B$-point for $\varphi$, 
\begin{equation}\label{equationforcorollary4.14}
\lim_{r\rightarrow 1^{-}} \dfrac{1-\lvert \varphi(r\tau)\rvert}{1-r} = \liminf_{\substack{\lambda\rightarrow \tau \\ \lambda\in\mathbb{D}^{d}}}\dfrac{1-\lvert \varphi(\lambda)\rvert}{1-\|\lambda\|_{\infty}}.
\end{equation}
\end{corollary}
\begin{proof}
It is easily seen that (3) $\implies$ (2) $\implies$(1).

(1) $\implies$ (3). Suppose $\tau$ is a $B$-point of $\varphi$ and consider a set $S\subset \mathbb{D}^{d}$ that approaches $\tau$ in such a way that $ \|\lambda - \tau \|_{\infty}\leq c ( 1 - \|\lambda\|_{\infty})$, for some $c\geq 1$. Let 
$$\alpha = \liminf_{\substack{\lambda\rightarrow \tau \\ \lambda\in\mathbb{D}^{d}}} \dfrac{1-\lvert \varphi(\lambda)\rvert}{1-\|\lambda\|_{\infty}},$$
then by the first assumption, $\alpha$ is finite. By Corollary \ref{horocyclecorollary}, there exists an $\omega\in\mathbb{T}$ such that, for all $R>0$, 
\begin{equation}\label{varphiequationinthisproof}
\varphi(E(\tau,R))\subset E(\omega,\alpha R) = D\bigg(\dfrac{\omega}{\alpha R+1}, \dfrac{\alpha R}{\alpha R+1}\bigg).
\end{equation}
Pick any $\lambda\in S$ and $\epsilon>0$ and let 
$$R_1=(1+\epsilon)\dfrac{\|\lambda-\tau\|_{\infty}^{2}}{1-\|\lambda\|_{\infty}^{2}}.$$
Then $\lambda\in E(\tau, R_1)$ and, since $\lambda \in S$, 
\begin{equation}\label{stemR1}
0<R_1 \leq (1+\epsilon)\dfrac{c^{2}(1-\|\lambda\|_{\infty})^{2}}{1-\|\lambda\|_{\infty}^{2}}\leq (1+\epsilon)c^{2}(1-\|\lambda\|_{\infty}).
\end{equation}
By the relation (\ref{varphiequationinthisproof}),
\begin{align*}
1-\lvert \varphi(\lambda)\rvert & \leq \lvert \varphi(\lambda)-\omega\rvert \\
& \leq  \lvert \varphi(\lambda)-\dfrac{\omega}{\alpha R_1+1} \rvert + 
\lvert \dfrac{\omega}{\alpha R_1+1} - \omega \rvert\\
                                   & \leq \dfrac{2\alpha R_1}{\alpha R_1+1}.
                                    \end{align*}                            
By inequality \eqref{stemR1},
$$1-\lvert \varphi(\lambda)\rvert \leq \dfrac{2\alpha R_1}{\alpha R_1+1} \leq 2\alpha R_1 
                                    \leq 2(1+\epsilon)\alpha c^{2}(1-\|\lambda\|_{\infty}).$$
Since the inequality above holds for all $\epsilon>0$, by dividing through by $(1-\|\lambda\|_{\infty})$, we have
$$\dfrac{1-\lvert \varphi(\lambda) \rvert}{1-\|\lambda\|_{\infty}} \leq 2\alpha c^{2}.$$
Hence (1)$\implies$(3).

 We now prove equation (\ref{equationforcorollary4.14}). By Theorem \ref{theorem4.9generalized}, there exists $\omega\in\mathbb{T}$ such that the inequality (\ref{alphamaxequtionforgeneralized}) holds for all $\lambda\in\mathbb{D}^{d}$. In particular, choosing $\lambda=r\tau$, we have
$$\dfrac{\lvert \varphi(r\tau)-\omega\rvert^{2}}{1-\lvert \varphi(r\tau)\rvert^{2}}\leq\alpha \dfrac{1-r}{1+r},$$
for $0\leq r<1$. This means that $\varphi(r\tau)$ is in the horocycle $E(\omega,R)$ where $R=\alpha(1-r)/(1+r)$. Since $1-\lvert\varphi(r\tau)\rvert\leq\lvert \varphi(r\tau)-\omega\rvert\leq2R/(R+1)$, 
\begin{equation}\label{inequalitywith2alpha}
\dfrac{1-\lvert\varphi(r\tau)\rvert}{1-r}\leq \dfrac{\lvert\varphi(r\tau)-\omega\rvert}{1-r}\leq\dfrac{2\alpha}{1+\alpha+r(1-\alpha)}.
\end{equation}
Let 
$$f(r) = \dfrac{2\alpha}{1+\alpha+r(1-\alpha)}.$$
We claim that $f(r)$ is bounded for $0\leq r< 1$. Note that $(1+\alpha)+r(1-\alpha)=0$ if and only if $r=(1+\alpha)/(\alpha -1)$.  If $\alpha>1$, then 
$$\dfrac{1+\alpha}{\alpha-1} = \dfrac{\alpha-1+2}{\alpha-1}=1+\dfrac{2}{\alpha-1}>1.$$
If $0\leq\alpha<1$, then $\alpha-1<0$ and so
$$\dfrac{1+\alpha}{\alpha-1}<0.$$
Hence, for $0\leq r <1$, $f(r)$ is bounded, and inequality (\ref{inequalitywith2alpha}) is properly defined, since the denominator is non-zero. Taking limits of the right hand side of inequality (\ref{inequalitywith2alpha}),
$$\lim_{r\rightarrow 1^{-}}\dfrac{2\alpha}{1+r+\alpha(1-r)}=\alpha.$$
Thus we have
$$\alpha=\liminf_{\substack{\lambda\rightarrow\tau \\ \lambda\in\mathbb{D}^{d}}}\dfrac{1-\lvert \varphi(\lambda)\rvert}{1-\|\lambda\|_{\infty}}\leq \liminf_{r\rightarrow 1^{-}} \dfrac{1-\lvert \varphi(r\tau)\rvert}{1-r}\leq \limsup_{r\rightarrow 1^{-}} \dfrac{1-\lvert\varphi(r\tau)\rvert}{1-r}\leq \alpha.$$
Hence we prove equation (\ref{equationforcorollary4.14}). \end{proof}

\section{Connection between $B$-points for $\varphi$ and realizations}\label{B-point-realization}
Proposition \ref{prop5.2cara} will give us a connection between a $B$-point for a function $\varphi\in\mathcal{SA}_{d}$ and a $B$-point  of a model $(\mathcal{M},u)$ for $\varphi$.
Recall Definition \ref{B-point-phi},  $\tau\in\partial \mathbb{D}^{d}$ is a  $B$-point for $\varphi$
 if there exists a sequence $\{\lambda^{(n)}\}$ in $\mathbb{D}^{d}$ converging to $\tau$ such that 
 \begin{equation} \label{2.3}
\sup_{n}\dfrac{1-\lvert{\varphi(\lambda^{(n)})}\rvert}{1-\|\lambda^{(n)}\|_{\infty}}<\infty.
\end{equation} 

\begin{lemma}\label{limit-mod}
Let $\varphi\in\mathcal{SA}_{d}$, let $ \ \tau \in \partial \mathbb{D}^{d}$ and let $\tau$ be a $B$-point for $\varphi$.
The following statements hold:\\
{\rm (i)}
Suppose $\{\lambda^{(n)}\}$ is a sequence in $\mathbb{D}^d$ such that $\lambda^{(n)} \rightarrow \tau$ and Condition \eqref{2.3} is satisfied. Then $\lim_{n\to\infty} |\ph(\lambda^{(n)})| =1$.\\
{\rm (ii)} For every subset $S$ of $\D^d$ that approaches $\tau$ nontangentially,  $\lim_{\la\to\tau, \la\in S} |\varphi(\lambda)| =1$,
\end{lemma}
\begin{proof}
(i).  Condition \eqref{2.3} implies that, for some $M\geq 0,$
\[
1-|\ph(\lambda^{(n)}| \leq M(1-\|\lambda^{(n)}\|_{\infty}) \quad \text{for all} \quad n\in\mathbb{N}.
\]
Since $\lambda^{(n)} \to \tau \in \partial \mathbb{D}^d$, we have $\|\lambda^{(n)}\|_{\infty} \to 1$, and so $1-|\ph(\lambda^{(n)})| \to 0$ as $n\to\infty$.

(ii). By Corollary \ref{coro4.14cara},
the following conditions are equivalent:
\begin{enumerate}
    \item $\tau$ is a $B$-point for $\varphi$;
    \item $\dfrac{1-\lvert \varphi(\lambda) \rvert}{1-\|\lambda\|_{\infty}}$ is bounded on every set in $\mathbb{D}^{d}$ that approaches $\tau$ non-tangentially.
\end{enumerate}
Therefore, for a subset $S$ of $\D^d$ that approaches $\tau$ nontangentially, there is $M\geq 0,$ such that
$${1-\lvert \varphi(\lambda) \rvert} \leq M ({1-\|\lambda\|_{\infty}})\quad \text{for all} \quad \lambda \in S.$$ 
Since $\lambda \to \tau, \la\in S$, we have $\|\lambda\|_{\infty} \to 1$, and so $1-|\ph(\lambda)| \to 0$ as 
$\la\to\tau, \la\in S$.
\end{proof}

\begin{proposition}\label{prop5.2cara}
Let $\varphi\in\mathcal{SA}_{d}$,  let  $ \ \tau \in \partial \mathbb{D}^d$, and let
 $\{\lambda^{(n)}\}$ be a sequence in $\mathbb{D}^{d}$ such that $\lambda^{(n)}{\overset{nt}\rightarrow}\tau$.\\
{\rm (i)} If $ \ \tau \in \partial \mathbb{D}^d \setminus \mathbb{T}^d$,  then  the following statements are equivalent:
 \begin{enumerate}
    \item $\sup_{n} \dfrac{1-\lvert \varphi(\lambda^{(n)})\rvert}{1-\|\lambda^{(n)}\|_{\infty}}<\infty;$
    \item  $ \lim_{n \to \infty} |\ph(\lambda^{(n)})| =1$ and  there exists a model $(\mathcal{M},u)$ of $\varphi$ such that $(u(\lambda^{(n)}))_{n=1}^{\infty}$ is bounded in $\mathcal{M}$;
    \item  $ \lim_{n \to \infty} |\ph(\lambda^{(n)})| =1$, and  for every model $(\mathcal{M},u)$ of $\varphi$, $(u(\lambda^{(n)}))_{n=1}^{\infty}$ is bounded in $\mathcal{M}$.
\end{enumerate}
 {\rm(ii)} If $ \ \tau \in \T^d$, then  the following statements are equivalent:
\begin{enumerate}
    \item $\sup_{n} \dfrac{1-\lvert \varphi(\lambda^{(n)})\rvert}{1-\|\lambda^{(n)}\|_{\infty}}<\infty;$
    \item there exists a model $(\mathcal{M},u)$ of $\varphi$ such that $(u(\lambda^{(n)}))_{n=1}^{\infty}$ is bounded in $\mathcal{M}$;
    \item for every model $(\mathcal{M},u)$ of $\varphi$, $(u(\lambda^{(n)}))_{n=1}^{\infty}$ is bounded in $\mathcal{M}$.
\end{enumerate}
\end{proposition}
\begin{proof}
Consider  a sequence $\{\lambda^{(n)}\}$ in $\mathbb{D}^{d}$ such that $\lambda^{(n)}{\overset{nt}\rightarrow}\tau$
and choose a $c$ such that 
\begin{equation} \label{c-nontang}
\|\tau-\lambda^{(n)}\|_{\infty}\leq c(1-\|\lambda^{(n)}\|_{\infty}) \;\; \text{ for all} \; \; n. 
\end{equation}
Clearly, in both cases (i) and (ii), (3) implies (2).  

Let us show that in both cases (i) and (ii), (1)$\Rightarrow$(3). Assume (1): say the inequality  (\ref{2.3}) holds with bound $M$:  
\begin{equation}\label{sup=M}
\sup_{n}\dfrac{1-\lvert \varphi(\lambda^{(n)})\rvert}{1-\|\lambda^{(n)}\|_{\infty}}=M<\infty.
\end{equation}
Observe that, for $\lambda^{(n)}\in\mathbb{D}^{d}$, 
$$ \lvert \lambda^{(n)}_{j}\rvert^{2} \leq \lvert \lambda^{(n)}_{j} \rvert \leq \|\lambda^{(n)}\|_{\infty}~~\text{and}~~ 1-\lvert \lambda^{(n)}_{j}\rvert^{2} \geq 1- \lvert \lambda^{(n)}_{j} \rvert \geq 1- \|\lambda^{(n)}\|_{\infty}.$$
Thus, if $(\mathcal{M},u)$ is a model of $\varphi$,
\begin{align}
(1-\|\lambda^{(n)}\|_{\infty})\|u(\lambda^{(n)})\|_{\mathcal{M}}^{2} &= (1-\|\lambda^{(n)}\|_{\infty})\sum_{j=1}^{d}\|u_{j}(\lambda^{(n)})\|_{\mathcal{M}_{j}}^{2} \nonumber\\
&\leq \sum_{j=1}^{d}(1-\lvert \lambda^{(n)}_{j}\rvert^{2})\|u_{j}(\lambda^{(n)})\|_{\mathcal{M}_{j}}^{2} \nonumber\\
&=\bigg\langle(1_{\mathcal{M}}-(\lambda^{(n)})_{P}^{*}(\lambda^{(n)})_{P}u(\lambda^{(n)}),u(\lambda^{(n)})\bigg\rangle_{\mathcal{M}}. \label{4.3.2}
\end{align} 
By equation (\ref{whenlambdaequalsmuinpolymodel}),
\begin{equation} \label{model-lambda}
 1-\lvert \varphi(\lambda^{(n)})\rvert^{2}= \bigg\langle(1_{\mathcal{M}}-(\lambda^{(n)})_{P}^{*}(\lambda^{(n)})_{P}u(\lambda^{(n)}),u(\lambda^{(n)})\bigg\rangle_{\mathcal{M}}.
 \end{equation}
Therefore, by inequality \eqref{4.3.2}, equation \eqref{model-lambda} and the hypothesis \eqref{sup=M}, we obtain
$$(1-\|\lambda^{(n)}\|_{\infty})\|u(\lambda^{(n)})\|_{\mathcal{M}}^{2} \leq 1-\lvert \varphi(\lambda^{(n)})\rvert^{2}
\leq 2M(1-\|\lambda^{(n)}\|_{\infty}).$$ 
Hence, if we divide through by $(1-\|\lambda^{(n)}\|_{\infty})$ in the last inequality above, we deduce that in both Cases,  (1) implies that, for all models $(\mathcal{M},u)$ of  $ \varphi$, $u(\lambda^{(n)})$ is  bounded. By Lemma \ref{limit-mod},   $ \lim_{n \to \infty} |\varphi(\lambda^{(n)})| =1$. Therefore  (1)$\Rightarrow$(3). 

It remains to show that (2)$\Rightarrow$ (1).
In Case (ii), assume (2):  there exists  a  model $(\mathcal{M},u)$ of  $ \varphi$   such  that $u(\lambda^{(n)})$ is  bounded. 
Let $\|u(\lambda^{(n)})\|_{\mathcal{M}}\leq M$ for some model $(\mathcal{M},u)$ of $\varphi$ and all $n\geq 1$ and suppose that $\tau\in\mathbb{T}^{d}$. Clearly, for $j=1,\hdots,d$ and for all $\lambda^{(n)}\in\mathbb{D}^{d}$,
\begin{align}
\lvert \lambda^{(n)}_{j} \rvert < 1 &\iff 1+\lvert \lambda^{(n)}_{j} \rvert < 2 \nonumber\\
    & \iff (1-\lvert \lambda^{(n)}_{j}\rvert )(1+\lvert \lambda^{(n)}_{j} \rvert) < 2(1-\lvert \lambda^{(n)}_{j}\rvert) \nonumber\\
    & \iff 1-\lvert \lambda^{(n)}_{j} \rvert^{2} < 2(1-\lvert \lambda^{(n)}_{j}\rvert)\label{inequalitypart1}.
\end{align}
Note that
\begin{equation}\label{inequalitypart2}
1-\lvert \lambda^{(n)}_{j}\rvert = \lvert \tau_{j}\rvert-\lvert \lambda^{(n)}_{j}\rvert \leq \lvert \tau_{j}-\lambda^{(n)}_{j}\rvert.
\end{equation}
By equation (\ref{whenlambdaequalsmuinpolymodel}),
\begin{align*}
1-\lvert \varphi(\lambda^{(n)})\rvert &\leq 1-\lvert\varphi(\lambda^{(n)})\rvert^{2}\\
    &=\bigg\langle(1_{\mathcal{M}}-(\lambda^{(n)})_{P}^{*}(\lambda^{(n)})_{P})u(\lambda^{(n)}),u(\lambda^{(n)})\bigg\rangle_{\mathcal{M}}\\
     &=\sum_{j=1}^{d}(1-\lvert \lambda^{(n)}_{j}\rvert^{2})\|u_{j}(\lambda^{(n)})\|_{\mathcal{M}_{j}}^{2}.
\end{align*}
By inequality (\ref{inequalitypart1}),
$$1-\lvert \varphi(\lambda^{(n)})\rvert \leq \sum_{j=1}^{d}2(1-\lvert \lambda^{(n)}_{j}\rvert)\|u_{j}(\lambda^{(n)})\|_{\mathcal{M}_{j}}^{2}.
$$
By inequality (\ref{inequalitypart2}),
$$1-\lvert \varphi(\lambda^{(n)})\rvert     \leq 2\sum_{j=1}^{d}\lvert \tau_{j}-\lambda^{(n)}_{j}\rvert\|u_{j}(\lambda^{(n)})\|_{\mathcal{M}_{j}}^{2}
    \leq 2\|\tau-\lambda^{(n)}\|_{\infty}\sum_{j=1}^{d}\|u_{j}(\lambda^{(n)})\|_{\mathcal{M}_{j}}^{2}.$$
By assumption, $\lambda^{(n)}{\overset{nt}\rightarrow}\tau$. Hence by inequality (\ref{c-nontang}),
$$1-\lvert \varphi(\lambda^{(n)})\rvert     \leq 2c(1-\|\lambda^{(n)}\|_{\infty})\sum_{j=1}^{d}\|u_{j}(\lambda^{(n)})\|_{\mathcal{M}_{j}}^{2}\\
    \leq 2cM^{2}(1-\|\lambda^{(n)}\|_{\infty}).$$
Therefore (2)$\implies$(1) when $\tau\in\mathbb{T}^{d}$.

Let us prove that (2)$\implies$(1) in Case (i) when $\tau\in\partial\mathbb{D}^{d}\setminus\mathbb{T}^{d}$. Since any permutation of co-ordinates is an automorphism of $\overline{\mathbb{D}}^{d}$, it is enough to prove that (2)$\implies$(1) under the assumption that $\tau\in\mathbb{T}^{a}\times \mathbb{D}^{b}$ for some positive integers $a,b$ such that $a+b=d$.
Assume that statement (2) holds:
 $ \lim_{n \to \infty} |\varphi(\lambda^{(n)})| =1$ and  there exists  a  model $(\mathcal{M},u)$ of  $ \varphi$   such  that $u(\lambda^{(n)})$ is  bounded.

Note that $ \lim_{n \to \infty} |\varphi(\lambda^{(n)})| =1$, and so there exist a subsequence of $\{\lambda^{(n)}\}$ and $\omega\in\T$  such that 
\[
 \lim_{n\to\infty} \ph(\lambda^{(n)}) =\omega.
\]
For the model $(\mathcal{M},u)$ of $\varphi$, let $\mu=\lambda^{(n)}$, $n= 1,2, \dots$, in equation (\ref{modeleqn}) to obtain, for $\la\in\D^d$, 
\begin{equation} \label{2.8b}
1 - \overline{\varphi(\lambda^{(n)})} \varphi(\lambda) = \sum_{j=1}^{d}(1-\overline{\lambda^{(n)}_{j}}\lambda_{j})
\langle u_{j}(\lambda),u_{j}( \lambda^{(n)})  \rangle_{\mathcal{M}_{j}}.
\end{equation} 
Since  $u(\lambda^{(n)})$ is  bounded in $\M$, we may pass to a subsequence of $\{\lambda^{(n)}\}$ for which there exists $x\in Y_\tau$ such that $u(\lambda^{(n)})$ converges weakly to the point $x=(x_1, x_2, \dots, x_d)$ as $ n \to \infty$.
Take limits in equation \eqref{2.8b} as $\lambda^{(n)} \rightarrow \tau$, to obtain, for $\la\in\D^d$,
\begin{equation}\label{2.8c}
1-\overline{\omega}\varphi(\lambda)=\sum_{j=1}^{d}(1-\overline{\tau}_{j}\lambda_{j})\langle u_{j}(\lambda),x_{j}\rangle_{\mathcal{M}_{j}}.
\end{equation}

If $\tau\in\mathbb{T}^{a}\times \mathbb{D}^{b}$ is  such that  $|\tau_j| < 1$, for some $j=1,\hdots,d$, then on setting $\la=\lambda^{(n)}$ in equation (\ref{2.8c}) and letting $n\to \infty$ we find 
$$1-\lvert\omega\rvert^{2}=0=\sum_{j=1}^{d}(1-\lvert\tau_{j}\rvert^{2})\| x_{j}\|^{2}_{\mathcal{M}_{j}}.$$
Hence $x_{j}=0$ whenever $\tau_{j}\in\mathbb{D}$. 
Therefore, there exists an $\omega\in\mathbb{T}$ and an $x\in\mathcal{M}$ such that, for all $\lambda\in\mathbb{D}^{d}$,
$$1-\overline{\omega}\varphi(\lambda^{(n)})=\sum_{j=1}^{a}(1-\overline{\tau}_{j}\lambda^{(n)}_{j})\langle u_{j}(\lambda),x_{j}\rangle_{\mathcal{M}_{j}},$$
and 
$$\overline{\omega}\varphi(\lambda^{(n)})= 1+\sum_{j=1}^{a}(\overline{\tau}_{j})(\lambda^{(n)}_{j}-\tau_{j})\langle u_{j}(\lambda^{(n)}),x_{j}\rangle_{\mathcal{M}_{j}}.$$
Thus,
\begin{align*}
\lvert \varphi(\lambda^{(n)})\rvert^{2} &=\bigg\lvert 1+\sum_{j=1}^{a}(\overline{\tau}_{j})(\lambda^{(n)}_{j}-\tau_{j})\langle u_{j}(\lambda^{(n)}),x_{j}\rangle_{\mathcal{M}_{j}}\bigg\rvert^{2}\\
&=1+2\mathrm{Re}\bigg(\sum_{j=1}^{a}\overline{\tau}_{j}(\lambda^{(n)}_{j}-\tau_{j})\langle u_{j}(\lambda^{(n)}),x_{j}\rangle_{\mathcal{M}_{j}}\bigg)+\bigg\lvert \sum_{j=1}^{a}\overline{\tau}_{j}(\lambda^{(n)}_{j}-\tau_{j})\langle u_{j}(\lambda^{(n)}),x_{j}\rangle_{\mathcal{M}_{j}} \bigg\rvert^{2}.
\end{align*}
Since $\lvert 1-\lvert\varphi(\lambda^{(n)})\rvert^{2}\rvert  = 1-\lvert\varphi(\lambda^{(n)})\rvert^{2}$, we get
\begin{align*}
    & 1-\lvert \varphi(\lambda^{(n)})\rvert^{2} \\
    &=\Bigg\lvert 2\mathrm{Re}\bigg(\sum_{j=1}^{a}\overline{\tau}_{j}(\lambda^{(n)}_{j}-\tau_{j})\langle u_{j}(\lambda^{(n)}),x_{j}\rangle_{\mathcal{M}_{j}}\bigg)+\bigg\lvert \sum_{j=1}^{a}\overline{\tau}_{j}(\lambda^{(n)}_{j}-\tau_{j})\langle u_{j}(\lambda^{(n)}),x_{j}\rangle_{\mathcal{M}_{j}} \bigg\rvert^{2}\Bigg\rvert.
\end{align*}
Next, divide the above equation by $1-\|\lambda^{(n)}\|_{\infty}^{2}$ to find the following inequality
\begin{align}
    &\dfrac{1-\lvert \varphi(\lambda^{(n)})\rvert^{2}}{1-\|\lambda^{(n)}\|_{\infty}^{2}} \nonumber\\
    &= \dfrac{1}{1-\|\lambda^{(n)}\|_{\infty}^{2}}\Bigg\lvert 2\mathrm{Re}\bigg(\sum_{j=1}^{a}\overline{\tau}_{j}(\lambda^{(n)}_{j}-\tau_{j})\langle u_{j}(\lambda^{(n)}),x_{j}\rangle_{\mathcal{M}_{j}}\bigg)\nonumber\\  
&\hspace{6.5cm}+\bigg\lvert \sum_{j=1}^{a}\overline{\tau}_{j}(\lambda^{(n)}_{j}-\tau_{j})\langle u_{j}(\lambda^{(n)}),x_{j}\rangle_{\mathcal{M}_{j}} \bigg\rvert^{2}\Bigg) \nonumber \\       
    &\leq \dfrac{1}{1-\|\lambda^{(n)}\|_{\infty}}\Bigg(2\bigg\lvert\sum_{j=1}^{a}\overline{\tau}_{j}(\lambda^{(n)}_{j}-\tau_{j})\langle u_{j}(\lambda^{(n)}),x_{j}\rangle_{\mathcal{M}_{j}}\bigg\rvert\nonumber\\
    &\hspace{6.5cm}+\bigg\lvert \sum_{j=1}^{a}\overline{\tau}_{j}(\lambda^{(n)}_{j}-\tau_{j})\langle u_{j}(\lambda^{(n)}),x_{j}\rangle_{\mathcal{M}_{j}} \bigg\rvert^{2}\Bigg).\label{lastinequalityhereg}\\\nonumber
\end{align} 
Define the vector-function $v:\mathbb{D}^{d}\rightarrow\mathcal{M}$ by 
\begin{equation*}
v(\lambda)=\begin{bmatrix}
	\overline{\tau}_{1}(\lambda_{1}-\tau_{1})u_{1}(\lambda) \\
    \vdots \\
    \overline{\tau}_{a}(\lambda_{a}-\tau_{a})u_{a}(\lambda)  \\
\end{bmatrix}~~\text{for all}~~\lambda\in\mathbb{D}^{d},
\end{equation*}
and observe that 
$$    \| v(\lambda)\|_{\mathcal{M}} \leq \| \lambda - \tau \|_\infty  \| u(\lambda)\|_{\mathcal{M}} \; \text{ for all } \; \lambda \in \mathbb{D}^{d}.$$ 
Then the inequality (\ref{lastinequalityhereg}) can be written in the following way:
\begin{align*}
    &\dfrac{1-\lvert \varphi(\lambda^{(n)})\rvert^{2}}{1-\|\lambda^{(n)}\|_{\infty}^{2}}\\
    &\leq\dfrac{1}{1-\|\lambda^{(n)}\|_{\infty}}\Bigg(2\Bigg\lvert\Bigg\langle v(\lambda^{(n)}),\begin{bmatrix}
	x_{1} \\
    \vdots \\
    x_{a}  \\
\end{bmatrix}\Bigg\rangle_{\mathcal{M}}\Bigg\rvert+\Bigg\lvert\Bigg\langle v(\lambda^{(n)}),\begin{bmatrix}
	x_{1} \\
    \vdots \\
    x_{a}  \\
\end{bmatrix}\Bigg\rangle_{\mathcal{M}}\Bigg\rvert^{2}\Bigg)\\
&\leq \dfrac{1}{1-\|\lambda^{(n)}\|_{\infty}}\Bigg(2\|v(\lambda^{(n)})\|_{\mathcal{M}}\Bigg\|\begin{bmatrix}
	x_{1} \\
    \vdots \\
    x_{a}  \\
\end{bmatrix}\Bigg\|_{\mathcal{M}}+\|v(\lambda^{(n)})\|_{\mathcal{M}}^{2}\Bigg\|\begin{bmatrix}
	x_{1} \\
    \vdots \\
    x_{a}  \\
\end{bmatrix}\Bigg\|_{\mathcal{M}}^{2}\Bigg), 
\end{align*}
and so
\begin{align*}
    &\dfrac{1-\lvert \varphi(\lambda^{(n)})\rvert^{2}}{1-\|\lambda^{(n)}\|_{\infty}^{2}}\\
&\leq \dfrac{1}{1-\|\lambda^{(n)}\|_{\infty}}\Bigg(2\|\lambda^{(n)}-\tau\|_{\infty}\|u(\lambda^{(n)})\|_{\mathcal{M}}\|x\|_{\mathcal{M}}\\
&\hspace{5.2cm}+\|x\|^{2}_{\mathcal{M}}\sum_{j=1}^{a}\lvert\lambda^{(n)}_{j}-\tau_{j}\rvert^{2}\|u_{j}(\lambda^{(n)})\|^{2}_{\mathcal{M}_{j}}\Bigg)\\
&\leq \dfrac{1}{1-\|\lambda^{(n)}\|_{\infty}}\Bigg(2\|\lambda^{(n)}-\tau\|_{\infty}\|u(\lambda^{(n)})\|_{\mathcal{M}}\|x\|_{\mathcal{M}}\\
&\hspace{5.2cm}+\|x\|^{2}_{\mathcal{M}}\|u(\lambda^{(n)})\|^{2}_{\mathcal{M}}\sum_{j=1}^{a}\lvert\lambda^{(n)}_{j}-\tau_{j}\rvert^{2}\Bigg).\\
\end{align*}
By assumption (2), $\|u(\lambda^{(n)})\|_{\mathcal{M}}\leq M$. Therefore
\begin{align}
    &\dfrac{1-\lvert \varphi(\lambda^{(n)})\rvert^{2}}{1-\|\lambda^{(n)}\|_{\infty}^{2}} \nonumber\\
    &\leq\dfrac{1}{1-\|\lambda^{(n)}\|_{\infty}}\Bigg(2M\|\lambda^{(n)}-\tau\|_{\infty}\|x\|_{\mathcal{M}}+M^{2}\|x\|^{2}_{\mathcal{M}}\sum_{j=1}^{a}\lvert\lambda^{(n)}_{j}-\tau_{j}\rvert\cdot\lvert\lambda^{(n)}_{j}-\tau_{j}\rvert\Bigg) \nonumber\\
    &\leq 2M\|x\|_{\mathcal{M}}\dfrac{\|\lambda^{(n)}-\tau\|_{\infty}}{1-\|\lambda^{(n)}\|_{\infty}}+M^{2}\|x\|^{2}_{\mathcal{M}}\dfrac{\|\lambda^{(n)}-\tau\|_{\infty}}{1-\|\lambda^{(n)}\|_{\infty}}\sum_{j=1}^{a}\lvert\lambda^{(n)}_{j}-\tau_{j}\rvert. \label{bound4}
\end{align} 
Consider a set $S\subset \mathbb{D}^{d}$ that approaches $\tau$ in such a way that, for some $c\geq 1$, 
$$\|\tau-\lambda^{(n)}\|_{\infty}\leq c(1-\|\lambda^{(n)}\|_{\infty})~~\text{for all}~n. $$ 
Thus, for $(\lambda^{(n)})$ a sequence in $S$ tending to $\tau$, we also have that, for large enough $n$, $\lvert \lambda^{(n)}_{j}-\tau_{j}\rvert <1$ for all $j=1,\hdots, a$. Use both of these inequalities in  inequality \eqref{bound4} to see that if $(\lambda^{(n)})$ is a sequence in $S$ tending to $\tau$, for large enough $n$,
    $$\dfrac{1-\lvert \varphi(\lambda^{(n)})\rvert^{2}}{1-\|\lambda^{(n)}\|_{\infty}^{2}} \leq 2M\|x\|_{\mathcal{M}}c+M^{2}\|x\|^{2}_{\mathcal{M}}ca.$$   
Therefore  
$$\sup_{n} \dfrac{1-\lvert \varphi(\lambda^{(n)})\rvert}{1-\|\lambda^{(n)}\|_{\infty}} \leq 2 \sup_{n} \dfrac{1-\lvert \varphi(\lambda^{(n)})\rvert^{2}}{1-\|\lambda^{(n)}\|_{\infty}^{2}}  <\infty.$$
Thus, when $\tau\in\partial\mathbb{D}^{d}\setminus\mathbb{T}^{d}$, (2)$\implies$(1). Hence statements (1), (2) and (3) are equivalent.\end{proof}

\begin{remark}
\rm The following result  \cite[Remark 5.3d]{AMYcara} illustrates 
 the difference between parts (i) and (ii) of Proposition \ref{prop5.2cara}. Let $\varphi \in \mathcal{S}(\D^2)$. Whereas for $\tau\in\T^2$, $\tau$ is a B-point for $\ph$ if and only if $u(\lambda^{(n)})$ is bounded for some model $(\M,u)$ of $\ph$ and every sequence $(\lambda^{(n)})$ in $\D^2$ that converges nontangentially to $\tau$; for $\tau\in \partial \D^2 \setminus \T^2$ this equivalence does not hold.  For consider the example $\ph(\lambda) = \lambda_2$ for $\lambda=(\lambda_1,\lambda_2) \in \D^2$.
A model for $\ph$ is $(\M,u)$ where $\M={0}\oplus\C$ and $u:\D^2 \to \M$ is given by $u(\la)=(u_1(\la),u_2(\la))$ where $u_1(\la)=0$ and $u_2(\la)=1$ for all $\la\in\D^2$. The model relation becomes
\[
1-\overline{\ph(\mu)}\ph(\la) = (1-\overline{\mu_1} \la_1)0+ (1-\overline{\mu_2} \la_2) \langle{1},{1} \rangle_\C,
\]
which is clearly valid.  Here $u(\la)$ is bounded, but the point $\tau=(1,0)$ is not a B-point of $\ph$, for if $\la^{(n)}$ is any sequence in $\D^2$ that converges to $\tau$ we have $\la^{(n)}_1 \to 1, \la^{(n)}_2 \to 0$ and
\[
\frac{1-|\ph(\la^{(n)})|}{1-\|\la^{(n)}\|} = \frac{1-|\la^{(n)}_2|}{1-\max(|\la^{(n)}_1|,|\la^{(n)}_2|)},
\]
which is unbounded, since the numerator tends to $1$ and the denominator tends to $0$. 
\end{remark}

The following statement gives us a bound on the model $(\mathcal{M},u)$ of $\varphi\in\mathcal{SA}_{d}$ in terms of a function-theoretic property of $\varphi$.

\begin{corollary}\label{coro5.7cara}
Let $\varphi\in\mathcal{SA}_{d}$, let $\tau\in\partial\mathbb{D}^{d}$ and let $(\mathcal{M},u)$ be a model of $\varphi$. \\
{\rm (i)} If $ \ \tau \in \partial \D^d \setminus\T^d$, then  the following conditions are equivalent: 
\begin{enumerate}
    \item $\dfrac{1-\lvert \varphi(\lambda) \rvert}{1-\|\lambda\|_{\infty}}$ is bounded on some sequence $\{\lambda^{(n)}\}$ that converges to $\tau$;
    \item $\lim_{n \to \infty} |\varphi(\lambda^{(n)})| =1$ and $u(\lambda^{(n)})$ is bounded on some sequence $\{\lambda^{(n)}\}$ that approaches $\tau$ nontangentially;
    \item $\dfrac{1-\lvert \varphi(\lambda) \rvert}{1-\|\lambda\|_{\infty}}$ is bounded on every subset of $\mathbb{D}^{d}$ that approaches $\tau$ nontangentially;
    \item for every subset $S$ of $\D^d$ that approaches $\tau$ nontangentially,  $\lim_{\la\to\tau, \la\in S} |\varphi(\lambda)| =1$, and  $ u(\lambda) $ is  bounded  on $S$.
\end{enumerate} 
{\rm (ii)} If $ \ \tau \in \T^d$, then  the following conditions are equivalent: 
\begin{enumerate}
    \item $\dfrac{1-\lvert \varphi(\lambda) \rvert}{1-\|\lambda\|_{\infty}}$ is bounded on some sequence $\{\lambda^{(n)}\}$ that converges to $\tau$;
    \item $u(\lambda^{(n)})$ is bounded on some sequence $\{\lambda^{(n)}\}$ that approaches $\tau$ nontangentially;
    \item $\dfrac{1-\lvert \varphi(\lambda) \rvert}{1-\|\lambda\|_{\infty}}$ is bounded on every subset of $\mathbb{D}^{d}$ that approaches $\tau$ nontangentially;
    \item $u(\lambda)$ is bounded on every subset of $\mathbb{D}^{d}$ that approaches $\tau$ nontangentially.
\end{enumerate}
\end{corollary}
\begin{proof} By Corollary \ref{coro4.14cara}, (1)$\iff$(3).
By Corollary \ref{coro4.14cara} and  Proposition \ref{prop5.2cara},  (1)$\implies$(2). By  Proposition \ref{prop5.2cara},  (2)$\implies$(1). Thus (1)$\iff$(2).
By  Proposition \ref{prop5.4} and Lemma \ref{limit-mod},  (3)$\implies$(4). By  Proposition \ref{prop5.2cara},  (4)$\implies$(3). Thus (3)$\iff$(4).
It is clear that  (4)$\implies$(2) trivially. Finally, (1) is the definition that $\tau$ is a $B$-point for $\varphi$, so that (1)$\implies$(4) follows from  Proposition \ref{prop5.4} and Lemma \ref{limit-mod}.  \end{proof}

\begin{lemma}\label{ujnotxj}
Let $\varphi\in\mathcal{SA}_{d}$, let $(\mathcal{M},u)$ be a model of $\varphi$, where  $\mathcal{M}$ is a separable complex Hilbert space with an orthogonal decomposition $\mathcal{M}=\mathcal{M}_{1}\oplus\hdots\oplus\mathcal{M}_{d}$ and $u=(u_1, \dots, u_d)$ is an analytic map  $u_j: \mathbb{D}^{d} \to \mathcal{M}_{j} $ for $j =1, \dots, d$, such that condition 
\eqref{modeleqn}    holds, and let $\tau\in\partial\mathbb{D}^{d}$ be a $B$-point for $\varphi$. Suppose that $\{\lambda^{(n)}\}$ converges to $\tau$ nontangentially in $\mathbb{D}^{d}$ and $\{u(\lambda^{(n)})\}$ converges weakly to $x$ in $\mathcal{M}$ as $n\rightarrow \infty$.
Suppose also that, for some index $j$,  $\{u_j(\lambda^{(n)})\}$ does not tend to $x_j$  in norm as $n\rightarrow \infty$.
Then there exists a   subsequence $(\lambda^{(n_k)})$ of $(\lambda^{(n)})$ such that, 
\begin{equation} \label{lemma-uj-0}  \|u_{j}(\lambda^{(n_k)})\|_{\mathcal{M}_{j}}^{2}-\mathrm{Re} \langle u_{j}(\lambda^{(n_k)}),x_{j}\rangle_{\mathcal{M}_{j}} \geq 0 \;\ \text{for all} \; \ k \ge 1.
\end{equation}
\end{lemma}
\begin{proof}
Since
$\|u_{j}(\lambda^{(n)}) - x_j\|_{\mathcal{M}_{j}} \; \text{ does not tend to} \; 0 \; \text{as} \; n\rightarrow \infty,$
 there exist a strictly increasing   sequence $(n_k)_{k \geq 1}$ in $\mathbb{N}$ and a strictly positive number $C$ such that 
$$\|u_{j}(\lambda^{(n_k)})- x_j \|_{\mathcal{M}_{j}}^{2}
 \geq C \;\ \text{for all} \; \ k \ge 1.$$ 
Thus, for all $  k \ge 1$,
\begin{equation} \label{lemma-uj-1} 
 \|u_{j}(\lambda^{(n_k)})-x_{j}\|_{\mathcal{M}_{j}}^{2}=\|u_{j}(\lambda^{(n_k)})\|_{\mathcal{M}_{j}}^{2}-2 \mathrm{Re}
 \langle u_{j}(\lambda^{(n_k)}),x_{j}\rangle_{\mathcal{M}_{j}} +\|x_{j}\|_{\mathcal{M}}^{2} \geq C>0.
\end{equation} 
Since $\{u(\lambda^{(n)})\} \to x$  weakly in $\mathcal{M}$ as $n\rightarrow \infty$, it follows that, for each index $j= 1,2, \dots, d$,  $\{u_{j}(\lambda^{(n)})\} \to x_{j}$ weakly in $\mathcal{M}_{j}$ as $n\rightarrow \infty$, and so 
$$\langle u_{j}(\lambda^{(n_k)}),x_{j}\rangle_{\mathcal{M}_{j}} \to \|x_{j}\|_{\mathcal{M}_{j}}^{2}  \; \text { as} \; k \to \infty.
$$
Hence, there exists a $N_0 \in  \mathbb{N}$ such that 
\begin{equation} \label{lemma-uj-2} 
\Big| \|x_{j}\|_{\mathcal{M}_{j}}^{2} - \langle u_{j}(\lambda^{(n_k)}),x_{j}\rangle_{\mathcal{M}_{j}}  \Big| \leq \half C \; \text{ for all } \; k \ge N_0 .
\end{equation} 
Combine the inequalities \eqref{lemma-uj-1} and \eqref{lemma-uj-2} to obtain, for all  $k \ge N_0 $,
\begin{align}\label{lemma-uj-3} 
&\|u_{j}(\lambda^{(n_k)})\|_{\mathcal{M}_{j}}^{2} -\mathrm{Re} \langle u_{j}(\lambda^{(n_k)}),x_{j}\rangle_{\mathcal{M}_{j}}\nonumber\\  
&\hspace{3cm}=\|u_{j}(\lambda^{(n_k)})\|_{\mathcal{M}_{j}}^{2}-2 \mathrm{Re} \langle u_{j}(\lambda^{(n_k)}),x_{j}\rangle_{\mathcal{M}_{j}} +\|x_{j}\|_{\mathcal{M}}^{2} \nonumber\\ 
&\hspace{6cm}+\mathrm{Re} \langle u_{j}(\lambda^{(n_k)}),x_{j}\rangle_{\mathcal{M}_{j}} -\|x_{j}\|_{\mathcal{M}}^{2} \nonumber\\
&\hspace{3cm} \geq C -\half C>0.
\end{align}
Thus the subsequence $(\lambda^{(n_{k + N_0})})_{k \ge 1}$ of $(\lambda^{(n)})_{n \geq 1}$ has the stated property 
\eqref{lemma-uj-0}.
\end{proof}

\begin{proposition}\label{prop5.8cara}
Let $\varphi\in\mathcal{SA}_{d}$, let $(\mathcal{M},u)$ be a model of $\varphi$,  where  $\mathcal{M}$ is a separable complex Hilbert space with an orthogonal decomposition $\mathcal{M}=\mathcal{M}_{1}\oplus\hdots\oplus\mathcal{M}_{d}$, 
 and let $\tau\in\partial\mathbb{D}^{d}$ be a $B$-point for $\varphi$. Suppose that $\{\lambda^{(n)}\}$ converges to $\tau$ nontangentially in $\mathbb{D}^{d}$. If $\{u(\lambda^{(n)})\}$ converges weakly in $\mathcal{M}$ as $n\rightarrow \infty$, then $\{u(\lambda^{(n)})\}$ converges in norm as $n\rightarrow \infty$.
\end{proposition}
\begin{proof} Suppose  $u(\lambda^{(n)})\rightarrow x$ weakly in $\mathcal{M}$.  The equality
$$\|u(\lambda^{(n)})-x\|_{\mathcal{M}}^{2} = \langle u(\lambda^{(n)})-x,u(\lambda^{(n)})-x\rangle_{\mathcal{M}}$$
implies that
\begin{align}
\|u(\lambda^{(n)})-x\|_{\mathcal{M}}^{2}&=\|u(\lambda^{(n)})\|_{\mathcal{M}}^{2}-2\mathrm{Re} \langle u(\lambda^{(n)}),x\rangle_{\mathcal{M}} +\|x\|_{\mathcal{M}}^{2}\nonumber\\
&=\bigg(\|u(\lambda^{(n)})\|_{\mathcal{M}}^{2}-\mathrm{Re} \langle u(\lambda^{(n)}),x\rangle_{\mathcal{M}} \bigg)\nonumber\\
&\hspace{3cm}+\bigg(\|x\|_{\mathcal{M}}^{2}-\mathrm{Re} \langle u(\lambda^{(n)}),x\rangle_{\mathcal{M}} \bigg)\label{2ndequationinthisalign}.
\end{align}
Since $u(\lambda^{(n)})\rightarrow x$ weakly as $n\rightarrow\infty$, 
$$\langle u(\lambda^{(n)}),x \rangle_{\mathcal{M}} \rightarrow \langle x,x\rangle_{\mathcal{M}}~~\text{as}~n\rightarrow \infty.$$
Thus
\begin{equation}\label{equivalenceofweakconvergence}
 \|x\|_{\mathcal{M}}^{2}-\mathrm{Re} \langle u(\lambda^{(n)}),x\rangle_{\mathcal{M}} \rightarrow 0~\text{as}~n\rightarrow \infty.
\end{equation}
To show that $u(\lambda^{(n)})$ converges in norm to $x$ as $n\rightarrow \infty$, all that remains is to show that the first term on the right hand side of equation (\ref{2ndequationinthisalign}), that is, 
\begin{equation}\label{ujtoxj}
\sum_{j=1}^{d}\Big(\|u_{j}(\lambda^{(n)})\|_{\mathcal{M}_{j}}^{2}-\mathrm{Re} \langle u_{j}(\lambda^{(n)}),x_{j}\rangle \Big) \to 0
\end{equation}
as  $n\rightarrow \infty$.
Since $\{u(\lambda^{(n)})\} \to x$ weakly in $\mathcal{M}$ as $n\rightarrow \infty$, it follows that, for each index $j= 1,2, \dots, d$,  $\{u_{j}(\lambda^{(n)})\} \to x_{j}$ weakly in $\mathcal{M}_{j}$ as $n\rightarrow \infty$.
 If for each $j$, $\|u_{j}(\lambda^{(n)}) - x_j\| \to 0$  as $n\rightarrow \infty$, then $u(\lambda^{(n)})$ converges in norm to $x$ as $n\rightarrow \infty$. Suppose that there are some $j$ such that 
$\|u_{j}(\lambda^{(n)}) - x_j\|$ does not tend to $0$  as $n\rightarrow \infty$. Let us show that in this case we will get a contradiction.

 By Proposition \ref{phi-nt-omega},
there exists $\omega\in\mathbb{T}$ such that $\varphi(\lambda)\rightarrow \omega~~\text{as}~~\lambda{\overset{nt}\rightarrow}\tau.$ By Proposition \ref{prop4.2generalized} (3), 
\begin{equation}\label{omega-model}
  1-\overline{\omega}\varphi(\lambda)=\sum_{\lvert \tau_{j} \rvert = 1}\Big(1-\overline{\tau}_{j}\lambda_{j}\Big)\langle u_{j}(\lambda),x_{j}\rangle_{\mathcal{M}_{j}}
  \end{equation}
holds for all $\lambda\in\mathbb{D}^{d}$.
By Proposition \ref{prop4.2generalized} (2), if $(x_{1},\hdots,x_{d})=x\in Y_{\tau}$ where $\tau\in\partial\mathbb{D}^{d}$ is such that  $(\tau_{1},\hdots,\tau_{d})=\tau\in\mathbb{T}^{a}\times\mathbb{D}^{b}$ where $a+b=d$, then for each $\tau_{j}\in\mathbb{D}$, $x_{j}=0$.

Thus,  by equation \eqref{omega-model},
   for $\lambda^{(n)}\in\mathbb{D}^{d}$, $n=1,2,\dots$,
\begin{equation}\label{prop4.2equationhere}
1-\overline{\omega}\varphi(\lambda^{(n)})=\sum_{j=1}^{d}\Big(1-\overline{\tau}_{j}\lambda^{(n)}_{j}\Big)\langle u_{j}(\lambda^{(n)}),x_{j}\rangle_{\mathcal{M}_{j}}.
\end{equation}
Hence,
\begin{equation}\label{REprop4.2equationhere}
\mathrm{Re}\bigg(1-\overline{\omega}\varphi(\lambda^{(n)})\bigg)=\mathrm{Re}\bigg(\sum_{j=1}^{d}\Big(1-\overline{\tau}_{j}\lambda^{(n)}_{j}\Big)\langle u_{j}(\lambda^{(n)}),x_{j}\rangle_{\mathcal{M}_{j}}\bigg).
\end{equation}
By equations (\ref{whenlambdametmu}) and (\ref{REprop4.2equationhere}),
\begin{align}
&\sum_{j=1}^{d}(1-\lvert \lambda_{j}^{(n)} \rvert^{2})\|u_{j}(\lambda^{(n)})\|_{\mathcal{M}_{j}}^{2}\nonumber\\
&= 1-\lvert \varphi(\lambda^{(n)}) \rvert^{2} \nonumber\\
&= 1-\lvert \varphi(\lambda^{(n)}) \rvert^{2}-2\mathrm{Re}\big(1-\overline{\omega}\varphi(\lambda^{(n)})\big)\nonumber\\
&\hspace{4cm}+ \mathrm{Re}\bigg(\sum_{j=1}^{d}2\Big(1-\overline{\tau}_{j}\lambda^{(n)}_{j}\Big)\langle u_{j}(\lambda^{(n)}),x_{j}\rangle_{\mathcal{M}_{j}}\bigg).
\label{rearrangmentofthemodelhere}
\end{align}
Subtraction of  $\; \sum_{j=1}^{d}(1-\lvert\lambda^{(n)}_{j}\rvert^{2})\mathrm{Re} \langle u_{j}(\lambda^{(n)}),x_{j}\rangle_{\mathcal{M}_{j}} \;$ from both sides of equation (\ref{rearrangmentofthemodelhere}) gives us 
\begin{align}
&\sum_{j=1}^{d}(1-\lvert \lambda^{(n)}_{j} \rvert^{2})\Big(\|u_{j}(\lambda^{(n)})\|_{\mathcal{M}_{j}}^{2}-\mathrm{Re}\langle u_{j}(\lambda^{(n)}),x_{j}\rangle  \Big)\nonumber\\
&=1-\lvert \varphi(\lambda^{(n)}) \rvert^{2}-2\mathrm{Re}\big(1-\overline{\omega}\varphi(\lambda^{(n)})\big)\nonumber\\
&\hspace{2cm}+ \mathrm{Re}\bigg(\sum_{j=1}^{d}2\Big(1-\overline{\tau}_{j}\lambda^{(n)}_{j}\Big)\langle u_{j}(\lambda^{(n)}),x_{j}\rangle_{\mathcal{M}_{j}}\bigg) - \sum_{j=1}^{d}(1-\lvert\lambda^{(n)}_{j}\rvert^{2})
\mathrm{Re} \langle u_{j}(\lambda^{(n)}),x_{j}\rangle_{\mathcal{M}_{j}} \nonumber \\
&=-1+2\mathrm{Re}\big(\overline{\omega}\varphi(\lambda^{(n)})\big)-\lvert \varphi(\lambda^{(n)}) \rvert^{2}\nonumber \\
&\hspace{2cm}+ \mathrm{Re}\bigg[\sum_{j=1}^{d}\Big(2(1-\overline{\tau}_{j}\lambda^{(n)}_{j})-(1-\lvert \lambda^{(n)}_{j}\rvert^{2})\Big)\langle u_{j}(\lambda^{(n)}),x_{j}\rangle_{\mathcal{M}_{j}}\bigg]. \label{Prop4.5-01}
\end{align} 
Note that
$$\lvert1-\overline{\omega}\varphi(\lambda^{(n)})\rvert^{2}= 1-2\mathrm{Re}\big(\overline{\omega}\varphi(\lambda^{(n)})\big)+\lvert \varphi(\lambda^{(n)}) \rvert^{2}, $$
and so by equation \eqref{Prop4.5-01},
\begin{align}
&\sum_{j=1}^{d}(1-\lvert \lambda^{(n)}_{j} \rvert^{2})\Big(\|u_{j}(\lambda^{(n)})\|_{\mathcal{M}_{j}}^{2}-
\mathrm{Re} \langle u_{j}(\lambda^{(n)}),x_{j}\rangle \Big)\nonumber\\
&=-\lvert1-\overline{\omega}\varphi(\lambda^{(n)})\rvert^{2}+\mathrm{Re}\Big[\sum_{j=1}^{d}\Big(1-2\overline{\tau}_{j}\lambda^{(n)}_{j}+\lvert\lambda^{(n)}_{j}\rvert^{2}\Big)\|x_{j}\|_{\mathcal{M}_{j}}^{2}\Big]\nonumber\\
&\hspace{2.5cm}+\mathrm{Re}\bigg[\sum_{j=1}^{d}\Big(1-2\overline{\tau}_{j}\lambda^{(n)}_{j}+\lvert\lambda^{(n)}_{j}\rvert^{2}\Big)\langle u_{j}(\lambda^{(n)})-x_{j},x_{j}\rangle_{\mathcal{M}_{j}}\bigg].\label{lastequationhereinthispart}
\end{align} 
Note that, for $\lambda\in\mathbb{D}^{d}$, $\lvert \lambda_{j}\rvert^{2}\leq \lvert \lambda_{j} \rvert \leq \|\lambda\|_{\infty} <1$ for $j=1,\hdots, d$, and so
$$ 0< \dfrac{1}{1-\lvert \lambda_{j}\rvert^{2}} \leq \dfrac{1}{1-\|\lambda\|_{\infty}} < \infty.$$
By assumption, $\lambda^{(n)}{\overset{nt}\rightarrow}\tau$ as $n\rightarrow \infty$. Hence there exists a $c>0$ such that
$$\|\lambda^{(n)} - \tau\|_{\infty}\leq c(1-\|\lambda^{(n)}\|_{\infty}) \;\; \text{for all } \; n \in \mathbb{N}.$$
Consequently, for all $ n \in \mathbb{N}$,
$$\dfrac{\|\lambda^{(n)}-\tau\|_{\infty}}{1-\lvert \lambda^{(n)}_{j} \rvert^{2}}\leq \dfrac{\|\lambda^{(n)}-\tau\|_{\infty}}{1-\|\lambda^{(n)}\|_{\infty}}\leq c,$$
so that
\begin{equation}\label{inequalitywithc}
\dfrac{c \big(1-\lvert \lambda^{(n)}_{j}\rvert^{2}\big)}{\|\lambda^{(n)}-\tau\|_{\infty}}\geq 1.
\end{equation}
It is easy to see that, for all $\lambda\in\mathbb{D}^{d}$ and $\tau \in\partial\mathbb{D}^{d}$,
$$1-\|\lambda\|_{\infty}\leq \|\lambda-\tau\|_{\infty},$$
and so
\begin{equation}\label{inequalitywithc-5}
\dfrac{1}{\|\lambda-\tau\|_{\infty}} \leq \dfrac{1}{1-\|\lambda\|_{\infty}}.
\end{equation}
Note that, for $j=1,\hdots, d$,
$$1-\lvert \lambda_{j}\rvert^{2} \leq \|\tau+\lambda\|_{\infty}\| \tau-\lambda\|_{\infty}.$$
Let $c_{1}=\max_{i=1,\hdots d} \lvert \tau_{i}+\lambda_{i}\rvert = \|\tau+\lambda\|_{\infty}$. Then, for $j=1,\hdots, d$, we have the inequality
\begin{equation}\label{ineqc1}
1-\lvert \lambda^{(n)}_{j}\rvert^{2} \leq c_{1}\|\tau -\lambda^{(n)}\|_{\infty}.
\end{equation}
Equivalently, for $j=1,\hdots, d$,
\begin{equation}\label{inequalityforc1}
\dfrac{1}{\|\tau-\lambda^{(n)}\|_{\infty}} \leq \dfrac{c_{1}}{1-\lvert \lambda^{(n)}_{j} \rvert^{2}}.
\end{equation}
Multiply equation \eqref{lastequationhereinthispart} by $\dfrac{c}{\|\lambda^{(n)}-\tau\|_{\infty}}$, to obtain the following equality, 
\begin{align}\label{equation527}
& \dfrac{c}{\|\lambda^{(n)}-\tau\|_{\infty}}\cdot\sum_{j=1}^{d}(1-\lvert \lambda^{(n)}_{j} \rvert^{2})\Big(\|u_{j}(\lambda^{(n)})\|_{\mathcal{M}_{j}}^{2}-\mathrm{Re} \langle u_{j}(\lambda^{(n)}),x_{j}\rangle_{\mathcal{M}_{j}} \Big)\nonumber\\
&=\dfrac{c}{\|\lambda^{(n)}-\tau\|_{\infty}} \Bigg(-\lvert1-\overline{\omega}\varphi(\lambda^{(n)})\rvert^{2}+\sum_{j=1}^{d}\mathrm{Re}\Big(1-2\overline{\tau}_{j}\lambda^{(n)}_{j}+\lvert\lambda^{(n)}_{j}\rvert^{2}\Big)\|x_{j}\|_{\mathcal{M}_{j}}^{2}
\nonumber\\
&\hspace{1cm}+\mathrm{Re}\bigg[\sum_{j=1}^{d}\Big(1-2\overline{\tau}_{j}\lambda^{(n)}_{j}+\lvert\lambda^{(n)}_{j}\rvert^{2}\Big)\langle u_{j}(\lambda^{(n)})-x_{j},x_{j}\rangle_{\mathcal{M}_{j}}\bigg]\Bigg)\nonumber\\
&=c\Bigg(\dfrac{-\lvert1-\overline{\omega}\varphi(\lambda^{(n)})\rvert^{2}}{\|\lambda^{(n)}-\tau\|_{\infty}}+\sum_{j=1}^{d}\mathrm{Re}\Big(1-2\overline{\tau}_{j}\lambda^{(n)}_{j}+\lvert\lambda^{(n)}_{j}\rvert^{2}\Big)\dfrac{\|x_{j}\|_{\mathcal{M}_{j}}^{2}}{\|\lambda^{(n)}-\tau\|_{\infty}}\nonumber\\
&\hspace{1cm}+\sum_{j=1}^{d}\dfrac{1}{\|\lambda^{(n)}-\tau\|_{\infty}}\mathrm{Re}\bigg[
\Big(1-2\overline{\tau}_{j}\lambda^{(n)}_{j}+\lvert\lambda^{(n)}_{j}\rvert^{2}\Big)\langle u_{j}(\lambda^{(n)})-x_{j},x_{j}\rangle_{\mathcal{M}_{j}}\bigg]\Bigg).
\end{align}

There are two cases: (i) $\tau\in\mathbb{T}^{d}$ and (ii) $\tau\in\partial\mathbb{D}^{d}\setminus\mathbb{T}^{d}$.

{\rm (i)}. Assume firstly that $\tau\in\mathbb{T}^{d}$. Then, for each $j$, 
\begin{equation}\label{inequalityfor527}
\mathrm{Re}\Big( 1-2\overline{\tau}_{j}\lambda^{(n)}_{j}+\lvert\lambda^{(n)}_{j}\rvert^{2} \Big)= |{\tau}_{j} -\lambda^{(n)}_{j}|^2,
\end{equation}
and 
\begin{align}\label{inequalityfor527-2}
&\Big|1-2\overline{\tau}_{j}\lambda^{(n)}_{j}+\lvert\lambda^{(n)}_{j}\rvert^{2}\Big| = \Big| \overline{\tau}_{j} {\tau}_{j}- \overline{{\tau}_{j}}\lambda^{(n)}_{j} -\overline{{\tau}_{j}}\lambda^{(n)}_{j} + \overline{\lambda^{(n)}_{j}}\lambda^{(n)}_{j}\Big|\nonumber\\
&\leq \Big|{\tau}_{j} -\lambda^{(n)}_{j}\Big| + \Big|\lambda_{j}^{(n)} \Big|\Big|\overline{{\tau}_{j}} -\overline{\lambda^{(n)}_{j}}\Big| = \Big|{\tau}_{j} -\lambda^{(n)}_{j}\Big| \Big(1 + \Big|\lambda_{j}^{(n)} \Big| \Big).
\end{align} 
Thus, by equations  (\ref{equation527}), \eqref{inequalityfor527} and inequalities \eqref{inequalityfor527-2} and \eqref{inequalitywithc-5},
\begin{align}\label{equation5271}
& \Bigg\lvert \dfrac{c}{\|\lambda^{(n)}-\tau\|_{\infty}}\cdot\sum_{j=1}^{d}(1-\lvert \lambda^{(n)}_{j} \rvert^{2})\Big(\|u_{j}(\lambda^{(n)})\|_{\mathcal{M}_{j}}^{2}-\mathrm{Re} \langle u_{j}(\lambda^{(n)}),x_{j}\rangle_{\mathcal{M}_{j}} \Big)
\Bigg\rvert \nonumber\\
&\leq c\Bigg(\dfrac{\lvert 1-\overline{\omega}\varphi(\lambda^{(n)})\rvert^{2}}{\|\lambda^{(n)}-\tau\|_{\infty}}+
\sum_{j=1}^{d}\dfrac{\Big[|{\tau}_{j} -\lambda^{(n)}_{j}|^2 
\|x_{j}\|_{\mathcal{M}_{j}}^{2}\Big]}{\|\lambda^{(n)}-\tau\|_{\infty}}\nonumber\\
&\hspace{4cm}+\sum_{j=1}^{d}\dfrac{\Big[
{\lvert \tau_{j}-\lambda^{(n)}_{j}\rvert}{(1+\lvert\lambda^{(n)}_{j})}
|\langle u_{j}(\lambda^{(n)})-x_{j},x_{j}\rangle_{\mathcal{M}_{j}}|\Big]}{\|\lambda^{(n)}-\tau\|_{\infty}}\nonumber\\
&\leq c\Bigg(\dfrac{\lvert 1-\overline{\omega}\varphi(\lambda^{(n)})\rvert^{2}}{1-\|\lambda^{(n)}\|_{\infty}}
+\|x\|_{\mathcal{M}}^{2}\sum_{j=1}^{d}\dfrac{\Big[|{\tau}_{j} -\lambda^{(n)}_{j}|^2 \Big]}{\|\lambda^{(n)}-\tau\|_{\infty}} 
\nonumber\\
&\hspace{4cm}+\sum_{j=1}^{d}\dfrac{\Big[
{\lvert \tau_{j}-\lambda^{(n)}_{j}\rvert}{(1+\lvert\lambda^{(n)}_{j})}
\Big]}{\|\lambda^{(n)}-\tau\|_{\infty}}\lvert\langle u_{j}(\lambda^{(n)})-x_{j},x_{j}\rangle_{\mathcal{M}}\rvert \Bigg).
\end{align} 
The first term on the right tends to zero because, firstly, by inequality \eqref{phi-omega} from Theorem \ref{theorem4.9generalized},
\begin{equation}\label{phi-omega-2}
 \lvert \varphi(\lambda)-\omega\rvert^{2} \leq  \|x\|_{\mathcal{M}}^{2}\bigg(\max_{j=1,\hdots,d}\dfrac{\lvert \tau_{j}-\lambda_{j}\rvert^{2}}{1-\lvert \lambda_{j}\rvert^{2}}\bigg)(1-\lvert \varphi(\lambda) \rvert^{2}),
\end{equation}
 secondly, by Corollary \ref{coro5.7cara},   $\dfrac{1-\lvert \varphi(\lambda) \rvert}{1-\|\lambda\|_{\infty}}$ is bounded on every subset of $\mathbb{D}^{d}$ that approaches $\tau$ nontangentially and because $\lambda^{(n)}{\overset{nt}\rightarrow}\tau$. 
The second term tends to zero because $\lambda^{(n)}{\overset{nt}\rightarrow}\tau$. The third term tends to zero because
each summand is the product of a bounded factor and a factor that tends to zero, since $u(\lambda^{(n)})\rightarrow x$ weakly in $\mathcal{M}$. Therefore,
\begin{equation}\label{equation5272}
 \dfrac{c}{\|\lambda^{(n)}-\tau\|_{\infty}}\cdot\sum_{j=1}^{d}(1-\lvert \lambda^{(n)}_{j} \rvert^{2})\Big(\|u_{j}(\lambda^{(n)})\|_{\mathcal{M}_{j}}^{2}-\mathrm{Re} \langle u_{j}(\lambda^{(n)}),x_{j}\rangle_{\mathcal{M}_{j}} \Big) \to 0
\end{equation}
as $n \rightarrow  \infty$.

Let $$\Lambda = \{ j: \ \text{such that} \ \|u_{j}(\lambda^{(n)}) - x_j\| \; \text{ does not tend to} \; 0 \; \text{as} \; n\rightarrow \infty \},$$
and let $$\Lambda'= \{ j: \ \text{such that} \ \|u_{j}(\lambda^{(n)}) - x_j\| \; \text{tends to} \; 0 \; \text{as} \; n \rightarrow \infty \}.$$
 By Lemma \ref{ujnotxj} and by repeatedly passing to subsequences, we can choose a strictly increasing sequence 
 $(n_k)_{k \ge 1}$ such that,  for each $j \in \Lambda$, 
$$\Big(\|u_{j}(\lambda^{(n_k)})\|_{\mathcal{M}_{j}}^{2}-\mathrm{Re} \langle u_{j}(\lambda^{(n_k)}),x_{j}\rangle_{\mathcal{M}_{j}} \Big) \geq 0  \; \text{ for all } \; k \ge1.$$
By inequality (\ref{inequalitywithc}), 
\begin{align}
& \sum_{j \in \Lambda}\Big(\|u_{j}(\lambda^{(n_k)})\|_{\mathcal{M}_{j}}^{2}-\mathrm{Re} \langle u_{j}(\lambda^{(n_k)}),x_{j}\rangle_{\mathcal{M}_{j}} \Big) \nonumber\\
& \leq \sum_{j \in \Lambda} \dfrac{c (1-\lvert \lambda^{(n_k)}_{j} \rvert^{2})}{\|\lambda^{(n_k)}-\tau\|_{\infty}} 
\Big(\|u_{j}(\lambda^{(n_k)})\|_{\mathcal{M}_{j}}^{2}-\mathrm{Re} \langle u_{j}(\lambda^{(n_k)}),x_{j}\rangle_{\mathcal{M}_{j}} \Big)\nonumber\\
&= \dfrac{c}{\|\lambda^{(n_k)}-\tau\|_{\infty}}\sum_{j \in \Lambda}(1-\lvert \lambda^{(n_k)}_{j} \rvert^{2})\Big(\|u_{j}(\lambda^{(n_k)})\|_{\mathcal{M}_{j}}^{2}-\mathrm{Re} \langle u_{j}(\lambda^{(n_k)}),x_{j}\rangle_{\mathcal{M}_{j}} \Big).\label{oneovercinequality}
\end{align}
Therefore, by inequality \eqref{oneovercinequality} and relation \eqref{equation5272},
\begin{align}
& \sum_{j \in \Lambda}\Big(\|u_{j}(\lambda^{(n_k)})\|_{\mathcal{M}_{j}}^{2}-\mathrm{Re} \langle u_{j}(\lambda^{(n_k)}),x_{j}\rangle_{\mathcal{M}_{j}} \Big) \nonumber\\
& \hspace{1cm} +\dfrac{c}{\|\lambda^{(n_k)}-\tau\|_{\infty}}\sum_{j \in \Lambda'}(1-\lvert \lambda^{(n_k)}_{j} \rvert^{2})\Big(\|u_{j}(\lambda^{(n_k)})\|_{\mathcal{M}_{j}}^{2}-\mathrm{Re} \langle u_{j}(\lambda^{(n_k)}),x_{j}\rangle_{\mathcal{M}_{j}} \Big)\nonumber\\
&\leq \dfrac{c}{\|\lambda^{(n_k)}-\tau\|_{\infty}}\cdot\sum_{j=1}^{d}(1-\lvert \lambda^{(n_k)}_{j} \rvert^{2})\Big(\|u_{j}(\lambda^{(n_k)})\|_{\mathcal{M}_{j}}^{2}-\mathrm{Re} \langle u_{j}(\lambda^{(n_k)}),x_{j}\rangle_{\mathcal{M}_{j}} \Big)
 \to 0  \label{oneovercinequality2}
\end{align}
 as $ n_k \rightarrow  \infty$.
Since, for each $j \in \Lambda'$, $\|u_{j}(\lambda^{(n)}) - x_j\| \rightarrow 0 $ as $n \rightarrow \infty,$
by inequality \eqref{ineqc1},
$$\dfrac{c}{\|\lambda^{(n_k)}-\tau\|_{\infty}}\sum_{j \in \Lambda'}(1-\lvert \lambda^{(n_k)}_{j} \rvert^{2})\Big(\|u_{j}(\lambda^{(n_k)})\|_{\mathcal{M}_{j}}^{2}-\mathrm{Re} \langle u_{j}(\lambda^{(n_k)}),x_{j}\rangle_{\mathcal{M}_{j}} \Big)
\rightarrow 0$$
 as  $n_k\rightarrow \infty$. Hence
$$ \sum_{j \in \Lambda}\Big(\|u_{j}(\lambda^{(n_k)})\|_{\mathcal{M}_{j}}^{2}-\mathrm{Re} \langle u_{j}(\lambda^{(n_k)}),x_{j}\rangle_{\mathcal{M}_{j}} \Big)\rightarrow 0 \; \text {as} \; n_k \rightarrow \infty,$$
contrary to our assumption. Therefore, for each $j$, 
$$\|u_{j}(\lambda^{(n)}) - x_j\| \rightarrow 0 \; \text{as} \; n \rightarrow \infty,$$
and so
 $$\|u(\lambda^{(n)})\|_{\mathcal{M}}^{2}-\mathrm{Re} \langle u(\lambda^{(n)}),x\rangle_{\mathcal{M}} \rightarrow 0 \; \text {as} \; n \rightarrow \infty.$$  
Therefore,  $u(\lambda^{(n)})$ converges in norm to $x$ as $n \rightarrow \infty$ when $\tau\in\mathbb{T}^{d}$.

 {\rm (ii)}. For the second case, let $a+b=d$ such that $\tau\in\mathbb{T}^{a}\times \mathbb{D}^{b}$. By Proposition \ref{prop4.2generalized}, we have for $j=b,\hdots, d$, $x_{j}=0$. Indeed, as 
$$1-\lvert \varphi(\lambda^{(n)})\rvert^{2}=\sum_{j=1}^{d}\big(1-\lvert \lambda_{j}^{(n)}\rvert^{2}\big)\|u_{j}(\lambda^{(n)})\|_{\mathcal{M}_{j}}^{2},$$
letting $\lambda^{(n)}$ tend to $\tau$ we get $1-\lvert \varphi(\lambda^{(n)})\rvert^{2}$ tends to $0$ and 
$\big(1-\lvert \lambda_{j}^{(n)}\rvert^{2}\big)$ tends to $1 -|\tau_j|^2 \neq 0$ for  all $j=b,\hdots, d$. Hence
 $\|u_{j}(\lambda^{(n)})\|_{\mathcal{M}_{j}}$ tends to zero for all $j=b,\hdots, d$. 
It is well known that in a Hilbert space, a sequence converges in norm if and only if it converges weakly and the norms converge in $\mathbb{R}$. Hence, for each $j$ such that $ |\tau_j| < 1$, $u_{j}(\lambda^{(n)})$ tends to  $x_{j}=0$ in norm in ${\mathcal{M}_{j}}$.

 Suppose again that there are some $j$ such that 
$\|u_{j}(\lambda^{(n)}) - x_j\|$ does not tend to $0$  as $n\rightarrow \infty$. Let us show that in this case we will also get a contradiction.
Let $$\Lambda = \{ j: \ \text{such that} \ \|u_{j}(\lambda^{(n)}) - x_j\| \; \text{ does not tend to} \; 0 \; \text{as} \; n\rightarrow \infty \},$$
and let $$\Lambda'= \{ j: \ \text{such that} \ \|u_{j}(\lambda^{(n)}) - x_j\| \; \text{tends to} \; 0 \; \text{as} \; n \rightarrow \infty \}.$$
We have shown that $\{ j : |\tau_j| < 1 \}$ is a subset of $\Lambda'$. By equality \eqref{equation527},
since $x_j=0$ for each $j$ such that $|\tau_j| < 1$, we have the following relation:
\begin{align}\label{equation527notTd}
& \dfrac{c}{\|\lambda^{(n)}-\tau\|_{\infty}}\cdot\sum_{j=1}^{d}(1-\lvert \lambda^{(n)}_{j} \rvert^{2})\Big(\|u_{j}(\lambda^{(n)})\|_{\mathcal{M}_{j}}^{2}-\mathrm{Re} \langle u_{j}(\lambda^{(n)}),x_{j}\rangle_{\mathcal{M}_{j}} \Big)\nonumber\\
&=c\Bigg(\dfrac{-\lvert1-\overline{\omega}\varphi(\lambda^{(n)})\rvert^{2}}{\|\lambda^{(n)}-\tau\|_{\infty}}+\sum_{j: |\tau_j| =1}\mathrm{Re}\Big(1-2\overline{\tau}_{j}\lambda^{(n)}_{j}+\lvert\lambda^{(n)}_{j}\rvert^{2}\Big)\dfrac{\|x_{j}\|_{\mathcal{M}_{j}}^{2}}{\|\lambda^{(n)}-\tau\|_{\infty}}\nonumber\\
&\hspace{1cm}+\sum_{j: |\tau_j| =1}\dfrac{1}{\|\lambda^{(n)}-\tau\|_{\infty}}\mathrm{Re}\bigg[
\Big(1-2\overline{\tau}_{j}\lambda^{(n)}_{j}+\lvert\lambda^{(n)}_{j}\rvert^{2}\Big)\langle u_{j}(\lambda^{(n)})-x_{j},x_{j}\rangle_{\mathcal{M}_{j}}\bigg]\Bigg).
\end{align}
If we simply repeat the calculations as before in  case (i), we produce the same result, that is, for $\tau\in\mathbb{T}^{a}\times\mathbb{D}^{b}$, for each $j$,
$$\|u_{j}(\lambda^{(n)})- x_{j}\|_{\mathcal{M}_{j}}\rightarrow 0~~\text{as}~~\lambda^{(n)}{\overset{nt}\rightarrow}\tau.$$
Therefore, for $\tau\in\partial\mathbb{D}^{d}$, $u(\lambda^{(n)})$ converges in norm to $x$ as $n \rightarrow \infty$.
\end{proof}

\begin{theorem}\label{theorem5.12cara}
Let $\varphi\in\mathcal{SA}_{d}$, let $\tau\in\partial\mathbb{D}^{d}$ be a $B$-point for $\varphi$ and let $(\mathcal{M},u)$ be a model of $\varphi$. For any realization $(a,\beta,\gamma, D)$ of $(\mathcal{M},u)$ there exists a unique vector $u(\tau)\in\mathcal{M}$ such that $(1_{\mathcal{M}}-D\tau_{P})u(\tau)=\gamma$ and $u(\tau)\perp\mathrm{Ker}(1_{\mathcal{M}}-D\tau_{P})$ where $P$ is the $d$-tuple $P=(P_{1},\hdots,P_{d})$ on $\mathcal{M}$ such that $P_{j}:\mathcal{M}\rightarrow\mathcal{M}$ is the orthogonal projection onto $\mathcal{M}_{j}$ for $j=1,\hdots,d$. Furthermore, if $S\subset \mathbb{D}$ approaches $1$ nontangentially, then
\begin{equation}\label{limitequationin5.12} 
\lim_{\substack{z\rightarrow 1 \\ z\in S}}u(z\tau)=u(\tau).
\end{equation}
Consequently, 
\begin{equation}\label{thisequationforliminf}
\|u(\tau)\|_{\mathcal{M}}^{2}=\liminf_{\lambda\rightarrow \tau}\dfrac{1-\lvert \varphi(\lambda) \rvert^{2}}{1-\|\lambda\|^{2}_{\infty}}.
\end{equation}
\end{theorem}
\begin{proof}
Consider the sequence $\lambda^{(n)}=(1-1/n)\tau$ in $\mathbb{D}^{d}$ for $n\in\mathbb{N}$, so that $\lambda^{(n)}{\overset{nt}\rightarrow}\tau$. By Proposition \ref{prop5.2cara}, the sequence $(u(\lambda^{(n)}))$ is bounded in $\mathcal{M}$ by a constant $c$. Define the ball
$$\mathrm{Ball}(\mathcal{M}, c):=\{x\in\mathcal{M}:\|x\|_{\mathcal{M}}\leq c\}.$$
Thus
$$\{u(\lambda^{(n)}):n\in\mathbb{N}\}\subseteq \mathrm{Ball}(\mathcal{M}, c).$$
By \cite[Theorems 3.6.3 and  3.6.11]{BarrySimon}, any closed, bounded ball in a separable Hilbert space is not only compact for the weak topology on $\mathcal{M}$, but also metrizable for the weak topology on $\mathcal{M}$. By \cite[Theorem 2.3.6]{BarrySimon}, since $\mathrm{Ball}(\mathcal{M}, c)$ is metrizable and weakly compact, $\mathrm{Ball}(\mathcal{M}, c)$ is weakly sequentially compact. Recall that a metric space $X$ is said to be sequentially compact if every sequence in $X$ has a convergent subsequence.  Since $$\{u((1-1/n)\tau)\}\subseteq \mathrm{Ball}(\mathcal{M}, c),$$ there is a subsequence $(z_{n})$ of $(1-1/n)$ such that $z_{n}\tau\rightarrow\tau$ and $u(z_{n}\tau)$ converges weakly to some element $x\in\mathrm{Ball}(\mathcal{M}, c)$. Therefore there exists a sequence $z_{n}\rightarrow 1$ such that $u(z_{n}\tau)$ tends to $x$ weakly in $\mathcal{M}$. By Theorem \ref{realizationtheoremschuragler}, for all $\lambda \in \D^d$, 
$$ (1_{\mathcal{M}}-D\lambda_{P})u(\lambda)=\gamma.$$
We conclude that $(1_{\mathcal{M}}-D\tau_{P})x=\gamma$. Since $D\tau_{P}$ is a contraction, $\mathrm{Ker}(1_{\mathcal{M}}-D\tau_{P})\perp\mathrm{Ran}(1_{\mathcal{M}}-D\tau_{P})$. Hence $\gamma\perp\mathrm{Ker}(1_{\mathcal{M}}-D\tau_{P})$, and it follows from equation \eqref{gamma-D-u} that, for all $z\in\mathbb{D}$,
$$u(z\tau)=(1_{\mathcal{M}}-D\tau_{P})(1_{\mathcal{M}}-zD\tau_{P})^{-1}x,$$
and so $u(z\tau)\perp\mathrm{Ker}(1_{\mathcal{M}}-D\tau_{P})$. Hence $x\perp\mathrm{Ker}(1_{\mathcal{M}}-D\tau_{P})$. We have shown that there is a vector $x\in\mathcal{M}$ with the properties that $(1_{\mathcal{M}}-D\tau_{P})x=\gamma$ and $x\perp\mathrm{Ker}(1_{\mathcal{M}}-D\tau_{P})$. Since such a vector is unique, we deduce the first assertion of the theorem by taking $u(\tau)=x$. 

To see equation (\ref{limitequationin5.12}), suppose that $z_{n}\rightarrow 1$ and $u(z_{n}\tau)\rightarrow v$. Then $v\perp\mathrm{Ker}(1_{\mathcal{M}}-D\tau_{P})$ and $(1_{\mathcal{M}}-D\tau_{P})v=\gamma$. Hence $v=u(\tau)$.

For any $r\in(0,1)$, we have by  equation (\ref{whenlambdametmu}),
$$\dfrac{1-\lvert\varphi(r\tau)\rvert^{2}}{1-r^{2}}=\|u(r\tau)\|_{\mathcal{M}}^{2}.$$
As $r\rightarrow 1^{-}$, the right hand side of the above equation tends to $\|u(\tau)\|^{2}$, and, by equation (\ref{limitequationin5.12}), the left hand side tends to the lim inf in equation (\ref{thisequationforliminf}), by Corollary \ref{coro4.14cara} and Proposition \ref{prop5.8cara}. This establishes equation (\ref{thisequationforliminf}).\end{proof}

\begin{proposition}\label{prop5.17cara}
Let $(\mathcal{M},u)$ be a model of $\varphi\in\mathcal{SA}_{d}$ and let $(a,\beta,\gamma,D)$ be a realization of $(\mathcal{M},u)$. \\
{\rm (i)} The following conditions are equivalent for a point $ \ \tau \in \partial \mathbb{D}^d \setminus \mathbb{T}^d$, 
\begin{enumerate}
\item  $\tau$ is a $B$-point for $\ph$;
\item $\gamma\in\mathrm{Ran}(1_{\mathcal{M}}-D\tau_{P})$ and
$ \lim_{\la \nt \tau} |\varphi(\lambda)| =1$. 
\end{enumerate}
{\rm (ii)} The following conditions are equivalent for a point  $ \ \tau \in  \mathbb{T}^d$:
\begin{enumerate}
\item  $\tau$ is a $B$-point for $\ph$;
\item $\gamma\in\mathrm{Ran}(1_{\mathcal{M}}-D\tau_{P})$.
\end{enumerate}
\end{proposition}
\begin{proof}
 In both cases (1) implies (2), because, by Theorem \ref{theorem5.12cara},  if $\tau$ is a $B$-point then there exists
 a $u({\tau}) \in \mathcal{M}$ such that $(1_{\mathcal{M}}-D\tau_{P})u(\tau)=\gamma$, and hence $\gamma\in\mathrm{Ran}(1_{\mathcal{M}}-D\tau_{P})$.
Note that,
by Lemma \ref{limit-mod}, if $\tau$ is a  $B$-point for $\varphi$
 then $ \lim_{\la \nt \tau} |\varphi(\lambda)| =1$.

To prove that (2) implies (1), suppose that $(1_{\mathcal{M}}-D\tau_{P})x=\gamma$ for some $x\in\mathcal{M}$. 
 Then, by equation \eqref{gamma-D-u}, for $\lambda\in\mathbb{D}^{d}$, 
\begin{align*}
u(\lambda)&=(1_{\mathcal{M}}-D\lambda_{P})^{-1}\gamma=(1_{\mathcal{M}}-D\lambda_{P})^{-1}(1_{\mathcal{M}}-D\tau_{P})x\\
    &=(1_{\mathcal{M}}-D\lambda_{P})^{-1}(1_{\mathcal{M}}-D\lambda_{P}+D(\lambda_{P}-\tau_{P}))x\\
    &=x+(1_{\mathcal{M}}-D\lambda_{P})^{-1}D(\lambda_{P}-\tau_{P})x.
\end{align*}
Since $D$ is a contraction, 
$$\|(1_{\mathcal{M}}-D\lambda_{P})^{-1}D\|=\bigg\| \sum_{j=1}^{\infty}(D\lambda_{P})^{j}D\bigg\|\leq \sum_{j=1}^{\infty}\|\lambda\|_{\infty}^{j}=\dfrac{1}{1-\|\lambda\|_{\infty}},$$
and so
\begin{align*}
\|u(\lambda)\|_{\mathcal{M}}&\leq \|x\|_{\mathcal{M}}+\|(1_{\mathcal{M}}-D\lambda_{P})^{-1}\|_{\mathcal{B(M)}}\cdot \|(\lambda_{P}-\tau_{P})x\|_{\mathcal{M}}\\
&\leq \|x\|_{\mathcal{M}}+\dfrac{\|\lambda-\tau\|_{\infty}}{1-\|\lambda\|_{\infty}}\|x\|_{\mathcal{M}}.
\end{align*}
If   $\lambda\in S{\overset{nt}\rightarrow}\tau$,
 there exists a $c>0$ such that
$$\|\lambda - \tau\|_{\infty}\leq c(1-\|\lambda\|_{\infty}) \; \; \text{for all } \; \lambda \in S.$$
Hence, for all $\lambda\in S{\overset{nt}\rightarrow}\tau$, 
$$\|u(\lambda)\|_{\mathcal{M}}\leq (1+c)\|x\|_{\mathcal{M}}.$$
Thus $u(\lambda)$ is bounded on $S$.
In the case when  $ \ \tau \in \partial \D^d \setminus\T^d$,   by assumption, 
 $ \lim_{\la \nt \tau} |\varphi(\lambda)| =1$.  
By Corollary \ref{coro5.7cara}, $\tau$ is a $B$-point of $\varphi$.\end{proof}

\section{Interactions of  carapoints with  models of Schur-Agler class functions} \label{interaction}

Our next task is to show that if a function  $\vp\in\sad$ has a singularity at a $B$-point $\tau$, then we can construct a generalized model  $(\mathcal H, u, I)$
of $\ph$ in which the singularity of $\ph$ is encoded in an $I(\cdot)$, in such a way that the model has a $C$-point at $\tau$.  The device that leads to this conclusion is to write vectors in and operators on $\M$ in terms of the orthogonal decomposition $\M=\N\oplus\N^\perp$ where $\N= \ker(1-D\tau_P)$ and $D$ comes from a realization of $(\M,u)$.  The following observation is straightforward.

Let $\mathcal{M}$ be a separable Hilbert space with orthogonal decomposition 
$$\mathcal{M}=\mathcal{M}_{1}\oplus\mathcal{M}_{2}\hdots\oplus \mathcal{M}_{d}.$$
For $\tau\in\mathbb{T}^{d}$, let us define the operator
$$\tau_{P} = \tau_{1}P_{1}+\hdots+\tau_{d}P_{d},$$
where $P_{j}:\mathcal{M}\rightarrow\mathcal{M}$ is the orthogonal projection from $\mathcal{M}$ onto $\mathcal{M}_{j}$ for $j=1,\hdots,d$. It is easy to see that $\tau_{P}$ is a unitary operator, 
\begin{align*}
    \tau_{P}^{*}\tau_{P} &= \overline{\tau}_{P}\tau_{P} = \lvert \tau_{1} \rvert^{2}P_{1}+\hdots+\lvert \tau_{d}\rvert^{2}P_{d} = P_{1}+\hdots+P_{d} = 1_{\mathcal{M}}\\
    \tau_{P}\tau^{*}_{P} &= \tau_{P}\overline{\tau}_{P} = \lvert \tau_{1} \rvert^{2}P_{1}+\hdots+\lvert \tau_{d}\rvert^{2}P_{d} = P_{1}+\hdots+P_{d} = 1_{\mathcal{M}}.
\end{align*} 
\begin{definition}\label{lambdaincdastupleofoperators}
Let $\mathcal{M}$ be a complex Hilbert space and let $\lambda \in \mathbb{C}^{d}$. Then if $X = (X_{1}, \hdots, X_{d})$ is a $d$-tuple of operators on $\mathcal{M}$, we define the operator
$$ \lambda_{X} = \lambda_{1}X_{1}+\hdots+\lambda_{d}X_{d}.$$
\end{definition}

Realizations provide an effective tool for the study of boundary behavior.  Here is a preliminary observation.
\begin{lemma}{\normalfont{\cite[Lemma 3.3]{ATY2}}}\label{gammainranperp}
Suppose that $\varphi\in \mathcal{SA}_{d}$ has a model $(\mathcal{M},u)$ with realization $(a,\beta, \gamma, D)$. Let $\tau \in \mathbb{T}^{d}$ be a carapoint of $\varphi$, let $\mathcal{N} = \mathrm{Ker}(1_{\mathcal{M}}-D\tau_{P})$, $\mathcal{M}=\mathcal{N}\oplus \mathcal{N}^{\perp}$. Then 
$$\gamma \in \mathrm{Ran}(1_{\mathcal{M}}-D\tau_{P}) \subset \mathcal{N}^{\perp}~~\text{and}~~ \overline{\tau_{P}}\beta \in \mathcal{N}^{\perp}.$$
\end{lemma}
\begin{proof}
First we show that $\overline{\tau_{P}}\beta\in\mathcal{N}^{\perp}$. Let 
$$L=\begin{bmatrix}
    a       & 1\otimes \beta \\
    \gamma \otimes 1       & D 
\end{bmatrix}.$$
Choose any $x\in\mathcal{N}$. Then $x=D\tau_{P} x$ and so
$$L\begin{bmatrix}
    0\\
    \tau_{P} x
\end{bmatrix}=  \begin{bmatrix}
                a & 1\otimes\beta \\
                \gamma\otimes 1 & D \\
                \end{bmatrix}\begin{bmatrix}
    0\\
    \tau_{P} x
\end{bmatrix}=\begin{bmatrix}
    \langle \tau_{P} x, \beta\rangle_{\mathcal{N}}\\
    D\tau_{P} x
\end{bmatrix}=\begin{bmatrix}
    \langle x, \overline{\tau_{P}}\beta\rangle_{\mathcal{N}}\\
    x
\end{bmatrix}.
$$
Since $L$ is a contraction and $\tau_{P}$ is an isometry, 
$$\Bigg\|\begin{bmatrix}
    \langle x, \overline{\tau_{P}}\beta\rangle_{\mathcal{N}}\\
    x
\end{bmatrix}\Bigg\|=\Bigg\|L\begin{bmatrix}
    0\\
    \tau_{P} x
\end{bmatrix}\Bigg\|\leq\Bigg\|\begin{bmatrix}
    0\\
    \tau_{P} x
\end{bmatrix}\Bigg\|=\|\tau_{P} x\|_{\mathcal{N}}=\|x\|_{\mathcal{N}},$$
and so $\langle x,\overline{\tau_{P}}\beta\rangle_{\mathcal{N}}=0$. Since $x\in\mathcal{N}$ is arbitrary, $\overline{\tau_{P}}\beta\in\mathcal{N}^{\perp}$.

 Proposition \ref{prop5.17cara} asserts that $\tau$ is a carapoint for $\varphi$ if and only if $\gamma\in\mathrm{Ran}(1_{\mathcal{M}}-D\tau_{P})$. Since $D\tau_{P}$ is a contraction on $\mathcal{M}$, then by 
 \cite[Proposition 3.1]{FoiasNagy}, for an eigenvector $x$ of $D\tau_{P}$ with eigenvalue $\lambda\in\mathbb{T}$, then $x$ is also an eigenvector of $(D\tau_{P})^{*}$ with eigenvalue $\overline{\lambda}\in\mathbb{T}$. Hence $\mathrm{Ker}(1_{\mathcal{M}}-D\tau_{P})=\mathrm{Ker}(1_{\mathcal{M}}-\overline{\tau_{P}}D^{*})$. Furthermore, by \cite[Lemma 6.11]{RynneYoungson},
$$\mathrm{Ker}(1_{\mathcal{M}}-\overline{\tau_{P}}D^{*}) = \mathrm{Ran}(1_{\mathcal{M}}-D\tau_{P})^{\perp}.$$
Hence
$$\Big(\mathrm{Ran}(1_{\mathcal{M}}-D\tau_{P})^{\perp}\Big)^{\perp}=\overline{\mathrm{Ran}(1_{\mathcal{M}}-D\tau_{P})} = \mathrm{Ker}(1_{\mathcal{M}}-\overline{\tau_{P}}D^{*})^{\perp}.$$
Therefore
$$\gamma\in\mathrm{Ran}(1_{\mathcal{M}}-D\tau_{P}) \subset \overline{\mathrm{Ran}(1_{\mathcal{M}}-D\tau_{P})}=\mathrm{Ker}(1_{\mathcal{M}}-\overline{\tau_{P}}D^{*})^{\perp}=\mathrm{Ker}(1_{\mathcal{M}}-D\tau_{P})^{\perp}=\mathcal{N}^{\perp}.$$\end{proof}

\begin{lemma}\label{projectiononcdlemma} 
Let $\mathcal{M} = \mathcal{M}_{1}\oplus \hdots \oplus \mathcal{M}_{d}$ be a complex Hilbert space, $P_{j}:\mathcal{M}\rightarrow \mathcal{M}$ be the orthogonal projection onto $\mathcal{M}_{j}$ for $j=1,\hdots,d$ and let $\mathcal{N}$ be a closed subspace of $\mathcal{M}$. Then with respect to the decomposition $\mathcal{M} = \mathcal{N} \oplus \mathcal{N}^{\perp}$, the operator $P_{j}$ has the operator matrix form
\begin{equation}\label{projectiononmd}
P_{j} =\begin{bmatrix}
    X_{j}       & B_{j} \\
    B^{*}_{j}      & Y_{j} 
\end{bmatrix},
\end{equation}
for some $X_{j}\in \mathcal{B}(\mathcal{N})$, $Y_{j} \in \mathcal{B}(\mathcal{N}^{\perp})$ and $B_{j} \in \mathcal{B}(\mathcal{N}^{\perp},\mathcal{N})$ for all integers such that $1\leq j \leq d$. Moreover, for $j=1,\hdots, d$,
\begin{enumerate}
  \item  $0 \leq X_{j} \leq 1_{\mathcal{N}}$, $ 0 \leq Y_{j} \leq 1_{\mathcal{N}^{\perp}}$;
  \item $\sum_{j=1}^{d} X_{j} = 1_{\mathcal{N}},~\sum_{j=1}^{d} B_{j} = 0, ~ \sum_{j=1}^{d} Y_{j} = 1_{\mathcal{N}^{\perp}}$;
  \item $B_{i}Y_{j} =  \delta_{ij}B_{j}-X_{i}B_{j}$, $B_{i}^{*}X_{j} =  \delta_{ij}B_{j}^{*}-Y_{i}B_{j}^{*} $;
  \item $B_{i}B^{*}_{j} =  \delta_{ij}X_{j}-X_{i}X_{j}, ~B^{*}_{i}B_{j} =  \delta_{ij}Y_{j}-Y_{i}Y_{j}$,
\end{enumerate}
where $\delta_{ij}$ denotes the Kronecker delta.
\end{lemma}
\begin{proof}
(1) Let $P_{j}:\mathcal{M}\rightarrow\mathcal{M}$ be the orthogonal projection onto $\mathcal{M}_{j}$ for $j=1,\hdots,d$. We may represent the operator $P_{j}$ on $\mathcal{M}$ by an operator matrix with respect to the orthogonal decomposition $\mathcal{M}=\mathcal{N}\oplus \mathcal{N}^{\perp}$ as follows. Let
\begin{equation}\nonumber
P_{j} =\begin{bmatrix}
    X_{j}       & B_{j} \\
    C_{j}      & Y_{j} 
\end{bmatrix},
\end{equation}
for some operators $X_{j},B_{j},C_{j},Y_{j}$. Since $P_{j}$ is an orthogonal projection, $P_{j} = P_{j}^{*}$, and so $X_{j}=X^{*}_{j}, B^{*}_{j}=C_{j}$ and $Y_{j} = Y^{*}_{j}$. Define $P_{\mathcal{N}}:\mathcal{M}\rightarrow \mathcal{M}$ to be the orthogonal projection onto $\mathcal{N}$. By equation (\ref{projectiononmd}), for $j=1,\hdots,d$,
$$X_{j} = P_{\mathcal{N}}P_{j}P_{\mathcal{N}}^{*}.$$
Note $0\leq P_{j} \leq 1_{\mathcal{M}}$. Hence
$$0\leq P_{\mathcal{N}}P_{j}P_{\mathcal{N}}^{*}=X_{j}\leq P_{\mathcal{N}}P_{\mathcal{N}}^{*} = 1_{\mathcal{N}},$$
and similarly for $Y_{j}$.

For all $x \in \mathcal{M}$, $\| P_{j}x\| \leq \| x \|$. Choose a $u\in\mathcal{N}$, then
\begin{equation}\nonumber
\| P_{j} u \|^{2} = \Bigg\| \begin{bmatrix}
    X_{j}u \\
    B^{*}_{j}u \\ 
\end{bmatrix} \Bigg\|^{2} = \big\| X_{j}u \big\|^{2}_{\mathcal{N}} + \big\| B^{*}_{j}u \big\|^{2}_{\mathcal{N}^{\perp}} \leq \| u \|^{2}_{\mathcal{N}},
\end{equation}
and therefore $\| X_{j}\| \leq 1$ for $1\leq j \leq d$. We likewise produce a similar result for $Y_{j}$ with a choice of $v\in\mathcal{N}^{\perp}$.

 (2) Since $P_{j}$ is the orthogonal projection onto $\mathcal{M}_{j}$, $j=1,\hdots,d$, we have, for all $x\in\mathcal{M}$,
$$P_{1}x + \hdots + P_{d}x = x.$$	
Therefore
$$ \sum_{j=1}^{d} P_{j} = 1_{\mathcal{M}},$$
and thus
$$ 1_{\mathcal{M}} = \begin{bmatrix}
	1_{\mathcal{N}} & 0 \\
	0 & 1_{\mathcal{N}^{\perp}}
\end{bmatrix} = \sum_{j=1}^{d} \begin{bmatrix}
X_{j} & B_{j} \\
B^{*}_{j} & Y_{j}
\end{bmatrix} = \begin{bmatrix}
\sum_{j=1}^{d} X_{j} & \sum_{j=1}^{d} B_{j} \\
\sum_{j=1}^{d} B^{*}_{j} & \sum_{j=1}^{d} Y_{j}
\end{bmatrix}. $$
Hence
$$ \sum_{j=1}^{d} X_{j} = 1_{\mathcal{N}}, ~~\sum_{j=1}^{d} B_{j}  = 0 ~~\text{and}~~\sum_{j=1}^{d} Y_{j}  = 1_{\mathcal{N}^{\perp}}.$$

(3) Again by properties of projections, $P^{2}_{j} = P_{j}$, and so
\begin{align}  \label{system1projections}
\left\{
\begin{aligned}
	X^{2}_{j} +B_{j}B^{*}_{j} & = X_{j} \nonumber  \\
	X_{j}B_{j}+B_{j}Y_{j} &= B_{j} \\
	B^{*}_{j}X_{j} + Y_{j}B^{*}_{j} &= B^{*}_{j} \nonumber\\
	B^{*}_{j}B_{j}+Y_{j}^{2} &= Y_{j},\nonumber
\end{aligned}\right.\\
\end{align}
that is,
\begin{align}
\left\{\begin{aligned}
B_{j}B^{*}_{j} & = X_{j} - X_{j}^{2} \nonumber\\
	B_{j}Y_{j} &= B_{j} - X_{j}B_{j} \nonumber \\   
	Y_{j}B^{*}_{j} & = B^{*}_{j}-B^{*}_{j}X_{j}\nonumber \\
	B^{*}_{j}B_{j} &= Y_{j}-Y^{2}_{j}.  \nonumber
\end{aligned}\right.  
\end{align}
Whereas, for $i\neq j$, $P_{i}P_{j}=0$ results in 
\begin{align}
\left\{
\begin{aligned}
	X_{i}X_{j} +B_{i}B^{*}_{j} & = 0 \nonumber\\
	X_{i}B_{j}+B_{i}Y_{j} &= 0  \nonumber\\
     	B^{*}_{i}X_{j} + Y_{i}B^{*}_{j} &= 0 \nonumber\\
	B^{*}_{i}B_{j}+Y_{i}Y_{j} &= 0,\nonumber
\end{aligned}\right.
\end{align}
that is, 
\begin{align}
\left\{
\begin{aligned}
B_{i}B^{*}_{j} & = - X_{i}X_{j}\nonumber\\
	B_{i}Y_{j} &= - X_{i}B_{j} \label{system2projections}\\
	Y_{i}B^{*}_{j} & =- B^{*}_{i}X_{j}\nonumber\\
	B^{*}_{i}B_{j} &= -Y_{i}Y_{j}.\nonumber 
\end{aligned}\right. \\
\end{align}
Combing the two systems of equations (\ref{system1projections}) and (\ref{system2projections}), we see that
\begin{align}
B_{i}B^{*}_{j} =  \delta_{ij}X_{j}-X_{i}X_{j} \label{propertiesoftheprojectionoperators1} \\
B_{i}Y_{j} =  \delta_{ij}B_{j}-X_{i}B_{j} \label{propertiesoftheprojectionoperators2} \\
B_{i}^{*}X_{j} =  \delta_{ij}B_{j}^{*}-Y_{i}B_{j}^{*} \label{propertiesoftheprojectionoperators3}\\
B^{*}_{i}B_{j} = \delta_{ij}Y_{j} -Y_{i}Y_{j}, \label{propertiesoftheprojectionoperators4}
\end{align}
where $ \delta_{ij}$ denotes the Kronecker delta. \end{proof}

These next results contain the essential tools needed to prove one of our main results, Theorem \ref{generalizedmodeltheorem}.

\begin{proposition}\label{mulambdaidentitiesforgeneralizedmodels}
Let $\mathcal{M} = \mathcal{M}_{1}\oplus \hdots \oplus \mathcal{M}_{d}$ be a complex Hilbert space and let $\mathcal{N}$ be a closed subspace of $\mathcal{M}$.  Let $X=(X_{1},\hdots, X_{d}), B=(B_{1},\hdots, B_{d})$ and $ Y=(Y_{1},\hdots, Y_{d})$ be  $d$-tuples of operators as in Lemma \ref{projectiononcdlemma}, where $X_{j}\in \mathcal{B}(\mathcal{N})$, $Y_{j} \in \mathcal{B}(\mathcal{N}^{\perp})$ and $B_{j} \in \mathcal{B}(\mathcal{N}^{\perp},\mathcal{N})$ for all integers such that $1\leq j \leq d$, and the following conditions are satisfied
\begin{enumerate}
  \item  $0 \leq X_{j} \leq 1_{\mathcal{N}}$, $ 0 \leq Y_{j} \leq 1_{\mathcal{N}^{\perp}}$;
  \item $\sum_{j=1}^{d} X_{j} = 1_{\mathcal{N}},~\sum_{j=1}^{d} B_{j} = 0, ~ \sum_{j=1}^{d} Y_{j} = 1_{\mathcal{N}^{\perp}}$;
  \item $B_{i}Y_{j} =  \delta_{ij}B_{j}-X_{i}B_{j}$, $B_{i}^{*}X_{j} =  \delta_{ij}B_{j}^{*}-Y_{i}B_{j}^{*} $;
  \item $B_{i}B^{*}_{j} =  \delta_{ij}X_{j}-X_{i}X_{j}, ~B^{*}_{i}B_{j} =  \delta_{ij}Y_{j}-Y_{i}Y_{j}$.
\end{enumerate}
Then for all $\lambda, \mu \in \mathbb{C}^{d}$,
\begin{align*}
\mu_{B}\lambda_{B^{*}} & = (\mu\lambda)_{X} - \mu_{X}\lambda_{X},\\
\mu_{B^{*}}\lambda_{X} & = (\mu\lambda)_{B^{*}} - \mu_{Y}\lambda_{B^{*}},\\
\mu_{B}\lambda_{Y} & = (\mu\lambda)_{B} - \mu_{X}\lambda_{B},\\
\mu_{B^{*}}\lambda_{B}& = (\mu\lambda)_{Y} - \mu_{Y}\lambda_{Y}.
\end{align*}
\end{proposition}
\begin{proof}
For all $\lambda, \mu \in\mathbb{C}^{d}$, 
$$\mu_{B}\lambda_{B^{*}} = \mu_{B}(\lambda_{1}B^{*}_{1} + \hdots + \lambda_{d}B^{*}_{d})  = \lambda_{1}\mu_{B}B^{*}_{1}+\hdots+\lambda_{d}\mu_{B}B_{d}^{*}.$$
Note that, by equation (\ref{propertiesoftheprojectionoperators1}), 
\begin{equation}\label{mubequation}
\mu_{B}B^{*}_{j} = (\mu_{1}B_{1}+\hdots + \mu_{d}B_{d})B_{j}^{*} = \mu_{1}( \delta_{1j}X_{j}-X_{1}X_{j})+\hdots+\mu_{d}( \delta_{dj}X_{j}-X_{d}X_{j}).
\end{equation}
Therefore
\begin{align*}
\mu_{B}\lambda_{B^{*}}
&= \mu_{B}(\lambda_{1}B_{1}^{*}+\hdots+\lambda_{d}B^{*}_{d})  \\	
&= \lambda_{1}\big(\mu_{1}(X_{1}-X^{2}_{1})- \hdots - \mu_{d}X_{d}X_{1})\\
&\hspace{4.5cm}+\hdots+\lambda_{d}\big(-\mu_{1}X_{1}X_{d}+\hdots+\mu_{d}(X_{d}-X^{2}_{d})\big) \\
& = \mu_{1}\lambda_{1}X_{1}+\hdots +\mu_{d}\lambda_{d}X_{d}-(\mu_{1}\lambda_{1}X^{2}_{1}+\hdots+\mu_{d}\lambda_{1}X_{d}X_{1})\\
&\hspace{6cm}-\hdots-(\mu_{1}\lambda_{d}X_{1}X_{d}+\hdots+\mu_{d}\lambda_{d}X_{d}^{2}) \\
& = (\mu\lambda)_{X} - \big(\mu_{X}\lambda_{1}X_{1}+\hdots+\mu_{X}\lambda_{d}X_{d}\big) \\
&= (\mu\lambda)_{X}-\mu_{X}\lambda_{X}.
\end{align*}
Similar calculations yield all the other identities.\end{proof}

\begin{theorem}\label{invertibleoperatorsforlambda}
Let $T=(T_{1},\hdots, T_{d})$ be a d-tuple of operators on a Hilbert space $\mathcal{H}$ such that $\sum_{j=1}^{d}T_{j}=1_{\mathcal{H}}$ and $0\leq T_{j}\leq 1$ for $j=1,\hdots, d$. For all $\lambda=(\lambda_{1},\hdots,\lambda_{d})\in\mathbb{C}^{d}$ such that $\mathrm{Re}(\lambda_{j})<1$ for $j=1,\hdots,d$, the operators 
$$(\mathds{1}-\lambda)_{T}~\text{and}~\bigg(\dfrac{\mathds{1}}{\mathds{1}-\lambda}\bigg)_{T}$$
on $\mathcal{B(H)}$ are invertible and
\begin{align}
\Big\| (\mathds{1}-\lambda)_{T}^{-1} \Big\|_{\mathcal{B(H)}} & \leq \dfrac{1}{1-\underset{j=1,\hdots, d}{\max}\mathrm{Re}(\lambda_{j})},\label{boundlabelforlambdax} \\
\Bigg\| \bigg(\dfrac{\mathds{1}}{\mathds{1}-\lambda}\bigg)_{T}^{-1} \Bigg\|_{\mathcal{B(H)}} &\leq ~ \max_{j=1,\hdots,d}\dfrac{\lvert 1- \lambda_{j} \rvert^{2}}{1-\mathrm{Re}(\lambda_{j})}.\label{boundlabelforlambdax2}
\end{align}
\end{theorem}
\begin{proof}
Let $u\in\mathcal{H}$ be such that $\| u \|_{\mathcal{H}} =1$. By the Cauchy-Schwarz inequality \index{Cauchy-Schwarz inequality}
$$ \lvert \langle x,u \rangle_{\mathcal{H}} \rvert \leq \| x \|_{\mathcal{H}} \| u \|_{\mathcal{H}}~~\text{for all}~x\in\mathcal{H}.$$
By assumption, $\sum_{j=1}^{d}T_{j}=1_{\mathcal{H}}$ and $\mathrm{Re}(\lambda_{j})<1$ for $j=1,\hdots,d$ and so $\mathrm{Re}(1-\lambda_{j})>0$ for $j=1,\hdots,d$. For $u\in\mathcal{H}$ with $\|u\|_{\mathcal{H}}=1$, we have
\begin{align}
\| (\mathds{1}-\lambda)_{T}u\|_{\mathcal{H}} 	& = \| (\mathds{1}-\lambda)_{T}u\|_{\mathcal{H}} \| u \|_{\mathcal{H}} \nonumber\\
				& \geq \big\lvert \langle (\mathds{1}-\lambda)_{T}u,u \rangle_{\mathcal{H}} \big\rvert 	\nonumber\\
				& \geq \mathrm{Re}\Big(\langle(\mathds{1}-\lambda)_{T}u,u\rangle_{\mathcal{H}}\Big) \nonumber\\
				& = \mathrm{Re}\Bigg(\Big\langle \sum_{j=1}^{d}(1-\lambda_{j})T_{j}u,u\Big\rangle_{\mathcal{H}}\Bigg) \nonumber\\
				& = \sum_{j=1}^{d} \mathrm{Re}\Bigg((1-\lambda_{j})\Big\langle T_{j}u,u\Big\rangle_{\mathcal{H}}\Bigg). \label{eq5star}
\end{align}
By assumption, $0\leq T_{j} \leq 1_{\mathcal{H}}$, that is, $0\leq\langle T_{j}u,u\rangle_{\mathcal{H}}\leq 1$ for all $u\in\mathcal{H}$, and so
$$\mathrm{Re}\Big((1-\lambda_{j})\big\langle T_{j}u,u\big\rangle_{\mathcal{H}}\Big) = \mathrm{Re}(1-\lambda_{j})\big\langle T_{j}u,u\big\rangle_{\mathcal{H}}.$$
Therefore, for $u\in\mathcal{H}$ with $\|u\|_{\mathcal{H}}=1$, equation (\ref{eq5star}) implies
\begin{align}
\| (\mathds{1}-\lambda)_{T}u\|_{\mathcal{H}} & \geq \sum_{j=1}^{d} \mathrm{Re}(1-\lambda_{j})\Big\langle T_{j}u,u\Big\rangle_{\mathcal{H}} \label{anequationwithrealandinner}\\
				& \geq \underset{j=1,\hdots, d}{\min}\mathrm{Re}(1-\lambda_{j}) \sum_{j=1}^{d}\langle T_{j}u,u\rangle_{\mathcal{H}}\nonumber.
\end{align}
Again, $\sum_{j=1}^{d}T_{j}=1_{\mathcal{H}}$, so that, for $u\in\mathcal{H}$ with $\|u\|_{\mathcal{H}}=1$
\begin{equation}\label{anequationwithrealandinner2}
    \sum_{j=1}^{d}\big\langle T_{j}u,u\big\rangle_{\mathcal{H}} = \bigg\langle \sum_{j=1}^{d}T_{j}u,u \bigg\rangle_{\mathcal{H}} = \langle u,u \rangle_{\mathcal{H}} = \| u \|_{\mathcal{H}}^{2}=1.
\end{equation}
Now inequality (\ref{anequationwithrealandinner}) and equation (\ref{anequationwithrealandinner2}) imply 
\begin{equation}\label{firstlineininvertibilityargument}
\| (\mathds{1}-\lambda)_{T}u\| \geq \underset{j=1,\hdots, d}{\min} \mathrm{Re}(1-\lambda_{j}),
\end{equation}
For all $v\in\mathcal{H}$, let $u = v/\|v\|_{\mathcal{H}}$. By equation (\ref{firstlineininvertibilityargument}),
\begin{equation}\label{modifiedfirstlineininvertibilityargument}
\| (\mathds{1}-\lambda)_{T} v \|_{\mathcal{H}}\geq \underset{j=1,\hdots, d}{\min}\mathrm{Re}(1-\lambda_{j}) \| v\|_{\mathcal{H}}.
\end{equation} 

Hence, by the open mapping theorem \cite[Theorem 4.43]{RynneYoungson},  $(\mathds{1}-\lambda)_{T}$ has a bounded left inverse. Furthermore, since $T_{j}=T_{j}^{*}$ for $j=1,\hdots,d$, we have $(\mathds{1}-\lambda)_{T}^{*} = (\mathds{1}-\overline{\lambda})_{T}$. Therefore a similar calculation to the above shows there is a bounded left inverse for $(\mathds{1}-\lambda)_{T}^{*}$. Hence $(\mathds{1}-\lambda)_{T}$ is invertible. To prove inequality (\ref{boundlabelforlambdax}), note that, by equation (\ref{modifiedfirstlineininvertibilityargument}), for all $v \in\mathcal{H}$, there exists a unique $w \in\mathcal{H}$ such that $v=(\mathds{1}-\lambda)_{T}^{-1}w$ and
\begin{align}
\| (\mathds{1}-\lambda)_{T} v \|_{\mathcal{H}} 	& = \| (\mathds{1}-\lambda)_{T} (\mathds{1}-\lambda)_{T}^{-1}w \|_{\mathcal{H}} \nonumber \\
						& = \| w \|_{\mathcal{H}} \nonumber \\
						&\geq \underset{j=1,\hdots, d}{\min}\mathrm{Re}(1-\lambda_{j}) \| (\mathds{1}-\lambda)_{T}^{-1}w\|_{\mathcal{H}} \label{invertibleoperatorboundsforprojection}. \\ \nonumber
\end{align}
Since by assumption, $\mathrm{Re}(\lambda_{j})<1$ for $j=1,\hdots,d$, we have $\mathrm{Re}(1-\lambda_{j}) >0$ for $j=1,\hdots,d$. Rearranging equation (\ref{invertibleoperatorboundsforprojection}), we find that 
$$ \| (\mathds{1}-\lambda)_{T}^{-1}w \| \leq \dfrac{\|w\|_{\mathcal{H}}}{\underset{j=1,\hdots,d}{\min}\mathrm{Re}(1-\lambda_{j})} = \dfrac{\|w\|_{\mathcal{H}}}{1-\underset{j=1,\hdots,d}{\max}\mathrm{Re}(\lambda_{j})}~~\text{for all}~~w\in\mathcal{H}.$$
Thus inequality (\ref{boundlabelforlambdax}) is proved.

 For the operator $\big(\frac{\mathds{1}}{\mathds{1}-\lambda}\big)_{T}$, we follow a similar argument to the one above. Let $u\in\mathcal{H}$ be such that $\| u \|_{\mathcal{H}} =1$. By the Cauchy-Schwarz inequality,
\begin{align*}
\bigg\| \bigg(\dfrac{\mathds{1}}{\mathds{1}-\lambda}\bigg)_{T}u \bigg\|_{\mathcal{H}} 	& = \bigg\| \bigg(\dfrac{\mathds{1}}{\mathds{1}-\lambda}\bigg)_{T}u \bigg\|_{\mathcal{H}} \| u \|_{\mathcal{H}} \\
				& \geq \bigg\lvert \bigg\langle \bigg(\dfrac{\mathds{1}}{\mathds{1}-\lambda}\bigg)_{T}u,u \bigg\rangle_{\mathcal{H}} \bigg\rvert 	\\
				& \geq \mathrm{Re}\Bigg(\bigg\langle\bigg(\dfrac{\mathds{1}}{\mathds{1}-\lambda}\bigg)_{T}u,u\bigg\rangle_{\mathcal{H}}\Bigg) \\
				& = \mathrm{Re}\Bigg(\bigg\langle \sum_{j=1}^{d}\bigg(\dfrac{1}{1-\lambda_{j}}\bigg)T_{j}u,u\bigg\rangle_{\mathcal{H}}\Bigg)\\
				& = \sum_{j=1}^{d}\mathrm{Re}\Bigg(\bigg(\dfrac{1}{1-\lambda_{j}}\bigg)\Big\langle T_{j}u,u\Big\rangle_{\mathcal{H}}\Bigg).\\
\end{align*}
By assumption, $0\leq T_{j} \leq 1_{\mathcal{H}}$, that is, $0\leq\langle T_{j}u,u\rangle_{\mathcal{H}}\leq 1$ for all $u\in\mathcal{H}$. Consequently,
$$\mathrm{Re}\Bigg(\bigg(\dfrac{1}{1-\lambda_{j}}\bigg)\Big\langle T_{j}u,u\Big\rangle_{\mathcal{H}^{\perp
}}\Bigg) = \mathrm{Re}\bigg(\dfrac{1}{1-\lambda_{j}}\bigg)\Big\langle T_{j}u,u\Big\rangle_{\mathcal{H}^{\perp
}}.$$
Therefore
\begin{align*}
\bigg\| \bigg(\dfrac{\mathds{1}}{\mathds{1}-\lambda}\bigg)_{T}u \bigg\|_{\mathcal{H}} & \geq \sum_{j=1}^{d} \mathrm{Re}\bigg(\dfrac{1}{1-\lambda_{j}}\bigg)  \Big\langle T_{j}u,u\Big\rangle_{\mathcal{H}} \\
				& \geq \underset{j=1,\hdots, d}{\min}\mathrm{Re}\bigg(\dfrac{1}{1-\lambda_{j}}\bigg) \sum_{j=1}^{d}\Big\langle T_{j}u,u\Big\rangle_{\mathcal{H}}.
\end{align*}
Again, $\sum_{j=1}^{d}T_{j}=1_{\mathcal{H}}$, and so
$$ \sum_{j=1}^{d}\big\langle T_{j}u,u\big\rangle_{\mathcal{H}} = \bigg\langle \sum_{j=1}^{d}T_{j}u,u \bigg\rangle_{\mathcal{H}} = \langle u,u \rangle_{\mathcal{H}} = \| u \|_{\mathcal{H}}^{2}=1.$$
Now
\begin{equation}\label{firstlineininvertibilityargumentxo}
\bigg\| \bigg(\dfrac{\mathds{1}}{\mathds{1}-\lambda}\bigg)_{T}u \bigg\|_{\mathcal{H}} \geq \underset{j=1,\hdots, d}{\min} \mathrm{Re}\bigg(\dfrac{1}{1-\lambda_{j}}\bigg),
\end{equation}
where by our assumption, $\mathrm{Re}(\lambda_{j})<1$ and so $\mathrm{Re}\bigg(\dfrac{1}{1-\lambda_{j}}\bigg) >0$. For $v\in\mathcal{H}$, let $u = v/\|v\|_{\mathcal{H}}$. By equation (\ref{firstlineininvertibilityargumentxo}),
\begin{equation}\label{modifiedfirstlineininvertibilityargument2}
\bigg\| \bigg(\dfrac{\mathds{1}}{\mathds{1}-\lambda}\bigg)_{T}v \bigg\|_{\mathcal{H}}\geq \underset{j=1,\hdots, d}{\min}\mathrm{Re}\bigg(\dfrac{1}{1-\lambda_{j}}\bigg) \| v\|_{\mathcal{H}}.
\end{equation} 
 By the open mapping theorem \cite[Theorem 4.43]{RynneYoungson}, 
 $\big(\frac{\mathds{1}}{\mathds{1}-\lambda}\big)_{T}$ has a bounded left inverse. Furthermore, since $T_{j}=T_{j}^{*}$ for $j=1,\hdots,d$, we have $\big(\frac{\mathds{1}}{\mathds{1}-\lambda}\big)_{T}^{*} = \big(\frac{\mathds{1}}{\mathds{1}-\overline{\lambda}}\big)_{T}$. Therefore a similar calculation to the above shows there is a bounded left inverse for $\big(\frac{\mathds{1}}{\mathds{1}-\lambda}\big)_{T}^{*}$. Hence $\big(\frac{\mathds{1}}{\mathds{1}-\lambda}\big)_{T}$ is invertible. To prove inequality (\ref{boundlabelforlambdax2}), note that, by inequality (\ref{modifiedfirstlineininvertibilityargument2}), for all $v \in\mathcal{H}$, there exists a unique $w \in\mathcal{H}$ such that $v=\big(\frac{\mathds{1}}{\mathds{1}-\lambda}\big)_{T}^{-1}w$ and
\begin{align}
\bigg\| \bigg(\dfrac{\mathds{1}}{\mathds{1}-\lambda}\bigg)_{T} v \bigg\|_{\mathcal{H}} 	& = \bigg\|\bigg(\dfrac{\mathds{1}}{\mathds{1}-\lambda}\bigg)_{T} \bigg(\dfrac{\mathds{1}}{\mathds{1}-\lambda}\bigg)_{T}^{-1}w \bigg\|_{\mathcal{H}} \nonumber \\
						& = \| w \|_{\mathcal{H}} \nonumber \\
						&\geq \underset{j=1,\hdots, d}{\min}\mathrm{Re}\bigg(\dfrac{1}{1-\lambda_{j}}\bigg)\bigg\| \bigg(\dfrac{\mathds{1}}{\mathds{1}-\lambda}\bigg)_{T}^{-1}w\bigg\|_{\mathcal{H}} \label{invertibleoperatorboundsforprojection2}. \\ \nonumber
\end{align}
By assumption, for $j=1,\hdots,d$, $\mathrm{Re}(\lambda_{j})<1$ and so $\mathrm{Re}\bigg(\dfrac{1}{1-\lambda_{j}}\bigg)>0$. For $j=1,\hdots,d$, 
\begin{equation}\label{realpartofoperator1minuslambda}
\mathrm{Re}\bigg(\dfrac{1}{1-\lambda_{j}}\bigg)= \mathrm{Re}\bigg(\dfrac{1-\overline{\lambda_{j}}}{\lvert 1-\lambda_{j}\rvert^{2}}\bigg)= \dfrac{\mathrm{Re}(1-\overline{\lambda_{j}})}{\lvert 1-\lambda_{j}\rvert^{2}}= \dfrac{\mathrm{Re}(1-\lambda_{j})}{\lvert 1-\lambda_{j}\rvert^{2}}.
\end{equation}
Thus, by the relations (\ref{invertibleoperatorboundsforprojection2}) and (\ref{realpartofoperator1minuslambda}) we obtain
\begin{align*}
    \bigg\| \bigg(\dfrac{\mathds{1}}{\mathds{1}-\lambda}\bigg)_{T}^{-1}w \bigg\|_{\mathcal{H}} &\leq \dfrac{\|w\|_{\mathcal{H}}}{\underset{j=1,\hdots,d}{\min}\mathrm{Re}\bigg(\dfrac{1}{1-\lambda_{j}}\bigg)} \\
    &\leq \underset{j=1,\hdots,d}{\max}\dfrac{\lvert 1-\lambda_{j}\rvert^{2}}{1-\mathrm{Re}(\lambda_{j})}\|w\|_{\mathcal{H}}.
\end{align*}
This implies inequality (\ref{boundlabelforlambdax2}).  \end{proof}

\begin{corollary}\label{invertibleoperatorsforlambdacorollary}
Let $T=(T_{1},\hdots, T_{d})$ be a d-tuple of operators on a Hilbert space $\mathcal{H}$ such that $\sum_{j=1}^{d}T_{j}=1_{\mathcal{H}}$ and $0\leq T_{j}\leq 1$ for $j=1,\hdots, d$. For all $\lambda=(\lambda_{1},\hdots,\lambda_{d})\in\mathbb{C}^{d}$ such that $\mathrm{Re}(\lambda_{j})>0$ for $j=1,\hdots,d$, the operators   
$$(\lambda)_{T}~\text{and}~\bigg(\dfrac{\mathds{1}}{\lambda}\bigg)_{T}$$
on $\mathcal{B(H)}$ are invertible and
\begin{align*}
\Big\| (\lambda)_{T}^{-1} \Big\|_{\mathcal{B(H)}} & \leq \dfrac{1}{\underset{j=1,\hdots, d}{\min}\mathrm{Re}(\lambda_{j})} \\
\Bigg\| \bigg(\dfrac{\mathds{1}}{\lambda}\bigg)_{T}^{-1} \Bigg\|_{\mathcal{B(H)}} &\leq ~ \max_{j=1,\hdots,d}\dfrac{\lvert \lambda_{j} \rvert^{2}}{\mathrm{Re}(\lambda_{j})}.
\end{align*}
\end{corollary}
\begin{proof}
By Theorem \ref{invertibleoperatorsforlambda}, for $\mu=(\mu_{1},\hdots,\mu_{d})\in\mathbb{C}^{d}$ such that $\mathrm{Re}(\mu_{j})<1$ for $j=1,\hdots,d$, the operators
$$(\mathds{1}-\mu)_{T}~\text{and}~\bigg(\dfrac{\mathds{1}}{\mathds{1}-\mu}\bigg)_{T}$$
on $\mathcal{B(H)}$ are invertible. Thus, for each $\lambda\in\mathbb{C}^{d}$ such that $\mathrm{Re}(\lambda_{j})>0$ for $j=1,\hdots,d$, an element $\mu\in\mathbb{C}^{d}$ defined by $\mathds{1}-\mu=\lambda$ has the property that $\mathrm{Re}(\mu_{j})<1$ for $j=1,\hdots,d$. The stated results follow.\end{proof}

\begin{proposition}\label{generalizedmodelprop}
Let $\mathcal{M} = \mathcal{M}_{1}\oplus \hdots \oplus \mathcal{M}_{d}$ be a complex Hilbert space and let $\mathcal{N}$ be a closed subspace of $\mathcal{M}$.  Let $X=(X_{1},\hdots, X_{d}), B=(B_{1},\hdots, B_{d})$ and $ Y=(Y_{1},\hdots, Y_{d})$ be  $d$-tuples of operators as in Lemma \ref{projectiononcdlemma}, where $X_{j}\in \mathcal{B}(\mathcal{N})$, $Y_{j} \in \mathcal{B}(\mathcal{N}^{\perp})$ and $B_{j} \in \mathcal{B}(\mathcal{N}^{\perp},\mathcal{N})$ for all integers such that $1\leq j \leq d$, and the following conditions are satisfied
\begin{enumerate}
  \item  $0 \leq X_{j} \leq 1$, $ 0 \leq Y_{j} \leq 1$;
  \item $\sum_{j=1}^{d} X_{j} = 1_{\mathcal{N}},~\sum_{j=1}^{d} B_{j} = 0, ~ \sum_{j=1}^{d} Y_{j} = 1_{\mathcal{N}^{\perp}}$;
  \item $B_{i}Y_{j} =  \delta_{ij}B_{j}-X_{i}B_{j}$, $B_{i}^{*}X_{j} =  \delta_{ij}B_{j}^{*}-Y_{i}B_{j}^{*} $;
  \item $B_{i}B^{*}_{j} =  \delta_{ij}X_{j}-X_{i}X_{j}, ~B^{*}_{i}B_{j} =  \delta_{ij}Y_{j}-Y_{i}Y_{j}$.
\end{enumerate}
Then for all $\lambda, \mu \in \mathbb{D}^{d}$,
\begin{equation}\label{mubadjointoperatoridentity}
\mu_{B^{*}}(1_{\mathcal{N}}-\lambda_{X})^{-1} = \bigg(\dfrac{\mu}{\mathds{1}-\lambda}\bigg)_{B^{*}} - \bigg( \dfrac{\mu}{\mathds{1}-\lambda}\bigg)_{Y}\bigg(\dfrac{\mathds{1}}{\mathds{1}-\lambda}\bigg)_{Y}^{-1}\bigg(\dfrac{\mathds{1}}{\mathds{1}-\lambda}\bigg)_{B^{*}}
\end{equation}
and
\begin{equation}\label{lambdaypluslambdabstar}
\lambda_{Y}+\lambda_{B^{*}}\big(1_{\mathcal{N}}-\lambda_{X}\big)^{-1}\lambda_{B} = 1_{\mathcal{N}^{\perp}}-\bigg(\dfrac{\mathds{1}}{\mathds{1}-\lambda}\bigg)_{Y}^{-1}.
\end{equation}
\end{proposition}
\begin{proof}
We firstly prove equation (\ref{mubadjointoperatoridentity}). By assumption  (2), $\sum_{j=1}^{d} X_{j} = 1_{\mathcal{N}}$,
\begin{align*}
\bigg(\dfrac{\mu}{\mathds{1}-\lambda}\bigg)_{B^{*}}\big(1_{\mathcal{N}}-\lambda_{X}\big) &= \bigg(\sum_{i=1}^{d} \dfrac{\mu_{i}}{1-\lambda_{i}} B_{i}^{*} \bigg)\bigg(1_{\mathcal{N}}-\sum_{j=1}^{d}\lambda_{j}X_{j}\bigg) \\
&= \bigg(\sum_{i=1}^{d} \dfrac{\mu_{i}}{1-\lambda_{i}} B_{i}^{*} \bigg)\bigg(\sum_{j=1}^{d}(1-\lambda_{j})X_{j}\bigg).
\end{align*}
Thus
\begin{align}
\bigg(\dfrac{\mu}{\mathds{1}-\lambda}\bigg)_{B^{*}}&\big(1_{\mathcal{N}}-\lambda_{X}\big)\nonumber\\
&= \dfrac{\mu_{1}}{1-\lambda_{1}}B_{1}^{*}\sum_{j=1}^{d}(1-\lambda_{j})X_{j}+\hdots+\dfrac{\mu_{d}}{1-\lambda_{d}}B_{d}^{*}\sum_{j=1}^{d}(1-\lambda_{j})X_{j} \nonumber\\
&= \dfrac{\mu_{1}}{1-\lambda_{1}}\sum_{j=1}^{d}(1-\lambda_{j})B_{1}^{*}X_{j}+\hdots+\dfrac{\mu_{d}}{1-\lambda_{d}}\sum_{j=1}^{d}(1-\lambda_{j})B_{d}^{*}X_{j}.\label{equation1234}  
\end{align}
By equation (\ref{propertiesoftheprojectionoperators3}), equation (\ref{equation1234}) can be written
\begin{align*}
&\bigg(\dfrac{\mu}{\mathds{1}-\lambda}\bigg)_{B^{*}}\big(1_{\mathcal{N}}-\lambda_{X}\big)\\
&= \dfrac{\mu_{1}}{1-\lambda_{1}}\sum_{j=1}^{d}(1-\lambda_{j})\Big[ \delta_{1j}B_{j}^{*}-Y_{1}B_{j}^{*}\Big]\\
&\hspace{5cm}+\hdots+\dfrac{\mu_{d}}{1-\lambda_{d}}\sum_{j=1}^{d}(1-\lambda_{j})\Big[ \delta_{dj}B_{j}^{*}-Y_{d}B_{j}^{*}\Big]\\
& = \dfrac{\mu_{1}}{1-\lambda_{1}}\bigg((1-\lambda_{1})(1-Y_{1})B_{1}^{*}-(1-\lambda_{2})Y_{1}B_{2}^{*}-\hdots-(1-\lambda_{d})Y_{1}B_{d}^{*}\bigg) \\
& +\hdots + \dfrac{\mu_{d}}{1-\lambda_{d}}\bigg(-(1-\lambda_{1})Y_{d}B_{1}^{*}-(1-\lambda_{2})Y_{d}B_{2}^{*}-\hdots+(1-\lambda_{d})(1-Y_{d})B_{d}^{*}\bigg) \\
&=\mu_{1}B_{1}^{*}-\dfrac{\mu_{1}}{1-\lambda_{1}}\bigg(Y_{1}(1-\lambda_{1})B_{1}^{*}+\hdots+Y_{1}(1-\lambda_{d})B_{d}^{*}\bigg)+\\
&~~~~~~~~~~~~~~~~~~\hdots+\mu_{d}B_{d}^{*}-\dfrac{\mu_{d}}{1-\lambda_{d}}\bigg((1-\lambda_{1})Y_{d}B_{1}^{*}+\hdots+(1-\lambda_{d})Y_{d}B_{d}^{*}\bigg). \\ 
\end{align*}
With Definition \ref{lambdaincdastupleofoperators} applied appropriately to each operator, we deduce that
\begin{align}
&\bigg(\dfrac{\mu}{\mathds{1}-\lambda}\bigg)_{B^{*}}\big(1_{\mathcal{N}}-\lambda_{X}\big)\nonumber \\ 
&=\mu_{B^{*}} - \Bigg[ \dfrac{\mu_{1}Y_{1}}{1-\lambda_{1}}\bigg((1-\lambda_{1})B_{1}^{*}+\hdots+(1-\lambda_{d})B_{d}^{*}\bigg)+\hdots \nonumber \\
&\hspace{5cm}+\dfrac{\mu_{d}Y_{d}}{1-\lambda_{d}}\bigg((1-\lambda_{1})B_{1}^{*}+\hdots+(1-\lambda_{d})B_{d}^{*}\bigg)\Bigg]\nonumber\\
& = \mu_{B^{*}} - \bigg[ \dfrac{\mu_{1}Y_{1}}{1-\lambda_{1}}+\hdots+\dfrac{\mu_{d}Y_{d}}{1-\lambda_{d}}\bigg](\mathds{1}-\lambda)_{B^{*}}\nonumber \\
& = \mu_{B^{*}} - \bigg(\dfrac{\mu}{\mathds{1}-\lambda}\bigg)_{Y}(\mathds{1}-\lambda)_{B^{*}}\label{nearlythelastequationinproof}.
\end{align}
By assumption (2), $\sum_{j=1}^{d}B_{j}=0$ so that 
$$(\mathds{1}-\lambda)_{B^{*}} = \sum_{j=1}^{d}(1-\lambda_{j})B_{j}^{*} = -\sum_{j=1}^{d}\lambda_{j}B_{j}^{*} = -\lambda_{B^{*}}. $$ 
Therefore, equation (\ref{nearlythelastequationinproof}) can be rearranged to 
\begin{equation}\label{verynearlythelastequation}
\mu_{B^{*}} = \bigg(\dfrac{\mu}{\mathds{1}-\lambda}\bigg)_{B^{*}}\big(1_{\mathcal{N}}-\lambda_{X}\big) - \bigg(\dfrac{\mu}{\mathds{1}-\lambda}\bigg)_{Y}\lambda_{B^{*}}.
\end{equation}
Take  $\mu=\mathds{1}$ in equation (\ref{verynearlythelastequation}) to obtain
$$\bigg(\dfrac{\mathds{1}}{\mathds{1}-\lambda}\bigg)_{B^{*}}(1_{\mathcal{N}}-\lambda_{X}) = \bigg(\dfrac{\mathds{1}}{\mathds{1}-\lambda}\bigg)_{Y}\lambda_{B^{*}}.$$
By Theorem \ref{invertibleoperatorsforlambda}, $(1_{\mathcal{N}}-\lambda_{X})^{-1}$ and $\Big(\dfrac{\mathds{1}}{\mathds{1}-\lambda}\Big)_{Y}^{-1}$ exist, therefore
\begin{equation}\label{invertibleidentityforthisproof}
\lambda_{B^{*}}(1_{\mathcal{N}}-\lambda_{X})^{-1} =\bigg(\dfrac{\mathds{1}}{\mathds{1}-\lambda}\bigg)_{Y}^{-1}\bigg(\dfrac{\mathds{1}}{\mathds{1}-\lambda}\bigg)_{B^{*}}.
\end{equation}
Multiply equation (\ref{verynearlythelastequation}) on the right by $(1-\lambda_{X})^{-1}$ and then use equation (\ref{invertibleidentityforthisproof}) to show the validity of equation (\ref{mubadjointoperatoridentity}), namely
\begin{align*}
\mu_{B^{*}}(1_{\mathcal{N}}-\lambda_{X})^{-1} &= \bigg(\dfrac{\mu}{\mathds{1}-\lambda}\bigg)_{B^{*}}- \bigg(\dfrac{\mu}{\mathds{1}-\lambda}\bigg)_{Y}\lambda_{B^{*}}(1_{\mathcal{N}}-\lambda_{X})^{-1}\\
			&= \bigg(\dfrac{\mu}{\mathds{1}-\lambda}\bigg)_{B^{*}}- \bigg(\dfrac{\mu}{\mathds{1}-\lambda}\bigg)_{Y}\bigg(\dfrac{\mathds{1}}{\mathds{1}-\lambda}\bigg)_{Y}^{-1}\bigg(\dfrac{\mathds{1}}{\mathds{1}-\lambda}\bigg)_{B^{*}}.
\end{align*}

To prove equation (\ref{lambdaypluslambdabstar}), let $\mu=\lambda$ in equation (\ref{mubadjointoperatoridentity}) to obtain
$$ \lambda_{B^{*}}(1_{\mathcal{N}}-\lambda_{X})^{-1} = \bigg(\dfrac{\lambda}{\mathds{1}-\lambda}\bigg)_{B^{*}}- \bigg(\dfrac{\lambda}{\mathds{1}-\lambda}\bigg)_{Y}\bigg(\dfrac{\mathds{1}}{\mathds{1}-\lambda}\bigg)_{Y}^{-1}\bigg(\dfrac{\mathds{1}}{\mathds{1}-\lambda}\bigg)_{B^{*}}.$$
Then, multiplying on the right by $\lambda_{B}$, we have
\begin{equation}\label{beginningequationforsecondidentity}
\lambda_{B^{*}}(1_{\mathcal{N}}-\lambda_{X})^{-1}\lambda_{B} = \bigg(\dfrac{\lambda}{\mathds{1}-\lambda}\bigg)_{B^{*}}\lambda_{B}- \bigg(\dfrac{\lambda}{\mathds{1}-\lambda}\bigg)_{Y}\bigg(\dfrac{\mathds{1}}{\mathds{1}-\lambda}\bigg)_{Y}^{-1}\bigg(\dfrac{\mathds{1}}{\mathds{1}-\lambda}\bigg)_{B^{*}}\lambda_{B}.
\end{equation}
By some subtle manipulations, we can reduce equation (\ref{beginningequationforsecondidentity}) as follows. Since
\begin{align*}
    &\lambda_{B^{*}}(1_{\mathcal{N}}-\lambda_{X})^{-1}\lambda_{B} \\
    &= \bigg(\dfrac{\lambda}{\mathds{1}-\lambda}\bigg)_{B^{*}}\lambda_{B}- \bigg(\dfrac{\mathds{1}-(\mathds{1}-\lambda)}{\mathds{1}-\lambda}\bigg)_{Y}\bigg(\dfrac{\mathds{1}}{\mathds{1}-\lambda}\bigg)_{Y}^{-1}\bigg(\dfrac{\mathds{1}}{\mathds{1}-\lambda}\bigg)_{B^{*}}\lambda_{B}, 
\end{align*}
and, by assumption (2), $\sum_{j=1}^{d}Y_{j}=1_{\mathcal{N}^{\perp}}$ and $\sum_{j=1}^{d} B_{j}=0$. Therefore,
\begin{align}
&\lambda_{B^{*}}(1_{\mathcal{N}}-\lambda_{X})^{-1}\lambda_{B} \nonumber \\
&= \bigg(\dfrac{\lambda}{\mathds{1}-\lambda}\bigg)_{B^{*}}\lambda_{B} - \bigg(\dfrac{\mathds{1}}{\mathds{1}-\lambda}\bigg)_{B^{*}}\lambda_{B} + \Big(1_{\mathcal{N}^{\perp}}\Big)\bigg(\dfrac{\mathds{1}}{\mathds{1}-\lambda}\bigg)_{Y}^{-1}\bigg(\dfrac{\mathds{1}}{\mathds{1}-\lambda}\bigg)_{B^{*}}\lambda_{B} \nonumber \\
& = \bigg(\dfrac{\mathds{1}-(\mathds{1}-\lambda)}{\mathds{1}-\lambda}\bigg)_{B^{*}}\lambda_{B} - \bigg(\dfrac{\mathds{1}}{\mathds{1}-\lambda}\bigg)_{B^{*}}\lambda_{B} +\bigg(\dfrac{\mathds{1}}{\mathds{1}-\lambda}\bigg)_{Y}^{-1}\bigg(\dfrac{\mathds{1}}{\mathds{1}-\lambda}\bigg)_{B^{*}}\lambda_{B} \nonumber\\
& = \mathds{1}_{B^{*}}\lambda_{B} +\bigg(\dfrac{\mathds{1}}{\mathds{1}-\lambda}\bigg)_{Y}^{-1}\bigg(\dfrac{\mathds{1}}{\mathds{1}-\lambda}\bigg)_{B^{*}}\lambda_{B}\nonumber \\
& = \bigg(\dfrac{\mathds{1}}{\mathds{1}-\lambda}\bigg)_{Y}^{-1}\bigg(\dfrac{\mathds{1}}{\mathds{1}-\lambda}\bigg)_{B^{*}}\lambda_{B}. \label{vnearlylastequationinthisproof}\\ \nonumber
\end{align}
Now, for each term in the operator $\mu_{B^{*}}$, let $\mu_{j}=\dfrac{1}{1-\lambda_{j}}$ for $j=1,\hdots,d$, so that
$$\mu_{B^{*}} = \sum_{j=1}^{d} \mu_{j}B_{j}^{*} = \sum_{j=1}^{d} \dfrac{1}{1-\lambda_{j}}B_{j}^{*} = \bigg(\dfrac{\mathds{1}}{\mathds{1}-\lambda}\bigg)_{B^{*}}.$$
From Proposition \ref{mulambdaidentitiesforgeneralizedmodels}, use the identity $\mu_{B^{*}}\lambda_{B} = (\mu\lambda)_{Y}-\mu_{Y}\lambda_{Y}$ with $\mu_{j}=\dfrac{1}{1-\lambda_{j}}$ for $j=1,\hdots,d$ to see that
$$ \bigg(\dfrac{\mathds{1}}{\mathds{1}-\lambda}\bigg)_{B^{*}}\lambda_{B} = \bigg(\dfrac{\lambda}{\mathds{1}-\lambda}\bigg)_{Y}-\bigg(\dfrac{\mathds{1}}{\mathds{1}-\lambda}\bigg)_{Y}\lambda_{Y}.$$
Thus, equation (\ref{vnearlylastequationinthisproof}) reduces to
\begin{align*}
\lambda_{B^{*}}(1_{\mathcal{N}}-\lambda_{X})^{-1}\lambda_{B} &= \bigg(\dfrac{\mathds{1}}{\mathds{1}-\lambda}\bigg)_{Y}^{-1}\Bigg[ \bigg(\dfrac{\lambda}{\mathds{1}-\lambda}\bigg)_{Y}-\bigg(\dfrac{\mathds{1}}{\mathds{1}-\lambda}\bigg)_{Y}\lambda_{Y}\Bigg] \\
 & = \bigg(\dfrac{\mathds{1}}{\mathds{1}-\lambda}\bigg)_{Y}^{-1}\bigg(\dfrac{\mathds{1}-(\mathds{1}-\lambda)}{1-\lambda}\bigg)_{Y}-\lambda_{Y} \\
 & = 1_{\mathcal{N}^{\perp}}-\bigg(\dfrac{\mathds{1}}{\mathds{1}-\lambda}\bigg)_{Y}^{-1} - \lambda_{Y},
\end{align*}
from which, after rearranging terms, we find equation (\ref{lambdaypluslambdabstar}).\end{proof}

\section{Analyticity of some operator-valued functions on $\mathbb{D}^{d}$}\label{analytic_I}

In Definition \ref{genmodel} of a generalization of models of functions from $\sad$ we consider operator-valued functions
$I(\cdot)$. The study of properties of  operator-valued functions $I(\cdot)$ shall be the objective of this section.
\begin{definition}
Let $\mathcal{H}$ be a Hilbert space. We say that 
$f:\mathbb{D}^{d} \to \mathcal{B}(\mathcal{H})$
is a contractive analytic $\mathcal{B}(\mathcal{H})$-valued  function on $\mathbb{D}^{d}$ if $f$ is analytic on $\mathbb{D}^{d}$ and, for all $\lambda\in\mathbb{D}^{d}$, 
$$ \| f(\lambda) \|_{\mathcal{B(H)}}\leq 1.$$ 
\end{definition}

We were introduced to our first result for operator-valued functions on $\mathbb{D}^{d}$ of the form $(f(\lambda))_{T}$ with Theorem \ref{invertibleoperatorsforlambda}. We now expand upon that theorem and include some essential results which will be used in later sections. For the following results, we define 
$$\mathbb{C}_{+}^{d} = \big\{(\lambda_{1},\hdots,\lambda_{d})\in\mathbb{C}^{d} : \mathrm{Re}(\lambda_{j})>0~~\text{for}~~j=1,\hdots,d\big\}.$$
Recall Definition \ref{lambdaincdastupleofoperators}. Let $\mathcal{M}$ be a complex Hilbert space and let $\lambda\in\mathbb{C}^{d}$. Then if $T=(T_{1},T_{2},\hdots,T_{d})$ is a $d$-tuple of operators on $\mathcal{M}$, we recall the operator 
\begin{equation}\label{lambdaXequation}
    \lambda_{T}=\lambda_{1}T_{1}+\hdots+\lambda_{d}T_{d}
\end{equation}.
\begin{proposition}\label{propforIbeinganalytic}
Let $T=(T_{1},\hdots, T_{d})$ be a d-tuple of operators on a Hilbert space $\mathcal{H}$ such that $\sum_{j=1}^{d}T_{j}=1_{\mathcal{H}}$ and $0\leq T_{j}\leq 1$ for $j=1,\hdots, d$. Then
\begin{enumerate}
    \item the map $\lambda\mapsto\lambda_{T}$ from $\mathbb{C}^{d}$ into $\mathcal{B}(\mathcal{H})$ is analytic in $\lambda$;\label{this}
    \item $\bigg(\dfrac{\mathds{1}}{\lambda}\bigg)_{T}$ is analytic in $\lambda$ for all $\lambda\in\mathbb{C}^{d}$ such that $\lambda_{j}\neq 0$ for $j=1,\hdots,d$;
    \item $\bigg(\dfrac{\mathds{1}}{\lambda}\bigg)^{-1}_{T}$ is analytic  in $\lambda$ for all $\lambda\in\mathbb{C}^{d}_{+}$;
    \item $\bigg(\dfrac{\mathds{1}}{\mathds{1}-\lambda}\bigg)^{-1}_{T}$ is analytic  in $\lambda$ for all $\lambda\in\mathbb{C}^{d}$ such that $\mathrm{Re}(\lambda_{j})<1$ for $j=1,\hdots,d$.
\end{enumerate}
\end{proposition}
\begin{proof}
We shall show that each of the operators above is analytic  in $\lambda$ in  their respective domains. To prove (1), let $f:\mathbb{C}^{d}\rightarrow\mathcal{B(H)}$ be defined by $f(\lambda) = \lambda_{T}$. We claim that, for every $\lambda\in\mathbb{C}^{d}$,
$$f'(\lambda):\mathbb{C}^{d}\rightarrow\mathcal{B(H)} : x\mapsto x_{T}$$
is the Fr\'{e}chet derivative. Note, for any $\lambda$ and $x$ in $\mathbb{C}^{d}$, 
$$(\lambda+x)_{T}=(\lambda_{1}+x_{1})T_{1}+\hdots+(\lambda_{d}+x_{d})T_{d},$$
and so $(\lambda+x)_{T}=\lambda_{T}+x_{T}.$ Indeed, $f$ is a bounded {\em linear} map from $\mathbb{C}^{d}$ to $\mathcal{B}(\mathcal{H})$, and so $f$ is its own Fr\'{e}chet derivative.

To prove (2), for all $\lambda\in\mathbb{C}^{d}$ such that $\lambda_{j}\neq 0$ for $j=1,\hdots,d$, the Fr\'{e}chet derivative $D\bigg(\dfrac{\mathds{1}}{\lambda}\bigg)_{T}:\mathbb{C}^{d}\rightarrow \mathcal{B(H)}$ is given by
$$D\bigg(\dfrac{\mathds{1}}{\lambda}\bigg)_{T}(x) = -\sum_{j=1}^{d}\dfrac{x_{j}}{\lambda_{j}^{2}}T_{j},~~\text{for}~x\in\mathbb{C}^{d}.$$
Note that, for $\lambda\in\mathbb{C}^{d}$ such that $\lambda_{j}\neq 0$ for $j=1,\hdots,d$ and for $x\in\mathbb{C}^{d}$ with $x\neq 0$ such that $x_{j}\neq -\lambda_{j}$ for $j=1,\hdots,d$,
\begin{align}
    &\dfrac{1}{\|x\|_{\mathbb{C}^{d}}}\bigg\|\bigg(\dfrac{\mathds{1}}{\lambda+x}\bigg)_{T}-\bigg(\dfrac{\mathds{1}}{\lambda}\bigg)_{T}-D\bigg(\dfrac{\mathds{1}}{\lambda}\bigg)_{T}(x)\bigg\|_{\mathcal{B(H)}}\nonumber \\
    &= \dfrac{1}{\|x\|_{\mathbb{C}^{d}}}\bigg\|\sum_{j=1}^{d}\dfrac{x_{j}}{(\lambda_{j}+x_{j})(\lambda_{j})}T_{j}-\sum_{j=1}^{d}\dfrac{x_{j}}{\lambda_{j}^{2}}T_{j}\bigg\|_{\mathcal{B(H)}}\nonumber\\
    &=\dfrac{1}{\|x\|_{\mathbb{C}^{d}}}\bigg\|\sum_{j=1}^{d}\dfrac{x_{j}}{\lambda_{j}}\bigg(\dfrac{1}{\lambda_{j}+x_{j}}-\dfrac{1}{\lambda_{j}}\bigg)T_{j}\bigg\|_{\mathcal{B(H)}}\nonumber\\
    &=\dfrac{1}{\|x\|_{\mathbb{C}^{d}}}\bigg\|\sum_{j=1}^{d}\dfrac{-x_{j}^{2}}{\lambda_{j}^{2}(\lambda_{j}+x_{j})}T_{j}\bigg\|_{\mathcal{B(H)}}\nonumber\\
    &\leq \dfrac{1}{\|x\|_{\mathbb{C}^{d}}}\sum_{j=1}^{d}\bigg\lvert\dfrac{x_{j}}{\lambda_{j}}\bigg\rvert^{2}\bigg\lvert\dfrac{1}{\lambda_{j}+x_{j}}\bigg\rvert \|T_{j}\|_{\mathcal{B(H)}}\nonumber\\
    &\leq \dfrac{1}{\|x\|_{\mathbb{C}^{d}}}\sum_{j=1}^{d}\bigg\lvert\dfrac{x_{j}}{\lambda_{j}}\bigg\rvert^{2}\bigg\lvert\dfrac{1}{\lambda_{j}+x_{j}}\bigg\rvert.\label{hellothere}
\end{align}
Since $\lambda,x\in\mathbb{C}^{d}$ are such that $x_{j}\neq -\lambda_{j}$ for $j=1,\hdots,d$ and $\lambda_{j}\neq 0$ for $j=1,\hdots,d$, we have $\lambda+x\neq 0$. Thus the summation of inequality (\ref{hellothere}) is bounded. For $j=1,\hdots,d$, 
\begin{align*}
    \lvert x_{j}\rvert &\leq \|x\|_{\infty}, \\
    \lvert \lambda_{j}\rvert -\lvert x_{j}\rvert&\geq \lvert \lambda_{j}\rvert-\|x\|_{\infty}.
\end{align*}
Thus, for $j=1,\hdots,d$, 
$$\dfrac{1}{\lvert \lambda_{j}\rvert-\lvert x_{j}\rvert}\leq \dfrac{1}{\lvert \lambda_{j}\rvert - \|x\|_{\infty}},$$
therefore we have
$$\bigg\lvert\dfrac{1}{\lambda_{j}+x_{j}}\bigg\rvert \leq \dfrac{1}{\lvert \lambda_{j} \rvert - \lvert x_{j} \rvert}\leq \dfrac{1}{\lvert \lambda_{j}\rvert - \|x\|_{\infty}}.$$
It is easy to see that $\|x\|_{\infty}\leq \|x\|_{\mathbb{C}^{d}}\leq \sqrt{d}\|x\|_{\infty}$. Thus, for any $x \in \mathbb{C}^d$ such that  $x_{j}\neq -\lambda_{j}$ for $j=1,\hdots,d$,
$$ \dfrac{1}{\lvert\lambda_{j}\rvert-\|x\|_{\infty}}\leq \dfrac{1}{\lvert \lambda_{j}\rvert-\|x\|_{\mathbb{C}^{d}}}.$$
Therefore from inequality (\ref{hellothere}),
\begin{align}
    &\lim_{\|x\|_{\mathbb{C}^{d}}\rightarrow 0}\dfrac{1}{\|x\|_{\mathbb{C}^{d}}}\bigg\|\bigg(\dfrac{\mathds{1}}{\lambda+x}\bigg)_{T}-\bigg(\dfrac{\mathds{1}}{\lambda}\bigg)_{T}-D\bigg(\dfrac{\mathds{1}}{\lambda}\bigg)_{T}(x)\bigg\|_{\mathcal{B(H)}}\nonumber \\
    &\leq \lim_{\|x\|_{\mathbb{C}^{d}}\rightarrow 0}\dfrac{\|x\|^{2}_{\infty}}{\|x\|_{\mathbb{C}^{d}}}\sum_{j=1}^{d}\bigg\lvert\dfrac{1}{\lambda_{j}}\bigg\rvert^{2}\dfrac{1}{\lvert \lambda_{j}\rvert - \|x\|_{\infty}} \nonumber\\
    &\leq \lim_{\|x\|_{\mathbb{C}^{d}}\rightarrow 0}\dfrac{\|x\|^{2}_{\mathbb{C}^{d}}}{\|x\|_{\mathbb{C}^{d}}}\sum_{j=1}^{d}\bigg\lvert\dfrac{1}{\lambda_{j}}\bigg\rvert^{2}\dfrac{1}{\lvert \lambda_{j}\rvert - \|x\|_{\mathbb{C}^{d}}} \nonumber\\
    &\leq \lim_{\|x\|_{\mathbb{C}^{d}}\rightarrow 0}\|x\|_{\mathbb{C}^{d}}\sum_{j=1}^{d}\bigg\lvert\dfrac{1}{\lambda_{j}}\bigg\rvert^{2}\dfrac{1}{\lvert \lambda_{j}\rvert - \|x\|_{\mathbb{C}^{d}}} \nonumber\\
    &=0.\nonumber
\end{align}
Therefore $\bigg(\dfrac{\mathds{1}}{\lambda}\bigg)_{T}$ is analytic for all $\lambda\in\mathbb{C}^{d}$ such that $\lambda_{j}\neq 0$ for $j=1,\hdots,d$ with
$$D\bigg(\dfrac{\mathds{1}}{\lambda}\bigg)_{T}(x) = -\bigg(\dfrac{x}{\lambda^{2}}\bigg)_{T}$$
for $x\in\mathbb{C}^{d}$ such that $x_{j}\neq -\lambda_{j}$ for $j=1,\hdots,d$.

Let us prove (3). By Corollary \ref{invertibleoperatorsforlambdacorollary}, for all $\lambda\in\mathbb{C}^{d}$ such that $\mathrm{Re}(\lambda_{j})>0$ for $j=1,\hdots,d$,  $\bigg(\dfrac{\mathds{1}}{\lambda}\bigg)_{T}$ is invertible. 
Note that, by
\cite[Proposition I.2.6]{completenormedalgebras}, the map 
$$g\footnote{The notation $\mathrm{Inv}(\mathcal{B(H)})$ denotes the group of regular elements of the Banach algebra $\mathcal{B(H)}$, with the topology defined by the operator norm.}
:\mathrm{Inv}(\mathcal{B(H)}) \rightarrow \mathrm{Inv}(\mathcal{B(H)})$$ 
such that $g(x)= x^{-1}$ is analytic and $Dg(x) h= -x^{-1} h x^{-1}$. 
Moreover, we have already shown in (2) that $\bigg(\dfrac{\mathds{1}}{\lambda}\bigg)_{T}$ is analytic for all $\lambda\in\mathbb{C}_{+}^{d}$. Observe that
$$\bigg(\dfrac{\mathds{1}}{\lambda}\bigg)^{-1}_{T} = g\circ \bigg(\dfrac{\mathds{1}}{\lambda}\bigg)_{T}.$$
Since the composition of analytic maps is also analytic,  $\bigg(\dfrac{\mathds{1}}{\lambda}\bigg)^{-1}_{T}$ is analytic. By the chain rule for Fr\'{e}chet derivatives
\cite[Proposition 5.3.7]{TheoreticalNumericalAnalysis},
$$D\bigg(\dfrac{\mathds{1}}{\lambda}\bigg)^{-1}_{T}(x) = \bigg(\dfrac{\mathds{1}}{\lambda}\bigg)^{-1}_{T}\bigg(\dfrac{x}{\lambda^{2}}\bigg)_{T}\bigg(\dfrac{\mathds{1}}{\lambda}\bigg)^{-1}_{T}$$
for all $x\in\mathbb{C}^{d}_{+}$.

Let us show (4). By Corollary \ref{invertibleoperatorsforlambdacorollary}, for $\mu=(\mu_{1},\hdots,\mu_{d})\in\mathbb{C}^{d}$ such that $\mathrm{Re}(\mu_{j})>0$ for $j=1,\hdots,d$, the operator $\bigg(\dfrac{\mathds{1}}{\mu}\bigg)_{T}$ is invertible and 
$\bigg(\dfrac{\mathds{1}}{\mu}\bigg)^{-1}_{T}$
is analytic in $\mu$ by (3). Thus, for every element $\lambda\in\mathbb{C}^{d}$ such that $\mathrm{Re}(\lambda_{j})<1$ for $j=1,\hdots,d$, $1-\lambda$ has the property that $\mathrm{Re}(1-\lambda_{j})>0$ for $j=1,\hdots,d$. Therefore the inverse of the operator $\bigg(\dfrac{\mathds{1}}{\mathds{1}-\lambda}\bigg)_{T}$ exists. The analyticity of the operator then follows similarly to (3) with 
\begin{equation}\label{frechetderivfor1over1minuslambda}
D\bigg(\dfrac{\mathds{1}}{\mathds{1}-\lambda}\bigg)^{-1}_{T}(x) = -\bigg(\dfrac{\mathds{1}}{\mathds{1}-\lambda}\bigg)^{-1}_{T}\bigg(\dfrac{x}{(\mathds{1}-\lambda)^{2}}\bigg)_{T}\bigg(\dfrac{\mathds{1}}{\mathds{1}-\lambda}\bigg)^{-1}_{T}
\end{equation}
for $x\in\mathbb{C}^{d}_{+}$.\end{proof}

\section{Properties of a special operator-valued function on $\mathbb{D}^{d}$}\label{Ioflambda}
In this section we shall show that the  inner operator-valued function $I(\cdot)$ defined by equation \eqref{defI}
is analytic and inner. This operator-valued function $I(\cdot)$ will be used in our construction of a desingularized model  for  $\varphi\in\mathcal{SA}_{d}$ with a carapoint $\tau\in\mathbb{T}^{d}$.
\begin{proposition}\label{Ilambdaproperties}
    Let $\mathcal{M}$ be a complex separable  Hilbert space and let $\tau\in\mathbb{T}^{d}$. Define
    \begin{equation}\label{defI}
I(\lambda)=1_{\mathcal{M}}-\bigg(\dfrac{\mathds{1}}{\mathds{1}-\overline{\tau}\lambda}\bigg)^{-1}_{Y}~~\text{for}~\lambda\in\overline{\mathbb{D}}^{d}~\text{such that}~\lambda_{j}\neq\tau_{j}~\text{for}~j=1,\hdots,d,
\end{equation}
where $(Y_{1},\hdots,Y_{d})$ is a $d$-tuple of positive contractions on $\mathcal{M}$ such that $\sum_{j=1}^{d}Y_{j}=1_{\mathcal{M}}$. Then the following properties are satisfied.
\begin{enumerate}
    \item $I$ is analytic on $\mathbb{D}^{d}$;
    \item $\|I(\lambda)\|_{\mathcal{B(M)}}<1$ for all $\lambda\in\mathbb{D}^{d}$;
    \item at every point $\lambda\in\mathbb{T}^{d}$ such that $\lambda_{j}\neq\tau_{j}$ for $j=1,\hdots,d$, $I(\lambda)$ exists, 
    $$I^{*}(\lambda)I(\lambda)=I(\lambda)I^{*}(\lambda)=1_{\mathcal{M}}~~\text{and}~~\|I(\lambda)\|_{\mathcal{B(M)}}=1;$$
    \item $I:\mathbb{D}^{d}\rightarrow \mathcal{B(M)}$ is an inner operator-valued function.
\end{enumerate}
\end{proposition}
\begin{proof}
    Let us prove statement (1). For all $\lambda\in\overline{\mathbb{D}}^{d}$ such that $\lambda_{j}\neq \tau_{j}$, $\mathrm{Re}(\overline{\tau_{j}}\lambda_{j})<1$ for $j=1,2,\hdots,d$. By Proposition \ref{propforIbeinganalytic} (iv), the operator $I(\lambda)$ is well defined on $\overline{\mathbb{D}}^{d}$ except when $\lambda_{j}=\tau_{j}$ for some $j=1,2,\hdots,d$ and is analytic on $\mathbb{D}^{d}$. By Propositon \ref{propforIbeinganalytic} (iv) and equation (\ref{frechetderivfor1over1minuslambda}), the Fr\"{e}chet derivative of $I(\lambda)$ is given by   
  \begin{align}
    DI(\lambda)(h) &= -\bigg(\dfrac{\mathds{1}}{\mathds{1}-\overline{\tau}\lambda}\bigg)^{-1}_{Y}\bigg(\sum_{j=1}^{d}\dfrac{\overline{\tau_{j}}h_{j}}{(1-\overline{\tau_{j}}\lambda_{j})^{2}}Y_{j}\bigg)\bigg(\dfrac{\mathds{1}}{\mathds{1}-\overline{\tau}\lambda}\bigg)^{-1}_{Y}\nonumber \\
    &= -\bigg(\dfrac{\mathds{1}}{\mathds{1}-\overline{\tau}\lambda}\bigg)^{-1}_{Y}\bigg(\dfrac{{\tau}^{*}h}{(\mathds{1}-\overline{\tau}\lambda)^{2}}\bigg)_{Y}\bigg(\dfrac{\mathds{1}}{\mathds{1}-\overline{\tau}\lambda}\bigg)^{-1}_{Y}.\label{thisequationrightherelolz}
\end{align}
Therefore by equation (\ref{thisequationrightherelolz}), $I(\lambda)$ is analytic on $\mathbb{D}^{d}$ and at each $\lambda\in\mathbb{T}^{d}$ such that $\lambda_{j} \neq \tau_{j}$ for all $j=1,\hdots,d.$. Thus $I$ is analytic on $\mathbb{D}$ and analytic almost everywhere on $\mathbb{T}^{d}$. 

We prove statement (2), that is, $\|I(\lambda)\|_{\mathcal{B(M)}}<1$ for all $\lambda\in\mathbb{D}^{d}$; note that
\begin{align}
    1_{\mathcal{M}}-I^{*}(\lambda)I(\lambda)&=1_{\mathcal{M}}-\Bigg[1_{\mathcal{M}}-\bigg(\dfrac{\mathds{1}}{\mathds{1}-\overline{\tau}\lambda}\bigg)_{Y}^{-1}\Bigg]^{*}\Bigg[1_{\mathcal{M}}-\bigg(\dfrac{\mathds{1}}{\mathds{1}-\overline{\tau}\lambda}\bigg)_{Y}^{-1}\Bigg]\nonumber \\
    &=\bigg(\dfrac{\mathds{1}}{\mathds{1}-\overline{\tau}\lambda}\bigg)_{Y}^{-1}+\Bigg[\bigg(\dfrac{\mathds{1}}{\mathds{1}-\overline{\tau}\lambda}\bigg)_{Y}^{-1}\Bigg]^{*}-\Bigg[\bigg(\dfrac{\mathds{1}}{\mathds{1}-\overline{\tau}\lambda}\bigg)_{Y}^{-1}\Bigg]^{*}\bigg(\dfrac{\mathds{1}}{\mathds{1}-\overline{\tau}\lambda}\bigg)_{Y}^{-1}\nonumber\\
    &=\bigg(\dfrac{\mathds{1}}{\mathds{1}-\overline{\tau}\lambda}\bigg)_{Y}^{-1}+\bigg(\dfrac{\mathds{1}}{\mathds{1}-\tau\overline{\lambda}}\bigg)_{Y}^{-1}-\bigg(\dfrac{\mathds{1}}{\mathds{1}-\tau\overline{\lambda}}\bigg)_{Y}^{-1}\bigg(\dfrac{\mathds{1}}{\mathds{1}-\overline{\tau}\lambda}\bigg)_{Y}^{-1}\nonumber\\
    &=\bigg(\dfrac{\mathds{1}}{\mathds{1}-\tau\overline{\lambda}}\bigg)_{Y}^{-1}\Bigg[1_{\mathcal{M}}+\bigg(\dfrac{\mathds{1}}{\mathds{1}-\tau\overline{\lambda}}\bigg)_{Y}\bigg(\dfrac{\mathds{1}}{\mathds{1}-\overline{\tau}\lambda}\bigg)_{Y}^{-1}-\bigg(\dfrac{\mathds{1}}{\mathds{1}-\tau\overline{\lambda}}\bigg)_{Y}^{-1}\Bigg]\nonumber\\
    &= \bigg(\dfrac{\mathds{1}}{\mathds{1}-\tau\overline{\lambda}}\bigg)_{Y}^{-1}\Bigg[\bigg(\dfrac{\mathds{1}}{\mathds{1}-\tau\overline{\lambda}}\bigg)_{Y}+\bigg(\dfrac{\mathds{1}}{\mathds{1}-\overline{\tau}\lambda}\bigg)_{Y}-1_{\mathcal{M}}\Bigg]\bigg(\dfrac{\mathds{1}}{\mathds{1}-\overline{\tau}\lambda}\bigg)_{Y}^{-1}.\label{CstarAC}
\end{align}
One can check that
\begin{align}\bigg(\dfrac{\mathds{1}}{\mathds{1}-\tau\overline{\lambda}}\bigg)_{Y}+\bigg(\dfrac{\mathds{1}}{\mathds{1}-\overline{\tau}\lambda}\bigg)_{Y}-1_{\mathcal{M}} &= \sum_{i=1}^{d}\dfrac{1}{1-\tau_{i}\overline{\lambda_{i}}}Y_{i}+\sum_{i=1}^{d}\dfrac{1}{1-\overline{\tau_{i}}\lambda_{i}}Y_{i}-\sum_{i=1}^{d}Y_{i}\nonumber\\
&=\sum_{i=1}^{d}\bigg(\dfrac{1}{1-\tau_{i}\overline{\lambda_{i}}}+\dfrac{1}{1-\overline{\tau_{i}}\lambda_{i}} -1\bigg)Y_{i} \nonumber\\
&=\sum_{i=1}^{d}\dfrac{1-\lvert \lambda_{i} \rvert^{2}}{\lvert 1-\overline{\tau_{i}}\lambda_{i}\rvert^{2}}Y_{i}=\bigg(\dfrac{\mathds{1}-\lvert \lambda \rvert^{2}}{\lvert \mathds{1}-\overline{\tau}\lambda\rvert^{2}}\bigg)_{Y}>0\label{CstarC2}
\end{align}
for all $\lambda\in\mathbb{D}^{d}$. Then equation (\ref{CstarAC}) and equation (\ref{CstarC2}) imply 
$$1_{\mathcal{M}}-I^{*}(\lambda)I(\lambda)=\bigg(\dfrac{\mathds{1}}{\mathds{1}-\tau\overline{\lambda}}\bigg)_{Y}^{-1}
\bigg(\dfrac{1-\lvert \lambda \rvert^{2}}{\lvert \mathds{1}-\overline{\tau}\lambda\rvert^{2}}\bigg)_{Y}\bigg(\dfrac{\mathds{1}}{\mathds{1}-\overline{\tau}\lambda}\bigg)_{Y}^{-1}.$$
Hence $1_{\mathcal{M}}-I^{*}(\lambda)I(\lambda)$ is of the form $X^{*}AX$, where $X, A\in\mathcal{B(\mathcal{M})}$, $X=\bigg(\dfrac{\mathds{1}}{\mathds{1}-\overline{\tau}\lambda}\bigg)_{Y}^{-1} $ is invertible and $A= \bigg(\dfrac{1-\lvert \lambda \rvert^{2}}{\lvert \mathds{1}-\overline{\tau}\lambda\rvert^{2}}\bigg)_{Y}$ is a strictly positive operator for all $\lambda\in\mathbb{D}^{d}$. Thus, by \cite[Lemma 2.81]{AMYOperatorbook}, $1_{\mathcal{M}}-I^{*}(\lambda)I(\lambda)>0$ for all $\lambda\in\mathbb{D}^{d}$.

Let us prove statement (3). Rearrange equation (\ref{defI}) to read
$$I(\lambda)\Big(\dfrac{\mathds{1}}{\mathds{1}-\overline{\tau}\lambda}\Big)_{Y} = \Big(\dfrac{\mathds{1}}{\mathds{1}-\overline{\tau}\lambda}\Big)_{Y} - 1_{\mathcal{M}},$$
and
$$\bigg(I(\lambda)\Big(\dfrac{\mathds{1}}{\mathds{1}-\overline{\tau}\lambda}\Big)_{Y}\bigg)^{*} = \Big(\dfrac{\mathds{1}}{\mathds{1}-\tau\overline{\lambda}}\Big)_{Y}I^{*}(\lambda).$$
Then 
\begin{align*}\Big(\dfrac{\mathds{1}}{\mathds{1}-\tau\overline{\lambda}}\Big)_{Y}I^{*}(\lambda)&I(\lambda)\Big(\dfrac{\mathds{1}}{\mathds{1}-\overline{\tau}\lambda}\Big)_{Y} = \bigg[\Big(\dfrac{\mathds{1}}{\mathds{1}-\tau\overline{\lambda}}\Big)_{Y}-1_{\mathcal{M}}\bigg]\bigg[\Big(\dfrac{\mathds{1}}{\mathds{1}-\overline{\tau}\lambda}\Big)_{Y}-1_{\mathcal{M}}\bigg]\\
    & = \Big(\dfrac{\mathds{1}}{\mathds{1}-\tau\overline{\lambda}}\Big)_{Y}\Big(\dfrac{\mathds{1}}{\mathds{1}-\overline{\tau}\lambda}\Big)_{Y}-\Big(\dfrac{\mathds{1}}{\mathds{1}-\tau\overline{\lambda}}\Big)_{Y}-\Big(\dfrac{\mathds{1}}{\mathds{1}-\overline{\tau}\lambda}\Big)_{Y}+1_{\mathcal{M}} \hspace{5cm}\\
    & = \Big(\dfrac{\mathds{1}}{\mathds{1}-\tau\overline{\lambda}}\Big)_{Y}\Big(\dfrac{\mathds{1}}{\mathds{1}-\overline{\tau}\lambda}\Big)_{Y} + 1_{\mathcal{M}} - \sum_{j=1}^{d}\bigg(\dfrac{1}{1-\tau_{j}\overline{\lambda_{j}}}+\dfrac{1}{1-\overline{\tau_{j}}\lambda_{j}}\bigg)Y_{j}\\
    & = \Big(\dfrac{\mathds{1}}{\mathds{1}-\tau\overline{\lambda}}\Big)_{Y}\Big(\dfrac{\mathds{1}}{\mathds{1}-\overline{\tau}\lambda}\Big)_{Y} + 1_{\mathcal{M}} - \sum_{j=1}^{d}\dfrac{2(1-\mathrm{Re}(\overline{\tau_{j}}\lambda_{j}))}{\lvert 1-\overline{\tau_{j}}\lambda_{j}\rvert^{2}}Y_{j}.
\end{align*}
Since $\tau\neq\lambda\in\mathbb{T}^{d}$ is such that $\tau_{j}\neq\lambda_{j}$ for $j=1,\hdots,d$, we have $\lvert 1-\overline{\tau}_{j}\lambda_{j}\rvert^{2}=2(1-\mathrm{Re}(\overline{\tau}_{j}\lambda_{j}))$. By assumption, $\sum_{j=1}^{d}Y_{j}=1_{\mathcal{M}}$. Hence, for $\lambda\in\mathbb{T}^{d}$ such that $\lambda_{j} \neq \tau_{j}$ for $j=1,\hdots,d$,
$$\Big(\dfrac{\mathds{1}}{\mathds{1}-\tau\overline{\lambda}}\Big)_{Y}I^{*}(\lambda)I(\lambda)\Big(\dfrac{\mathds{1}}{\mathds{1}-\overline{\tau}\lambda}\Big)_{Y} = \Big(\dfrac{\mathds{1}}{\mathds{1}-\tau\overline{\lambda}}\Big)_{Y}\Big(\dfrac{\mathds{1}}{\mathds{1}-\overline{\tau}\lambda}\Big)_{Y}$$
and so
\begin{equation}\label{innerfunctionproof}
    0=\Big(\dfrac{\mathds{1}}{\mathds{1}-\tau\overline{\lambda}}\Big)_{Y}\Big[1_{\mathcal{M}}-I^{*}(\lambda)I(\lambda)\Big]\Big(\dfrac{\mathds{1}}{\mathds{1}-\overline{\tau}\lambda}\Big)_{Y}.
\end{equation}
By Theorem \ref{invertibleoperatorsforlambda}, for $\lambda\in\mathbb{C}^{d}$ with $\mathrm{Re}(\overline{\tau}_{j}\lambda_{j})<1$ and $\mathrm{Re}(\tau_{j}\overline{\lambda_{j}})<1$ for $j=1,\hdots,d$, the operators
$$\Big(\dfrac{\mathds{1}}{\mathds{1}-\tau\overline{\lambda}}\Big)_{Y}~~\text{and}~~\Big(\dfrac{\mathds{1}}{\mathds{1}-\overline{\tau}\lambda}\Big)_{Y},$$
are invertible. We guarantee these conditions for all $\lambda\in\mathbb{T}^{d}$ such that $\lambda_{j}\neq\tau_{j}$ for $j=1,\hdots, d$. Therefore upon multiplying equation (\ref{innerfunctionproof}) by the appropriate invertible operators and rearranging it, we obtain
\begin{equation}\label{IstarIfinalformula}
I^{*}(\lambda)I(\lambda) = 1_{\mathcal{M}}.
\end{equation}
Likewise, if we intend to calculate 
\begin{equation}\label{IIstarcalcforuni}
I(\lambda)I^{*}(\lambda) = \bigg(1_{\mathcal{M}}-\Big(\dfrac{\mathds{1}}{\mathds{1}-\overline{\tau}\lambda}\Big)^{-1}_{Y}\bigg)\bigg(1_{\mathcal{M}}-\Big(\dfrac{\mathds{1}}{\mathds{1}-\tau\overline{\lambda}}\Big)^{-1}_{Y}\bigg),
\end{equation}
for $\lambda\in\mathbb{T}^{d}$ such that $\lambda_{j}\neq \tau_{j}$ for $j=1,\hdots,d$, then let us do so by firstly expanding equation (\ref{IIstarcalcforuni}) and then multiplying from the right and the left as below
\begin{equation*}
    I(\lambda)I^{*}(\lambda)\bigg(\dfrac{\mathds{1}}{\mathds{1}-\tau\overline{\lambda}}\bigg)_{Y}=\big(\dfrac{\mathds{1}}{\mathds{1}-\tau\overline{\lambda}}\Big)_{Y}-1_{\mathcal{M}}-\Big(\dfrac{\mathds{1}}{\mathds{1}-\overline{\tau}\lambda}\Big)^{-1}_{Y}\Big(\dfrac{\mathds{1}}{\mathds{1}-\tau\overline{\lambda}}\Big)_{Y}+\Big(\dfrac{\mathds{1}}{\mathds{1}-\overline{\tau}\lambda}\Big)_{Y}^{-1},
\end{equation*}
and
\begin{equation*}
    \Big(\dfrac{\mathds{1}}{\mathds{1}-\overline{\tau}\lambda}\Big)_{Y}I(\lambda)I^{*}(\lambda)\bigg(\dfrac{\mathds{1}}{\mathds{1}-\tau\overline{\lambda}}\bigg)_{Y}
    =\Big(\dfrac{\mathds{1}}{\mathds{1}-\overline{\tau}\lambda}\Big)_{Y}\big(\dfrac{\mathds{1}}{\mathds{1}-\tau\overline{\lambda}}\Big)_{Y}-\Big(\dfrac{\mathds{1}}{\mathds{1}-\overline{\tau}\lambda}\Big)_{Y}-\Big(\dfrac{\mathds{1}}{\mathds{1}-\tau\overline{\lambda}}\Big)_{Y}+1_{\mathcal{M}}.
\end{equation*}
Therefore, for $\lambda\in\mathbb{T}^{d}$ such that $\tau_{j}\neq \lambda_{j}$ for $j=1,\hdots,d$,
\begin{align*}
\Big(\dfrac{\mathds{1}}{\mathds{1}-\overline{\tau}\lambda}\Big)_{Y}I(\lambda)I^{*}(\lambda)\bigg(\dfrac{\mathds{1}}{\mathds{1}-\tau\overline{\lambda}}\bigg)_{Y}&=1_{\mathcal{M}}-\Big(\dfrac{\mathds{1}}{\mathds{1}-\overline{\tau}\lambda}\Big)_{Y}-\Big(\dfrac{\mathds{1}}{\mathds{1}-\tau\overline{\lambda}}\Big)_{Y}+\Big(\dfrac{\mathds{1}}{\mathds{1}-\overline{\tau}\lambda}\Big)_{Y}\big(\dfrac{\mathds{1}}{\mathds{1}-\tau\overline{\lambda}}\Big)_{Y}\\
    &=1_{\mathcal{M}}-\sum_{j=1}^{d}\dfrac{2(1-\mathrm{Re}(\overline{\tau}_{j}\lambda_{j}))}{\lvert1-\overline{\tau}_{j}\lambda{j}\rvert^{2}}Y_{j}+\Big(\dfrac{\mathds{1}}{\mathds{1}-\overline{\tau}\lambda}\Big)_{Y}\big(\dfrac{\mathds{1}}{\mathds{1}-\tau\overline{\lambda}}\Big)_{Y}.
\end{align*}
Since $\lvert 1-\tau_{j}\lambda_{j}\rvert^{2}=2(1-\mathrm{Re}(\overline{\tau}_{j}\lambda_{j}))$ and $\sum_{j=1}^{d}Y_{j}=1_{\mathcal{M}}$, we have
\begin{equation}\label{IIstarfinalformula}
\Big(\dfrac{\mathds{1}}{\mathds{1}-\overline{\tau}\lambda}\Big)_{Y}I(\lambda)I^{*}(\lambda)\bigg(\dfrac{\mathds{1}}{\mathds{1}-\tau\overline{\lambda}}\bigg)_{Y}=\Big(\dfrac{\mathds{1}}{\mathds{1}-\overline{\tau}\lambda}\Big)_{Y}\big(\dfrac{\mathds{1}}{\mathds{1}-\tau\overline{\lambda}}\Big)_{Y}.
\end{equation}
Thus, for all $\lambda\in\mathbb{T}^{d}$ such that $\lambda_{j}\neq\tau_{j}$ for all $j=1,\hdots,d$,
$$I(\lambda)I^{*}(\lambda)=1_{\mathcal{M}}.$$
Together, equations (\ref{IstarIfinalformula}) and (\ref{IIstarfinalformula}) prove that $I(\lambda)$ is unitary for $\lambda\in\mathbb{T}^{d}$ such that $\lambda_{j}\neq \tau_{j}$ for $j=1,\hdots,d$. 

Let us show that the radial limit $\lim_{\lambda\rightarrow\tau}I(r\tau)$ exists. Indeed,
\begin{align*}
    \|I(r\tau)-1_{\mathcal{M}}\|_{\mathcal{B(M)}}&=\Bigg\| 1_{\mathcal{M}}-\bigg(\dfrac{\mathds{1}}{\mathds{1}-\overline{\tau}\tau r}\bigg)_{Y}^{-1}-1_{\mathcal{M}}\Bigg\|_{\mathcal{B(M)}} \\
&=\Bigg\| \bigg(\dfrac{\mathds{1}}{\mathds{1}-\overline{\tau}\tau r}\bigg)_{Y}^{-1} \Bigg\|_{\mathcal{B(M)}} \\
&=\Bigg\| \bigg(\sum_{i=1}^{d}\dfrac{1}{1-r}Y_{i}\bigg)^{-1} \Bigg\|_{\mathcal{B(M)}} \\
&=\Bigg\| (1-r)\bigg(\sum_{i=1}^{d}Y_{i}\bigg)^{-1} \Bigg\|_{\mathcal{B(M)}}.
\end{align*}
By the assumption $\sum_{j=1}^{d}Y_{j}=1_{\mathcal{M}}$, and so
\begin{align*}
\|I(r\tau)-1_{\mathcal{M}}\|_{\mathcal{B(M)}}=1-r.
\end{align*}
Therefore we have
$$\lim_{r\rightarrow 1^{-}} \|I(r\tau)-1_{\mathcal{M}}\|_{\mathcal{B(M)}} = 0.$$
Thus, $I:\mathbb{D}^{d}\rightarrow \mathcal{B(M)}$ defined by
$$I(\lambda) = 1_{\mathcal{M}}-\bigg(\dfrac{\mathds{1}}{\mathds{1}-\overline{\tau}\lambda}\bigg)_{Y}$$
is an inner operator-valued function. \end{proof}

\section{Boundary behaviour of $\varphi\in\mathcal{SA}_{d}$ via desingularized models}\label{main_result}

In this section, for a general positive integer $d$ and for a function $\varphi\in\mathcal{SA}_{d}$ with a carapoint $\tau\in\mathbb{T}^{d}$, we construct a desingularized model $(\mathcal{M},u,I)$ of $\varphi$ relative to $\tau$, 
 see Definition \ref{desingularized}. 
 
 Let us recall the definition of a generalized model of an analytic function $\varphi:\mathbb{D}^{d}\rightarrow \mathbb{C}$.
\begin{definition}{\normalfont{\cite[Definition 3.1]{ATY2}}} 
Let $\varphi :\mathbb{D}^{d} \rightarrow \mathbb{C}$ be analytic. The triple $(\mathcal{M},u,I)$ is said to be a {\em{generalized model}} of $\varphi$ if
\begin{enumerate}
\item $\mathcal{M}$~\text{is a separable Hilbert space}; 
\item $u : \mathbb{D}^{d} \rightarrow \mathcal{M}$~\text{is an analytic mapping}; 
\item $I$~\text{is a contractive analytic}~$\mathcal{B(M)}$\text{-valued function on}~$\mathbb{D}^{d}$ such that $\|I(\lambda)\|<1$ for all $\lambda\in\mathbb{D}^{d}$;
\item \text{For all}~$\lambda, \mu \in \mathbb{D}^{d}$,
\begin{equation}\label{generalizedmodelindef}
1-\overline{\varphi(\mu)}\varphi(\lambda) = \Big\langle\big(1_{\mathcal{M}}-I(\mu)^{*}I(\lambda)\big)u({\lambda}),u({\mu}) \Big\rangle_{\mathcal{M}}.
\end{equation}
\end{enumerate}
\end{definition} 

Let $I:\mathbb{D}^{d}\rightarrow \mathcal{B(M)}$ be an analytic mapping.
We say $I(\cdot)$ {\em extends analytically} to the point $\lambda\in\mathbb{T}^{d}$ if there exists an open neighbourhood $U$ of $\lambda$ and an analytic $\mathcal{B}(\mathcal{M})$-valued function $F$ on $U$ such that $I$ agrees with $F$ on $U\cap\mathbb{D}^{d}$. 
\begin{theorem}\label{generalizedmodeltheorem}
Let $\tau \in\mathbb{T}^{d}$ be a carapoint for $\varphi \in\mathcal{SA}_{d}$. There exists a generalized model $(\mathcal{M},u,I)$ of $\varphi$ such that 
\begin{enumerate}
\item $\tau$ is a $C$-point for $(\mathcal{M},u,I)$;
\item at every point $\lambda\in\mathbb{T}^{d}$ such that $\lambda_{j}\neq\tau_{j}$ for $j=1,\hdots,d$, $I$ extends analytically to $\lambda$ and $\|I(\lambda)\|_{\mathcal{B(M)}}=1$;
\item $\tau$ is a carapoint for $I$ and $\lim_{\lambda{\overset{nt}\rightarrow}\tau}I(\lambda)=1_{\mathcal{M}}$;
\item $I:\mathbb{D}^{d}\rightarrow \mathcal{B(M)}$ is an inner $\mathcal{B(M)}$-valued function.
\end{enumerate}
Moreover, we can express $I$ in the form
\begin{equation}\label{Ioperatorgeneralizedmodel}
I(\lambda) = 1_{\mathcal{M}} - \bigg( \dfrac{\mathds{1}}{\mathds{1}-\overline{\tau}\lambda}\bigg)_{Y}^{-1},\; \text{for all} \;\lambda\in\mathbb{D}^{d},
\end{equation}
and for some $d$-tuple of positive contractions $(Y_{1},\hdots,Y_{d})$ on $\mathcal{M}$ such that $\sum_{j=1}^{d}Y_{j}=1_{\mathcal{M}}$.
\end{theorem}

\begin{proof}
By \cite[Theorem 9.102]{AMYOperatorbook}, every function $\varphi\in\mathcal{SA}_{d}$ has a model $(\mathcal{L},v)$. Choose any model $(\mathcal{L},v)$ of $\varphi\in\mathcal{SA}_{d}$ and any realization $(\alpha,\beta, \gamma, D)$ of $(\mathcal{L},v)$. By definition, $\mathcal{L}$ is a separable Hilbert space with orthogonal decomposition $\mathcal{L} = \mathcal{L}_{1}\oplus \hdots \oplus \mathcal{L}_{d}$. Let $P=(P_{1},\hdots, P_{d})$ be the $d$-tuple of operators on $\mathcal{L}$ where $P_{j}:\mathcal{L}\rightarrow \mathcal{L}$ is the orthogonal projection onto $\mathcal{L}_{j}$ for $j=1,\hdots,d$. \index{$\mathcal{N}=\mathrm{Ker}(1-D\tau_{P})$}\index{$\mathcal{L}=\mathcal{N}\oplus \mathcal{N}^{\perp}$}

Define $\mathcal{N}=\mathrm{Ker}(1-D\tau_{P})$, where $\tau_{P}=\tau_{1}P_{1}+\hdots+\tau_{d}P_{d}$. Since $\tau\in\mathbb{T}^{d}$ is a carapoint for $\varphi$, we may apply Lemma \ref{gammainranperp} to deduce that $\gamma\in\mathrm{Ran}(1-D\tau_P)\subset \mathcal{N}^{\perp}$ and $\overline{\tau}_{P}\beta\in\mathcal{N}^{\perp}$. Note that $\mathcal{L}=\mathcal{N}\oplus \mathcal{N}^{\perp}$.

Firstly, consider the possibility that $\mathcal{N}=\{0\}$. Then $1_{\mathcal{L}}-D\tau$ is injective, hence there exists a unique vector $v(\tau)\in\mathcal{L}=\mathcal{N}\oplus\mathcal{N}^{\perp}=0\oplus \mathcal{N}^{\perp}$ such that,
\begin{equation}\label{gammacontradiction}
(1_{\mathcal{L}}-D\tau_{P})v(\tau)=
	\gamma.
\end{equation}
Let $\{\lambda^{(n)}\}$ be any sequence in $\mathbb{D}^{d}$ that converges nontangentially to $\tau$ as $n\rightarrow \infty$. We claim that $v(\lambda^{(n)})\rightarrow v(\tau)$. Suppose that $v(\lambda^{(n)})\nrightarrow v(\tau)$ as $n\rightarrow \infty$. By Corollary \ref{coro5.7cara}, $\{v(\lambda^{(n)})\}$ is a bounded sequence in $\mathcal{L}$. Recall that the closed ball in a separable  Hilbert space is sequentially compact in the weak topology. Thus there exists a subsequence $\{\lambda^{(n_k)}\}$ of $\{\lambda^{(n)}\}$ such that $\{v(\lambda^{(n_k)})\}$ is weakly convergent in $\mathcal{L}$ to an $x\in\mathcal{L}$ such that $x\neq v(\tau)$ with 
$$\lim_{k\rightarrow \infty} v(\lambda^{(n_k)})=x \in \mathcal{L}.$$
Moreover, by Proposition \ref{prop5.8cara}, $\|v(\lambda^{(n_k)})-x\|_{\mathcal{L}}\rightarrow 0$ as $k\rightarrow \infty$. Now,
$$\lim_{k \rightarrow \infty}(1_{\mathcal{L}}-D\lambda_{P}^{(n_k)})v(\lambda^{(n_k)})=(1_{\mathcal{L}}-D\tau_{P})x= \begin{bmatrix} 
	0 \\
	\gamma \\
\end{bmatrix},$$
but $v(\tau)$ is the unique vector in $\mathcal{N}^{\perp}$ such that $(1_{\mathcal{L}}-D\tau_{P})v(\tau)=\gamma$ and here, $x\neq v(\tau)$, which contradicts the uniqueness of the vector $v(\tau)$ in equation (\ref{gammacontradiction}). Therefore $v(\lambda^{(n)})\rightarrow v(\tau)$ as $n\rightarrow \infty$. Moreover, $v(\lambda)$ extends continuously to $S\cup\{\tau\}$ for any set $S$ that tends to $\tau$ nontangentially. Hence $\tau$ is a $C$-point for $(\mathcal{L},v)$. Therefore the conclusion of the theorem for the case $\mathrm{Ker}(1_{\mathcal{L}}-D\tau_{P})=\{0\}$ follows if we take $\mathcal{M}=\mathcal{L}$, $u=v$ and $I(\lambda)=\lambda_{P}$.

 Assume now that $\mathrm{Ker}(1_{\mathcal{L}}-D\tau_{P})\neq \{0\}$. With respect to the orthogonal decomposition $\mathcal{L}=\mathcal{N}\oplus \mathcal{N}^{\perp}$, where $\mathcal{N}=\mathrm{Ker}(1-D\tau_{P})$, $D\tau_{P}$ has the operator matrix
$$D\tau_{P} = \begin{bmatrix}
    A_{1} & A_{2} \\
    A_{3} & A_{4}
\end{bmatrix},$$
for contractions $A_{1}\in\mathcal{B(N)}$, $A_{2}\in\mathcal{B(N^{\perp},N)}$, $A_{3}\in\mathcal{B(N,N^{\perp})}$ and $A_{4}\in \mathcal{B(N^{\perp})}$. Let $w\in\mathcal{N}$ and define $P_{\mathcal{N}}:\mathcal{L}\rightarrow\mathcal{L}$ to be the orthogonal projection onto $\mathcal{N}$ and $D\tau_{P}|_{\mathcal{N}}$ to be the restriction of $D\tau_{P}$ to $\mathcal{N}$, that is, $D\tau_{P}|_{\mathcal{N}}: \mathcal{N}\rightarrow \mathcal{L}$. Note that any $w\in\mathrm{Ker}(1-D\tau_{P})$ if and only if $D\tau_{P} w = w$. Therefore, for any $w\in\mathcal{N}$,
$$w = P_{\mathcal{N}}D\tau_{P}|_{\mathcal{N}}w = P_{\mathcal{N}}\begin{bmatrix}
    A_{1}w \\
    A_{3}w
\end{bmatrix} = A_{1}w.$$
Hence $A_{1}=1_{\mathcal{N}}$. Therefore,  in the contractive matrix  $D\tau_{P}$, since $1_{\mathcal{N}}$ is unitary, we have $A_{2}=A_{3}=0$. Let $Q=A_{4}$. so that $Q$ is contraction on $\mathcal{N}^{\perp}$.  Then 
$$D\tau_{P} = \begin{bmatrix}
    1_{\mathcal{N}} & 0 \\
    0 & Q
\end{bmatrix}.$$
We claim that $\mathrm{Ker}(1_{\mathcal{N}^{\perp}}-Q)=\{0\}$. Suppose $\mathrm{Ker}(1_{\mathcal{N}^{\perp}}-Q)\neq\{0\}$: there exists a $u_{0}\in\mathcal{N^{\perp}}$ such that $(1_{\mathcal{N}^{\perp}}-Q)u_{0}=0$. However, for any 
 $v\in\mathcal{L}$, $(1_{\mathcal{L}}-D\tau_{P})v=0$ if and only if
 $v\in\mathcal{N}$.  Choosing $v=u_0$ we deduce a contradiction, therefore $\mathrm{Ker}(1_{\mathcal{N}^{\perp}}-Q)=\{0\}$.
 Express $\lambda_{P}\in\mathcal{B(L)}$ as an operator matrix with respect to the decomposition $\mathcal{L}=\mathcal{N}\oplus\mathcal{N}^{\perp}$ as in Lemma \ref{projectiononcdlemma} by
\begin{equation}\label{lambdasanoperator123}
\lambda_{P}=\lambda_{1}P_{1}+\hdots+\lambda_{d}P_{d}= \begin{bmatrix}
        \sum_{j=1}^{d} \lambda_{j}X_{j} & \sum_{j=1}^{d}\lambda_{j}B_{j}\\
        \sum_{j=1}^{d}\lambda_{j}B_{j}^{*} & \sum_{j=1}^{d}\lambda_{j}Y_{j} \\
\end{bmatrix} = \begin{bmatrix}
\lambda_{X} & \lambda_{B}\\
\lambda_{B^{*}} & \lambda_{Y}
\end{bmatrix}.
\end{equation}
Likewise, 
$$\tau_{P} = \begin{bmatrix}
\tau_{X} & \tau_{B}\\
\tau_{B^{*}} & \tau_{Y}
\end{bmatrix}~\text{and}~\overline{\tau}_{P} = \begin{bmatrix}
        \overline{\tau}_{X} & \overline{\tau}_{B}\\
        \overline{\tau}_{B^{*}} & \overline{\tau}_{Y}
        \end{bmatrix},$$
so that, for any $\lambda\in\mathbb{D}^{d}$, 
\begin{align*}
    1_{\mathcal{L}}-D\lambda_{P}  &= 1_{\mathcal{L}}-D\tau_{P}\tau_{P}^{*}\lambda_{P}\\
                &=1_{\mathcal{L}}-\begin{bmatrix}
                    1_{\mathcal{N}} & 0 \\
                    0 & Q
                    \end{bmatrix}\begin{bmatrix}
\overline{\tau}_{X}\lambda_{X}+\overline{\tau}_{B}\lambda_{B^{*}} & \overline{\tau}_{X}\lambda_{B} + \overline{\tau}_{B}\lambda_{Y}\\
\overline{\tau}_{B^{*}}\lambda_{X}+\overline{\tau}_{Y}\lambda_{B^{*}} & \overline{\tau}_{B^{*}}\lambda_{B}+\overline{\tau}_{Y}\lambda_{Y}
                    \end{bmatrix}.
\end{align*}
By the identities from Proposition \ref{mulambdaidentitiesforgeneralizedmodels},
\begin{align*}
    1_{\mathcal{L}}-D\lambda_{P} &=1_{\mathcal{L}}-\begin{bmatrix}
                    1_{\mathcal{N}} & 0 \\
                    0 & Q \\
                    \end{bmatrix}\begin{bmatrix}
                    (\overline{\tau}\lambda)_{X} & (\overline{\tau}\lambda)_{B} \\
                    (\overline{\tau}\lambda)_{B^{*}} & (\overline{\tau}\lambda)_{Y}\\
                    \end{bmatrix}\\
    &=\begin{bmatrix}
    1_{\mathcal{N}}-(\overline{\tau}\lambda)_{X} & -(\overline{\tau}\lambda)_{B} \\
    -Q(\overline{\tau}\lambda)_{B^{*}} & 1_{\mathcal{N}^{\perp}}-Q(\overline{\tau}\lambda)_{Y}\\
    \end{bmatrix}.
\end{align*}
For the vector $v$ in the model $(\mathcal{L},v)$ of the function $\varphi$ with $\mathcal{L}=\mathcal{N}\oplus \mathcal{N}^{\perp}$, let $w(\lambda)=P_{\mathcal{N}}v(\lambda)$ and $u(\lambda)=P_{\mathcal{N}^{\perp}}v(\lambda)$ for any $\lambda\in\mathbb{D}^{d}$, so that $v(\lambda)=w(\lambda)+u(\lambda)$ for $\lambda\in\mathbb{D}^{d}$ and
\begin{align*}
    \begin{bmatrix}
    0 \\    
    \gamma \\
    \end{bmatrix}   &=(1_{\mathcal{L}}-D\lambda_{P})\begin{bmatrix}
                        w(\lambda)\\   
                        u(\lambda) 
                    \end{bmatrix}\\
                    &=\begin{bmatrix}
    1_{\mathcal{N}}-(\overline{\tau}\lambda)_{X} & -(\overline{\tau}\lambda)_{B} \\
    -Q(\overline{\tau}\lambda)_{B^{*}} & 1_{\mathcal{N}^{\perp}}-Q(\overline{\tau}\lambda)_{Y}\\
    \end{bmatrix}\begin{bmatrix}
                        w(\lambda)\\   
                        u(\lambda) 
                    \end{bmatrix}\\
                &=\begin{bmatrix}
    (1_{\mathcal{N}}-(\overline{\tau}\lambda)_{X})w(\lambda) -(\overline{\tau}\lambda)_{B}u(\lambda) \\
    -Q(\overline{\tau}\lambda)_{B^{*}}w(\lambda)+(1_{\mathcal{N}^{\perp}}-Q(\overline{\tau}\lambda)_{Y})u(\lambda)\\
    \end{bmatrix}
\end{align*}
for all $\lambda\in\mathbb{D}^{d}$. Hence the following system of equations holds,
\begin{align}
0 &= (1_{\mathcal{N}}-(\overline{\tau}\lambda)_{X})w(\lambda)-(\overline{\tau}\lambda)_{B}u(\lambda) \label{equationonegeneralizedmodel}\\
	\gamma &= -Q(\overline{\tau}\lambda)_{B^{*}}w(\lambda)+(1_{\mathcal{N}^{\perp}}-Q(\overline{\tau}\lambda)_{Y})u(\lambda)
        \label{equationtwogeneralizedmodel} \\ \nonumber
\end{align}
for all $\lambda\in\mathbb{D}^{d}$. From equation (\ref{equationonegeneralizedmodel}), 
\begin{equation}\label{importantequationlol}
(1_{\mathcal{N}}-(\overline{\tau}\lambda)_{X})w(\lambda) = (\overline{\tau}\lambda)_{B}u(\lambda)
\end{equation}
for all $\lambda\in\mathbb{D}^{d}$ and, furthermore, Theorem \ref{invertibleoperatorsforlambda} is applicable, with $\lambda$ replaced by $\overline{\tau}\lambda$, so that
\begin{equation}\label{wlambageneralized}
w(\lambda)=(1_{\mathcal{N}}-(\overline{\tau}\lambda)_{X})^{-1}(\overline{\tau}\lambda)_{B}u(\lambda)
\end{equation}
for all $\lambda\in\mathbb{D}^{d}$. Therefore, by equations (\ref{equationtwogeneralizedmodel}) and (\ref{wlambageneralized}), we obtain 
\begin{align}
\gamma  & = -Q(\overline{\tau}\lambda)_{B^{*}}w(\lambda)+(1_{\mathcal{N}^{\perp}}-Q(\overline{\tau}\lambda)_{Y})u(\lambda) \nonumber\\
        & = -Q(\overline{\tau}\lambda)_{B^{*}}(1_{\mathcal{N}}-(\overline{\tau}\lambda)_{X})^{-1}(\overline{\tau}\lambda)_{B}u(\lambda) + (1_{\mathcal{N}^{\perp}}-Q(\overline{\tau}\lambda)_{Y})u(\lambda)\nonumber\\
        & = \Big[ (1_{\mathcal{N}^{\perp}}-Q(\overline{\tau}\lambda)_{Y})-Q(\overline{\tau}\lambda)_{B^{*}}(1_{\mathcal{N}}-(\overline{\tau}\lambda)_{X})^{-1}(\overline{\tau}\lambda)_{B}\Big]u(\lambda)\nonumber\\ 
        & = \Big[ 1_{\mathcal{N}^{\perp}}-Q\Big((\overline{\tau}\lambda)_{Y}+(\overline{\tau}\lambda)_{B^{*}}(1_{\mathcal{N}}-(\overline{\tau}\lambda)_{X})^{-1}(\overline{\tau}\lambda)_{B}\Big)\Big]u(\lambda)
        \label{lastlineforgammageneralized}
\end{align}
for all $\lambda\in\mathbb{D}^{d}$. Again with $\lambda$ replaced by $\overline{\tau}\lambda$, we can use Proposition \ref{generalizedmodelprop}, equation (\ref{lambdaypluslambdabstar}), to show that 
\begin{equation}\label{lambdatauoverliney}
(\overline{\tau}\lambda)_{Y}+(\overline{\tau}\lambda)_{B^{*}}\big(1_{\mathcal{N}}-(\overline{\tau}\lambda)_{X}\big)^{-1}(\overline{\tau}\lambda)_{B} = 1_{\mathcal{N}^{\perp}}-\bigg(\dfrac{\mathds{1}}{\mathds{1}-\overline{\tau}\lambda}\bigg)_{Y}^{-1}.
\end{equation}
Therefore equation (\ref{lambdatauoverliney}) implies that equation (\ref{lastlineforgammageneralized}) takes the form
\begin{equation}\label{gammaequals1-qIu}
\gamma = \bigg[1_{\mathcal{N}^{\perp}}-Q\Big(1_{\mathcal{N}^{\perp}}-\Big(\dfrac{\mathds{1}}{\mathds{1}-\overline{\tau}\lambda}\Big)^{-1}_{Y}\Big)\bigg]u(\lambda) = \big(1_{\mathcal{N}^{\perp}}-QI(\lambda)\big)u(\lambda),
\end{equation}
where, for all $\lambda\in\overline{\mathbb{D}^{d}}$ such that $\lambda_{j}\neq \tau_{j}$ for $j=1,\hdots,d$, we define $I(\lambda)$ by
\begin{equation}\label{Ilambdaoperatorgeneralized}
I(\lambda) =1_{\mathcal{N}^{\perp}}-\bigg(\dfrac{\mathds{1}}{\mathds{1}-\overline{\tau}\lambda}\bigg)^{-1}_{Y}\in\mathcal{B(N^{\perp})}.
\end{equation}
By Proposition \ref{Ilambdaproperties}, $I$ is analytic on $\mathbb{D}^{d}$, $\|I(\lambda)\|_{\mathcal{B(\mathcal{N}^{\perp})}}<1$ for all $\lambda\in\mathbb{D}^{d}$ and $\|I(\lambda)\|_{\mathcal{B(\mathcal{N}^{\perp})}}=1$ at every point $\lambda\in\mathbb{T}^{d}$ such that, for all $j=1,\hdots,d$, $\lambda_{j}\neq\tau_{j}$.

Let us show that equation  (\ref{generalizedmodelindef}) holds for $\varphi$. For the vector $v(\lambda)$ in the model $(\mathcal{L},v)$ of the function $\varphi$ with $\mathcal{L}=\mathcal{N}\oplus \mathcal{N}^{\perp}$, let $w(\lambda)=P_{\mathcal{N}}v(\lambda)$ and $u(\lambda)=P_{\mathcal{N}^{\perp}}v(\lambda)$. Using equation
 (\ref{lambdasanoperator123}), we can perform the calculation
\begin{align}
    \overline{\tau}_{P}\lambda_{P}v(\lambda)   &= \overline{\tau}_{P}\begin{bmatrix}
                                            \lambda_{X} & \lambda_{B} \\
                                            \lambda_{B^{*}} & \lambda_{Y} \\
                                            \end{bmatrix}\begin{bmatrix}
                                                w(\lambda)\\
                                                u(\lambda)
                                                \end{bmatrix}\nonumber\\ 
                                        &=\begin{bmatrix}
                                            \overline{\tau}_{X} & \overline{\tau}_{B} \\
                                            \overline{\tau}_{B^{*}} & \overline{\tau}_{Y} \\
                                            \end{bmatrix}\begin{bmatrix}
                                            \lambda_{X}w(\lambda)+ \lambda_{B}u(\lambda) \\
                                            \lambda_{B^{*}}w(\lambda)+\lambda_{Y}u(\lambda) \\
                                            \end{bmatrix}\nonumber\\
                                        &=\begin{bmatrix}
                                           \overline{\tau}_{X} ( \lambda_{X}w(\lambda)+ \lambda_{B}u(\lambda) )+\overline{\tau}_{B}(\lambda_{B^{*}}w(\lambda)+\lambda_{Y}u(\lambda)) \\
                                            \overline{\tau}_{B^{*}} ( \lambda_{X}w(\lambda)+ \lambda_{B}u(\lambda) )+\overline{\tau}_{Y}(\lambda_{B^{*}}w(\lambda)+\lambda_{Y}u(\lambda))\\
                                            \end{bmatrix}\nonumber\\
                                        &=\begin{bmatrix}(\overline{\tau}_{X}\lambda_{X}+\overline{\tau}_{B}\lambda_{B^{*}})w(\lambda)+(\overline{\tau}_{X}\lambda_{B}+\overline{\tau}_{B}\lambda_{Y})u(\lambda)\\(\overline{\tau}_{B^{*}}\lambda_{X}+\overline{\tau}_{Y}\lambda_{B^{*}})w(\lambda)+(\overline{\tau}_{B^{*}}\lambda_{B}+\overline{\tau}_{Y}\lambda_{Y})u(\lambda)\label{finalmatrixinthisequation}                                        \end{bmatrix}.
\end{align}
Once more, use the identities from Proposition \ref{mulambdaidentitiesforgeneralizedmodels} with $\lambda$
replaced by $\overline{\tau}\lambda$ in equation (\ref{finalmatrixinthisequation}) to obtain
\begin{equation}\label{taustarlambdavlambda}
\overline{\tau}_{P}\lambda_{P}v(\lambda)=\begin{bmatrix}
        (\overline{\tau}\lambda)_{X}w(\lambda)+(\overline{\tau}\lambda)_{B}u(\lambda)\\
        (\overline{\tau}\lambda)_{B^{*}}w(\lambda)+(\overline{\tau}\lambda)_{Y}u(\lambda) \\
\end{bmatrix}_{\mathcal{N}\oplus{\mathcal{N}^{\perp}}}.
\end{equation}
Observe the following. By equation (\ref{importantequationlol}), 
$$(\overline{\tau}\lambda)_{B}u(\lambda)=w(\lambda)-(\overline{\tau}\lambda)_{X}w(\lambda).$$
Therefore
$$(\overline{\tau}\lambda)_{X}w(\lambda)+(\overline{\tau}\lambda)_{B}u(\lambda)=(\overline{\tau}\lambda)_{X}w(\lambda)+w(\lambda)-(\overline{\tau}\lambda)_{X}w(\lambda)=w(\lambda)$$
for $\lambda\in\mathbb{D}^{d}$. By equation (\ref{wlambageneralized}) and equation (\ref{lambdatauoverliney}), 
$$(\overline{\tau}\lambda)_{B^{*}}w(\lambda)+(\overline{\tau}\lambda)_{Y}u(\lambda) = \Big((\overline{\tau}\lambda)_{Y}+(\overline{\tau}\lambda)_{B^{*}}(1_{\mathcal{N}}-(\overline{\tau}\lambda)_{X})^{-1}(\overline{\tau}\lambda)_{B}\Big)u(\lambda)=I(\lambda)u(\lambda)$$
for $\lambda\in\mathbb{D}^{d}$. Hence, by equation (\ref{taustarlambdavlambda}),
\begin{equation}\label{importantforthegeneralizedmodel}
\overline{\tau}_{P}\lambda_{P}v(\lambda) =\begin{bmatrix}
w(\lambda) \\
I(\lambda)u(\lambda) \\
\end{bmatrix}
\end{equation}
for $\lambda\in\mathbb{D}^{d}$.
 With equation (\ref{importantforthegeneralizedmodel}), we can derive the generalized model formula from the simpler model formula of Definition \ref{modelpolydefinto}:
\begin{align}
    1-\overline{\varphi(\mu)}\varphi(\lambda) &= \langle (1_{\mathcal{L}}-\mu_{P}^{*}\lambda_{P})v(\lambda),v(\mu)\rangle_{\mathcal{L}}\nonumber\\
    &= \langle v(\lambda),v(\mu)\rangle_{\mathcal{L}} - \langle \lambda_{P}v(\lambda),\mu_{P}v(\mu)\rangle_{\mathcal{L}}\nonumber\\
    &= \Bigg\langle  \begin{bmatrix}
                w(\lambda) \\
                u(\lambda)
                \end{bmatrix},\begin{bmatrix}
                w(\mu) \\
                u(\mu)
                \end{bmatrix}\Bigg\rangle_{\mathcal{L}}- \langle \tau_{P}\overline{\tau}_{P}\lambda_{P}v(\lambda),\mu_{P}v(\mu)\rangle_{\mathcal{L}}\nonumber\\
    &=\langle w(\lambda),w(\mu)\rangle_{\mathcal{N}}+\langle u(\lambda),u(\mu)\rangle_{\mathcal{N}^{\perp}}-\langle \overline{\tau}_{P}\lambda_{P}v(\lambda),\overline{\tau}_{P}\mu_{P}v(\mu)\rangle_{\mathcal{L}}\nonumber\\
    &=\langle w(\lambda),w(\mu)\rangle_{\mathcal{N}}+\langle u(\lambda),u(\mu)\rangle_{\mathcal{N}^{\perp}}-\Bigg\langle  \begin{bmatrix}
                w(\lambda) \\
                I(\lambda)u(\lambda)
                \end{bmatrix},\begin{bmatrix}
                w(\mu) \\
                I(\mu)u(\mu)
                \end{bmatrix}\Bigg\rangle_{\mathcal{L}}\nonumber\\
    &=\langle w(\lambda),w(\mu)\rangle_{\mathcal{N}}+\langle u(\lambda),u(\mu)\rangle_{\mathcal{N}^{\perp}}-\langle w(\lambda),w(\mu)\rangle_{\mathcal{N}}-\langle I(\lambda)u(\lambda),I(\mu)u(\mu)\rangle_{\mathcal{N}^{\perp}}\nonumber\\
    &=\langle u(\lambda),u(\mu)\rangle_{\mathcal{N}^{\perp}}-\langle I(\lambda)u(\lambda),I(\mu)u(\mu)\rangle_{\mathcal{N}^{\perp}}\nonumber\\
    &=\langle (1_{\mathcal{N}^{\perp}}-I^{*}(\mu)I(\lambda))u(\lambda),u(\mu)\rangle_{\mathcal{N}^{\perp}}.\label{generalizedmodelfirstappeareance}
\end{align}
Hence $(\mathcal{N}^{\perp}, u, I)$ is an inner generalized model for $\varphi$. Let us take $\mathcal{M}=\mathcal{N}^{\perp}$, $u(\lambda)=P_{\mathcal{N}^{\perp}}v(\lambda)$ and $I(\lambda)$ as defined by equation (\ref{Ilambdaoperatorgeneralized}). Then $(\mathcal{M}, u, I)$ is an inner generalized model for $\varphi$.

 We now prove that $\tau$ is a $C$-point for the model $(\mathcal{N}^{\perp},u,I)$; that is, there exists a vector $u(\tau)$ such that, for every sequence $\{\lambda^{(n)}\}\in\mathbb{D}^{d}$ that converges to $\tau$ nontangentially, $u(\lambda^{(n)})\rightarrow u(\tau)$ in norm  as $n\rightarrow \infty$.
Once more, by Corollary \ref{coro5.7cara}, $\tau$ is a $B$-point for the model $(\mathcal{L},v)$ and $\gamma\in\mathrm{Ran}(1_{\mathcal{L}}-D\tau_{P})$. Let $u(\tau)$ be the unique element of smallest norm in the closed convex non-empty set 
$$\Big\{ x\in\mathcal{L}:(1_{\mathcal{L}}-D\tau_{P})x=\gamma\Big\}.$$
Then $u(\tau)\in\mathrm{Ker}(1_{\mathcal{L}}-D\tau_{P})^{\perp}$ and every element of $(1-D\tau_{P})^{-1}\gamma$ has the form $e\oplus u(\tau)$ for some $e\in\mathcal{N}$.

Let $X_{\tau}$ be the nontangential cluster set of $v(\lambda)$ at $\tau$ in the model $(\mathcal{L},v)$. By Proposition \ref{prop5.8cara}, since $\tau$ is a $B$-point for 
$\varphi$, a sequence $v(\lambda^{(n)})$ converges in norm as $n\rightarrow \infty$ if and only if it converges weakly in $\mathcal{L}$. With this, if $x\in X_{\tau}$ is the limit of the vector $v(\lambda^{(n)})$ as $\lambda^{(n)} {\overset{nt}\rightarrow}\tau$,  
$$(1_{\mathcal{L}}-D\lambda_{P}^{(n)})v(\lambda^{(n)})\rightarrow \begin{bmatrix}
    0 \\
    \gamma
\end{bmatrix},$$
then
$$ (1_{\mathcal{L}}-D\lambda_{P}^{(n)})v(\lambda^{(n)})\rightarrow (1_{\mathcal{L}}-D\tau_{P})x~~\text{as}~n\rightarrow \infty$$
and 
$$(1_{\mathcal{L}}-D\tau_{P})x
    =\begin{bmatrix}
    0 \\
    \gamma
\end{bmatrix}.$$
Thus 
$$x\in X_{\tau}\subset(1_{\mathcal{L}}-D\tau_{P})^{-1}\begin{bmatrix}
    0 \\
    \gamma
\end{bmatrix}\subset \bigg\{\begin{bmatrix}
    e \\
    u(\tau) \\
\end{bmatrix} : e\in\mathcal{N}\bigg\}.$$
We claim that $u(\lambda^{(n)})\rightarrow u(\tau)$ as $n\rightarrow \infty$ in norm for every sequence $\{\lambda^{(n)}\}$ such that $\lambda^{(n)} {\overset{nt}\rightarrow} \tau$. Suppose not, since $v(\lambda^{(n)})$ is bounded, then clearly $w(\lambda^{(n)})$ and $u(\lambda^{(n)})$ are bounded. By passing to a subsequence if necessary, we can assume that $u(\lambda^{(n)})\rightarrow y\neq u(\tau)$ as $n\rightarrow \infty$, and likewise, by passing to a further subsequence, we may suppose that $v(\lambda^{(n)})$ converges to some vector $ x\in X_{\tau}$ as $n\rightarrow \infty$. Then
$$v(\lambda^{(n)})=\begin{bmatrix}
    w(\lambda^{(n)}) \\
    u(\lambda^{(n)})
\end{bmatrix}\rightarrow x\in\bigg\{\begin{bmatrix}
    e \\
    u(\tau) \\
\end{bmatrix} : e\in\mathcal{N}\bigg\},$$
and therefore $u(\lambda^{(n)})\rightarrow u(\tau)$ as $n\rightarrow \infty$, which is a contradiction to our assumption. Hence, for every sequence $\{\lambda^{(n)}\}$ such that
 $\lambda^{(n)} {\overset{nt}\rightarrow} \tau$, $u(\lambda^{(n)})\rightarrow u(\tau)$ in norm as $n\rightarrow \infty$. Therefore $\tau$ is a $C$-point for the model $(\mathcal{N}^{\perp},u,I)$.

To prove $\tau$ is a carapoint for $I$, let $0<r<1$ and let $\lambda=r\tau$. Then
$$\lambda = (\lambda_{1},\hdots, \lambda_{d})=(r\tau_{1},\hdots,r\tau_{d}),$$
and 
$$\lambda\overline{\tau}=r\tau\overline{\tau}=(r\tau_{1}\overline{\tau}_{1},\hdots,r\tau_{d}\overline{\tau}_{d})=r\cdot\mathds{1}. $$
By definition of $I(\lambda)$, see equation (\ref{Ilambdaoperatorgeneralized}),
\begin{align*}
I(r\tau)&=1_{\mathcal{N^{\perp}}}-\bigg(\dfrac{\mathds{1}}{\mathds{1}-\overline{\tau}r\tau}\bigg)^{-1}_{Y}\\
        &=1_{\mathcal{N^{\perp}}}-\bigg(\dfrac{\mathds{1}}{\mathds{1}-r\cdot\mathds{1}}\bigg)^{-1}_{Y}\\
        &=1_{\mathcal{N^{\perp}}}-\Bigg[\dfrac{1}{1-r}\bigg(\sum_{j=1}^{d}Y_{j}\bigg)\Bigg]^{-1}\\
        &=1_{\mathcal{N^{\perp}}}-\Bigg[\dfrac{1}{1-r}\cdot 1_{\mathcal{N^{\perp}}}\Bigg]^{-1}=1_{\mathcal{N^{\perp}}}-(1-r)_{\mathcal{N^{\perp}}}=r\cdot 1_{\mathcal{N^{\perp}}}.
\end{align*}
Hence, 
$$\liminf_{\lambda\rightarrow \tau}\dfrac{1-\|I(\lambda)\|_{\mathcal{B(N^{\perp})}}}{1-\|\lambda\|_{\infty}}\leq \liminf_{r \rightarrow 1}\dfrac{1-\|I(r\tau)\|_{\mathcal{B(N^{\perp})}}}{1-\|r\cdot \mathds{1}\|_{\infty}}=\liminf_{r \rightarrow 1}\dfrac{1-r}{1-r}=1.$$
Therefore $\tau$ is a carapoint for $I$. We can also show $\lim_{\lambda{\overset{nt}\rightarrow}\tau}I(\lambda)=1_{\mathcal{M}}$, where $\mathcal{M}=\mathcal{N^{\perp}}$. Suppose that $\lambda{\overset{nt}\rightarrow}\tau$, then there exists a $c>0$ such that 
$$\|\lambda-\tau\|_{\infty}\leq c(1-\|\lambda\|_{\infty}).$$
By Theorem \ref{invertibleoperatorsforlambda}, 
\begin{align*}
    \|I(\lambda)-1_{\mathcal{M}}\|_{\mathcal{B(M)}} &= \bigg\|\bigg(\dfrac{\mathds{1}}{\mathds{1}-\overline{\tau}\lambda}\bigg)^{-1}_{Y}\bigg\|_{\mathcal{B(M)}}\\
    &\leq \max_{j=1,\hdots,d}\dfrac{\lvert 1-\overline{\tau_{j}}\lambda_{j}\rvert^{2}}{1-\mathrm{Re}(\overline{\tau_{j}}\lambda_{j})}\\
    &= \max_{j=1,\hdots,d}\dfrac{\lvert \tau_{j}-\lambda_{j}\rvert^{2}}{1-\mathrm{Re}(\overline{\tau_{j}}\lambda_{j})}.
\end{align*}
For $\tau_{j}\in\mathbb{T}$ and $\lambda_{j}\in\mathbb{D}$ for $j=1,\hdots,d$, 
$\mathrm{Re}(\overline{\tau_{j}}\lambda_{j}) \leq \lvert \overline{\tau_{j}}\lambda_{j}\rvert = \lvert \lambda_{j}\rvert.$
Thus,  for $j=1,\hdots,d$,
\begin{align*}
    1-\mathrm{Re}(\overline{\tau_{j}}\lambda_{j})&\geq 1-\lvert \lambda_{j} \rvert \\
    &\geq \min_{j=1,\hdots,d}(1-\lvert \lambda_{j}\rvert)\\
    &=1-\max_{j=1,\hdots,d}\lvert\lambda_{j}\rvert =1-\|\lambda\|_{\infty}.
\end{align*}
Hence, we have the following inequality
$$\|I(\lambda)-1_{\mathcal{M}}\|_{\mathcal{B(M)}}\leq \dfrac{\|\lambda-\tau\|^{2}_{\infty}}{1-\|\lambda\|_{\infty}}.$$
By assumption, $\lambda{\overset{nt}\rightarrow}\tau$. By Definition \ref{nontangetniallimitonDd}, there exists a $c>0$ such that $ \|\lambda - \tau \|_{\infty}\leq c ( 1 - \|\lambda\|_{\infty})$. Thus,
$$\|I(\lambda)-1_{\mathcal{M}}\|_{\mathcal{B(M)}}\leq \|\lambda-\tau\|_{\infty}\cdot\dfrac{\|\lambda-\tau\|_{\infty}}{1-\|\lambda\|_{\infty}}\leq  c \|\lambda-\tau\|_{\infty}.$$
By taking limits in the above equation, we find that $\lim_{\lambda{\overset{nt}\rightarrow}\tau}I(\lambda)=1_{\mathcal{M}}$.\end{proof}

The model $(\mathcal{M},u,I)$ constructed in Theorem \ref{generalizedmodeltheorem}
is often referred to as the {\em{desingularization}} of the model $(\mathcal{L},v)$ at $\tau$, relative to a particular realization of $(\mathcal{L},v)$, see the formal Definition \ref{desingularized}. 

\begin{corollary}\label{generalizedrealizationcorollary}
If $\tau\in\mathbb{T}^{d}$ is a carapoint for a function $\varphi\in\mathcal{SA}_{d}$, then $\varphi$ has a generalized realization
\begin{equation}\label{realizationforcarapoint}
\varphi(\lambda) = a+\langle I(\lambda)(1_{\mathcal{M}}-QI(\lambda))^{-1}\gamma,{\hat{\beta}}\rangle_{\mathcal{M}}~~\text{for all}~~\lambda\in\mathbb{D}^{d},
\end{equation}
for some separable Hilbert space $\mathcal{M}$,  scalar $a\in\mathbb{C}$, vectors ${\hat{\beta}}, \gamma\in\mathcal{M}$,  contraction $Q$ on $\mathcal{M}$ satisfying $\mathrm{Ker}(1_{\mathcal{M}}-Q)=\{0\}$ and  operator-valued analytic function $I:\mathbb{D}^{d}\rightarrow \mathcal{B(M)}$ such that $I$ is inner and
\begin{equation}\label{Ioperatoragainhere}
I(\lambda)=1_{\mathcal{M}}-\bigg(\dfrac{\mathds{1}}{\mathds{1}-\overline{\tau}\lambda}\bigg)^{-1}_{Y},
\end{equation}
for some $d$-tuple of positive contractions $Y=(Y_{1},\hdots,Y_{d})$ on $\mathcal{M}$ for which $\sum_{j=1}^{d}Y_{j}=1_{\mathcal{M}}$.
\end{corollary}
\begin{proof}
In the proof of Theorem \ref{generalizedmodeltheorem} we chose a  model $(\mathcal{L},v)$ for $\varphi\in\mathcal{SA}_{d}$ and any realization $(a,\beta,\gamma,D)$ of $(\mathcal{L},v)$ such that
\begin{equation}\label{realizationforphiinL}
\varphi(\lambda)=a+\langle\lambda_{P} v(\lambda),\beta\rangle_{\mathcal{L}},
\end{equation}
for all $\lambda\in\mathbb{D}^{d}$ where $\mathcal{L}$ is a separable Hilbert space with orthogonal decomposition $\mathcal{L} = \mathcal{L}_{1}\oplus \hdots \oplus \mathcal{L}_{d}$ and $\lambda_{P}=\lambda_{1}P_{1}+\hdots\lambda_{d}P_{d}$ such that $P=(P_{1},\hdots, P_{d})$ is the $d$-tuple of operators on $\mathcal{L}$ where $P_{j}:\mathcal{L}\rightarrow \mathcal{L}$ is the orthogonal projection onto $\mathcal{L}_{j}$ for $j=1,\hdots,d$.

 We defined $\mathcal{N}=\mathrm{Ker}(1_{\mathcal{L}}-D\tau_{P})$ and wrote $D\tau_{P}$ in the block matrix form relative to the orthogonal decomposition $\mathcal{L}=\mathcal{N}\oplus\mathcal{N}^{\perp}$,
\begin{equation}
    D\tau_{P} =\begin{bmatrix}
        1_{\mathcal{N}} & 0 \\
        0 & Q 
    \end{bmatrix}.
\end{equation}
It is clear from the definition of $Q$ that $Q$ is a contraction and $\mathrm{Ker}(1_{\mathcal{N}^{\perp}}-Q)=\{0\}$. 

Moreover, we proved that there exists a generalized model $(\mathcal{N}^{\perp},u,I)$ of $\varphi$, that is, $\mathcal{N}^{\perp}$ is a separable Hilbert space, $u:\mathbb{D}^{d}\rightarrow \mathcal{N}^{\perp}$ is an analytic mapping such that $u(\lambda)=P_{\mathcal{N}^{\perp}}v(\lambda)$ for all $\lambda\in\mathbb{D}^{d}$, $I$ is a contractive, analytic $\mathcal{B}(\mathcal{N}^{\perp})$-valued function on $\mathbb{D}^{d}$ such that $\|I(\lambda)\|_{\mathcal{N}^{\perp}}<1$ for all $\lambda\in\mathbb{D}^{d}$, $I$ is inner and $I$ is defined by equation (\ref{Ioperatoragainhere}). Thus, for all $\lambda,\mu\in\mathbb{D}^{d}$,
$$1-\overline{\varphi(\mu)}\varphi(\lambda) = \langle( 1_{\mathcal{N}^{\perp}}-I^{*}(\mu)I(\lambda))u(\lambda),u(\mu)\rangle_{\mathcal{N}^{\perp}}.$$
By equation (\ref{gammaequals1-qIu}), 
\begin{equation}\label{gammaequationagain}
\gamma=(1_{\mathcal{N}^{\perp}}-QI(\lambda))u(\lambda).
\end{equation}
By Lemma \ref{gammainranperp}, $\overline{\tau}_{P}\beta\in\mathcal{N}^{\perp}$. Define ${\hat{\beta}}=\overline{\tau}_{P}\beta$. Therefore by equations (\ref{realizationforphiinL}) and (\ref{gammaequationagain}),
\begin{align*}
\varphi(\lambda)    &= a+\langle \lambda_{P} v(\lambda),\beta\rangle_{\mathcal{L}} \\
                    &= a+\langle \overline{\tau}_{P}\lambda_{P} v(\lambda),\overline{\tau}_{P}\beta\rangle_{\mathcal{L}} \\
                    &= a+\langle P_{\mathcal{N}^{\perp}}\overline{\tau}_{P}\lambda_{P} v(\lambda),\overline{\tau}_{P}\beta\rangle_{\mathcal{N}^{\perp}} \\ 
                    &= a+\langle I(\lambda)u(\lambda),\overline{\tau}_{P}\beta\rangle_{\mathcal{N}^{\perp}}\\
                    &= a+\langle I(\lambda)u(\lambda),{\tilde{\beta}}\rangle_{\mathcal{N}^{\perp}}.
\end{align*}
Note that
$$\begin{bmatrix}
    a & 1\otimes{\hat{\beta}} \\
    \gamma\otimes 1 & Q \\
\end{bmatrix}\begin{bmatrix}
    1 \\
    I(\lambda)u(\lambda)
\end{bmatrix} = \begin{bmatrix}
    a+\langle I(\lambda)u(\lambda),{\hat{\beta}}\rangle_{\mathcal{N}^{\perp}} \\
    \gamma +QI(\lambda)u(\lambda)
\end{bmatrix} =\begin{bmatrix}
    \varphi(\lambda) \\
    u(\lambda)
\end{bmatrix}.$$
We show that the operator 
$$\begin{bmatrix}
    a & 1\otimes{\hat{\beta}} \\
    \gamma\otimes 1 & Q \\
\end{bmatrix}$$
is a unitary operator on $\mathbb{C}\oplus\mathcal{N}^{\perp}$. Indeed, define $P_{\mathbb{C}\oplus\{0\}\oplus\mathcal{N}^{\perp}}$ and $P_{\mathbb{C}\oplus\mathcal{N}^{\perp}}$ to be the orthogonal projections on $\mathbb{C}\oplus\mathcal{L}$ such that
$$P_{\mathbb{C}\oplus\{0\}\oplus\mathcal{N}^{\perp}} = \begin{bmatrix}
    1       & 0 \\
   0       & 0 \\
   0 & 1_{\mathcal{N}^{\perp}}\\
\end{bmatrix},~P_{\mathbb{C}\oplus\mathcal{N}^{\perp}} = \begin{bmatrix}
    1       & 0 & 0\\
   0       & 0 & 1_{\mathcal{N}^{\perp}}\\
\end{bmatrix},$$
where 
\begin{equation}\label{projectionsforthistheorem}
P^{*}_{\mathbb{C}\oplus\{0\}\oplus\mathcal{N}^{\perp}}P_{\mathbb{C}\oplus\{0\}\oplus\mathcal{N}^{\perp}} = 1_{\mathbb{C}\oplus\mathcal{N}^{\perp}}~\text{and}~~P^{*}_{\mathbb{C}\oplus\mathcal{N}^{\perp}}P_{\mathbb{C}\oplus\mathcal{N}^{\perp}} = 1_{\mathbb{C}\oplus\mathcal{L}}.
\end{equation}
By Lemma \ref{gammainranperp}, $\gamma\in\mathcal{N}^{\perp}$. Therefore the unitary operator
\begin{equation}\label{Lmatrix}
L=\begin{bmatrix}
    a       & 1\otimes \beta \\
    \gamma \otimes 1       & D 
\end{bmatrix}
\end{equation}
on $\mathbb{C}\oplus\mathcal{L}$ such that $\mathcal{L}=\mathcal{N}\oplus\mathcal{N}^{\perp}$ can be written as 
$$L = \begin{bmatrix}
    a & 1\otimes \beta|_{\mathcal{N}} & 1\otimes \beta|_{\mathcal{N}^{\perp}}\\
    0 & D_{11} & D_{12}\\
    \gamma\otimes 1 & D_{21} & D_{22}
\end{bmatrix}.$$
Thus
\begin{align*}
&P_{\mathbb{C}\oplus\mathcal{N}^{\perp}}\begin{bmatrix}
    a & 1\otimes \beta|_{\mathcal{N}} & 1\otimes \beta|_{\mathcal{N}^{\perp}}\\
    0 & D_{11} & D_{12}\\
    \gamma\otimes 1 & D_{21} & D_{22}
\end{bmatrix}\begin{bmatrix}
    1 & 0 & 0 \\
    0& \tau_{P}|_{\mathcal{N}} & 0 \\
    0 & 0 & \tau_{P}|_{\mathcal{N}^{\perp}}
\end{bmatrix}P_{\mathbb{C}\oplus\{0\}\oplus\mathcal{N}^{\perp}}\\
&=P_{\mathbb{C}\oplus\mathcal{N}^{\perp}}\begin{bmatrix}
    a & \tau_{P}\otimes \beta|_{\mathcal{N}} & \tau_{P}\otimes \beta|_{\mathcal{N}^{\perp}}\\
    0 & D_{11}\tau_{P} & D_{12}\tau_{P}\\
    \gamma\otimes 1 & D_{21}\tau_{P} & D_{22}\tau_{P}
\end{bmatrix}P_{\mathbb{C}\oplus\{0\}\oplus\mathcal{N}^{\perp}}\\
&=\begin{bmatrix}
    1       & 0 & 0\\
   0       & 0 & 1_{\mathcal{N}^{\perp}}\\
\end{bmatrix}\begin{bmatrix}
    a & 1\otimes \overline{\tau}_{P}\beta|_{\mathcal{N}} & 1\otimes \overline{\tau}_{P}\beta|_{\mathcal{N}^{\perp}}\\
    0 & D_{11}\tau_{P} & D_{12}\tau_{P}\\
    \gamma\otimes 1 & D_{21}\tau_{P} & D_{22}\tau_{P}
\end{bmatrix}\begin{bmatrix}
    1       & 0 \\
   0       & 0 \\
   0 & 1_{\mathcal{N}^{\perp}}\\
\end{bmatrix}.
\end{align*}
Recall that we have defined  ${\hat{\beta}}$ to be   $ \overline{\tau}_{P}\beta$ above.
By Lemma \ref{gammainranperp}, $\hat{\beta} \in\mathcal{N}^{\perp}$ and we can write the block matrix $D\tau_{P}$ as 
$$D\tau_{P} = \begin{bmatrix}
    1_{\mathcal{N}} & 0 \\
    0 & Q
\end{bmatrix},$$
 where $Q$ is a contraction and $\mathrm{Ker}(1-Q)=\{0\}$. Therefore, with ${\hat{\beta}} = \overline{\tau}_{P}\beta$, 
\begin{equation*}
    \begin{bmatrix}
    a & 1\otimes{\hat{\beta}}|_{\mathcal{N}} & 1\otimes {\hat{\beta}}|_{\mathcal{N}^{\perp}}\\
    0 & D_{11}\tau_{P} & D_{12}\tau_{P}\\
    \gamma\otimes 1 & D_{21}\tau_{P} & D_{22}\tau_{P}
\end{bmatrix} = \begin{bmatrix}
    a &  0 & 1\otimes {\hat{\beta}}\\
    0 & 1_{\mathcal{N}} & 0\\
    \gamma\otimes 1 &0 & Q
\end{bmatrix}.
\end{equation*}
Hence
\begin{align*}
P_{\mathbb{C}\oplus\mathcal{N}^{\perp}}&\begin{bmatrix}
    a       & 1\otimes \beta \\
    \gamma \otimes 1       & D 
\end{bmatrix}\begin{bmatrix}
    1_{\mathcal{N}} & 0 \\
    0& \tau_{P}
\end{bmatrix}P_{\mathbb{C}\oplus\{0\}\oplus\mathcal{N}^{\perp}}\\
&=\begin{bmatrix}
    1       & 0 & 0\\
   0       & 0 & 1\\
\end{bmatrix}\begin{bmatrix}
    a &  0 & 1\otimes {\hat{\beta}}\\
    0 & 1_{\mathcal{N}} & 0\\
    \gamma\otimes 1 &0 & Q
\end{bmatrix}\begin{bmatrix}
    1       & 0 \\
   0       & 0 \\
   0 & 1_{\mathcal{N}^{\perp}}\\
\end{bmatrix}\\
&=\begin{bmatrix}
    a & 1\otimes {\hat{\beta}} \\
    \gamma\otimes 1 & Q
\end{bmatrix}.
\end{align*}
Since $L$ from equation (\ref{Lmatrix}) and 
$\begin{bmatrix}
    1 & 0 \\
    0 & \tau_{P}
\end{bmatrix}$
are unitary on $\mathbb{C}\oplus\mathcal{L}$, by equation (\ref{projectionsforthistheorem}), 
\begin{align*}
    \begin{bmatrix}
    a & 1\otimes {\hat{\beta}} \\
    \gamma\otimes 1 & Q
\end{bmatrix}^{*}&\begin{bmatrix}
    a & 1\otimes{\hat{\beta}} \\
    \gamma\otimes 1 & Q
\end{bmatrix}\\
&=P^{*}_{\mathbb{C}\oplus\{0\}\oplus\mathcal{N}^{\perp}}L^{*}\begin{bmatrix}
    1 & 0 \\
    0& \overline{\tau}_{P}
\end{bmatrix}P^{*}_{\mathbb{C}\oplus\mathcal{N}^{\perp}}P_{\mathbb{C}\oplus\mathcal{N}^{\perp}}L\begin{bmatrix}
    1 & 0 \\
    0& \tau_{P}
\end{bmatrix}P_{\mathbb{C}\oplus\{0\}\oplus\mathcal{N}^{\perp}} \\
&= 1_{\mathbb{C}\oplus\mathcal{N}^{\perp}}
\end{align*}
and
\begin{align*}
    \begin{bmatrix}
    a & 1\otimes {\hat{\beta}} \\
    \gamma\otimes 1 & Q
\end{bmatrix}&\begin{bmatrix}
    a & 1\otimes {\hat{\beta}} \\
    \gamma\otimes 1 & Q
\end{bmatrix}^{*}\\
&=P_{\mathbb{C}\oplus\mathcal{N}^{\perp}}L\begin{bmatrix}
    1 & 0 \\
    0& \tau_{P}
\end{bmatrix}P_{\mathbb{C}\oplus\{0\}\oplus\mathcal{N}^{\perp}}P^{*}_{\mathbb{C}\oplus\{0\}\oplus\mathcal{N}^{\perp}}\begin{bmatrix}
    1 & 0 \\
    0& \overline{\tau}_{P}
\end{bmatrix}L^{*}P^{*}_{\mathbb{C}\oplus\mathcal{N}^{\perp}} \\
&= 1_{\mathbb{C}\oplus\mathcal{N}^{\perp}}.
\end{align*}
Thus 
$$\begin{bmatrix}
    a & 1\otimes {\hat{\beta}} \\
    \gamma\otimes 1 & Q
\end{bmatrix}$$
is unitary on $\mathbb{C}\oplus\mathcal{N}^{\perp}$. Hence $(a,{\hat{\beta}},\gamma,Q)$ is a generalized realization of the generalized model $(\mathcal{N}^{\perp},u,I)$ of $\varphi$. The statement of the theorem holds if we take $\mathcal{M}=\mathcal{N}^{\perp}$.\end{proof}

\begin{remark}\label{case-d=2} \rm Let us consider the case when $d=2$. Theorem \ref{card2} states that $\ph\in\mathcal{S}(\D^2)$ if and only if $\ph\in\mathcal{SA}_2$. Let us compare the known results for $\ph\in\mathcal{S}(\D^2)$ from \cite{ATY2} with our Theorem \ref{generalizedmodeltheorem}. In \cite[Theorem 3.6]{ATY2} the authors considered a function 
$\ph\in\mathcal{S}(\D^2)$ with a carapoint $\tau \in \T^2$  and used a model $(\mathcal{L},v)$ of $\varphi$ and a realization $(\alpha,\beta, \gamma, D)$ of $(\mathcal{L},v)$ to construct an inner generalized model $(\m,u,I_{ATY})$, where
\be\label{I1}
I_{ATY}(\lam) = \frac{\bar\tau_1\lam_1 Y_1+\bar\tau_2\lam_2 Y_2 - \bar\tau_1 \bar\tau_2 \lam_1\lam_2 1_\m}{1_\m-\bar\tau_1\lam_1 Y_2 - \bar\tau_2\lam_2 Y_1} \ \text{for}\ \lam\in\D^2,
\ee
for some pair of positive commuting contractions $Y_1,Y_2$ on $\m$ such that $Y_1+Y_2=1_\m$. In
Theorem \ref{generalizedmodeltheorem}, for $\ph\in\mathcal{SA}_2$ with a carapoint $\tau \in \T^2$, we can use the same model $(\mathcal{L},v)$ of $\varphi$ and the same realization $(\alpha,\beta, \gamma, D)$ of $(\mathcal{L},v)$ to prove
the existence of an inner generalized model $(\m,u,I)$ where $\m$ and $u$ are the same as before and  the $\mathcal{B}(\m)$-valued inner function $I$  is given by
\be\label{I2}
I(\lam) = 1_\m - \left(\frac{\1}{\1-\bar\tau \lam}\right)_Y^{-1} \ \text{for}\ \lam\in\D^2,
\ee
again for some pair $Y=(Y_1,Y_2)$ of positive contractions summing to $1_\m$.
We claim that $I=I_{ATY}$. The operators $Y_1$ and $Y_2=1-Y_1$ commute, and it is simple to verify that, for all $\lam\in\D^2 $, 
\be\label{I3}
I_{ATY}(\lam) = 1_\m - \frac{(1-\bar\tau_1\lam_1)(1-\bar\tau_2\lam_2)}{1_\m-\bar\tau_1\lam_1 Y_2 - \bar\tau_2\lam_2 Y_1}.
\ee
Therefore, for all $\lam\in\D^2 $,
\begin{align*}
I(\lam)=I_{ATY}(\lam) &\iff \left(\frac{\1}{\1-\bar\tau\lam}\right)_Y^{-1} = \frac{(1-\bar\tau_1\lam_1)(1-\bar\tau_2\lam_2)}{1_\m-\bar\tau_1\lam_1 Y_2 - \bar\tau_2\lam_2 Y_1}\ \\
&\iff 1_\m = \frac{(1-\bar\tau_1\lam_1)(1-\bar\tau_2\lam_2)}{1_\m-\bar\tau_1\lam_1 Y_2 - \bar\tau_2\lam_2 Y_1}\left(\frac{\1}{\1-\bar\tau\lam}\right)_Y\\
	&\iff 1_\m = \frac{(1-\bar\tau_1\lam_1)(1-\bar\tau_2\lam_2)}{1_\m-\bar\tau_1\lam_1 Y_2 - \bar\tau_2\lam_2 Y_1} \left(\frac{1}{1-\bar\tau_1\lam_1}Y_1 +\frac{1}{1-\bar\tau_2\lam_2}Y_2\right)\\
	&\iff  1_\m = \frac{1}{1_\m-\bar\tau_1\lam_1 Y_2 - \bar\tau_2\lam_2 Y_1} \left(
	(1-\bar\tau_2\lam_2)Y_1 +(1-\bar\tau_1\lam_1)Y_2\right).
\end{align*}
Since $Y_1+Y_2=1_\m$, the last equation holds.  Hence $I=I_{ATY}$ and the two generalized models coincide.
\end{remark}

\section{Directional derivatives and slope functions}\label{directional_der}

Let $\varphi\in\mathcal{SA}_{d}$ and let $\tau\in\mathbb{T}^{d}$ be a carapoint for $\varphi$. Let us describe a set of directions $\delta\in\mathbb{C}^{d}$ which  point into $\mathbb{D}^{d}$ at $\tau$.
For $\tau\in\mathbb{C}^{d}$, define
$$\mathbb{C}_{+}^{d}(\tau) = \big\{( z_{1},\hdots, z_{d})\in\mathbb{C}^{d} : \mathrm{Re}( z_{j}\overline{\tau}_{j})>0~~\text{for}~~j=1,\hdots,d\big\},$$
and 
$$\mathbb{C}_{+}(\tau_{j}) = \big\{ z\in\mathbb{C} : \mathrm{Re}(z\overline{\tau}_{j})>0\big\}~~\text{for}~j=1,\hdots,d.$$
Observe that $\mathbb{C}^{d}_{+}(\tau) = \mathbb{C}_{+}(\tau_{1})\times \hdots \times\mathbb{C}_{+}(\tau_{d})$. Note that for $t>0$ and $\delta\in\mathbb{C}^{d}$,
$$\lvert \tau_{j} + t\delta_{j}\rvert^{2} = 1+2t\mathrm{Re}(\overline{\tau}_{j}\delta_{j})+t^{2}\lvert \delta_{j}\rvert^{2}.$$
Therefore $\delta$ points into $\mathbb{D}^{d}$ at $\tau$ if and only if $\mathrm{Re}(\overline{\tau}_{j}\delta_{j})<0$ for $j=1,\hdots,d$, that is, if and only if $-\delta \in\mathbb{C}^{d}_{+}(\tau)$. Note that $\tau\in\mathbb{C}_{+}^{d}(\tau)$ and, if $ z\in\mathbb{C}_{+}^{d}(\tau)$, then $ \overline{z}\in\mathbb{C}^{d}_{+}(\overline{\tau})$, since
$$0<\mathrm{Re}(\overline{\tau_{j}} z_{j})=\mathrm{Re}(\tau_{j}\overline{ z_{j}}).$$ 

 Let $\varphi\in\mathcal{SA}_{d}$ and let $\tau\in\mathbb{T}^{d}$ be a carapoint for the function $\varphi$. Let $\delta\in\mathbb{C}^{d}$ point into $\mathbb{D}^{d}$ at $\tau$. $\varphi$ is differentiable in the direction $\delta\in\mathbb{C}^{d}$ at the point $\tau$ if there exists a complex number $\omega$ of unit modulus such that, whenever $\tau+t \delta\in\mathbb{D}^{d}$ for sufficiently small positive $t$, 
$$D_{\delta}\varphi(\tau)=\lim_{t\rightarrow 0^{+}} \dfrac{\varphi(\tau+t\delta)-\omega}{t}~~{\text{exists}}.$$
\begin{theorem}\label{directionalderivtheorem}
  Let $\tau\in\mathbb{T}^{d}$ be a carapoint for $\varphi\in\mathcal{SA}_{d}$. 
Then  $\varphi$ is directionally differentiable at $\tau$ and 
there exists an analytic function $h:\mathbb{C}_{+}^{d}(\tau)\rightarrow\mathbb{C}$ such that
    \begin{enumerate}
        \item $\mathrm{Re}(-h(z))>0$ for all $ z\in\mathbb{C}_{+}^{d}(\tau)$;
        \item \begin{equation}\label{hfunction}
            h(\tau) = -\liminf_{\lambda{\rightarrow}\tau} \dfrac{1-\lvert \varphi(\lambda)\rvert}{1-\|\lambda \|_{\infty}};
        \end{equation}
        \item for all $ \delta\in\mathbb{C}_{+}^{d}(\tau)$, 
        \begin{equation}\label{directionalderivativeD}
        D_{- \delta}\varphi(\tau) = \varphi(\tau)h( \delta),
\end{equation}
where $ \varphi(\tau) = \lim_{\lambda{\overset{nt}\rightarrow}\tau}\varphi(\lambda)$ and $|\varphi(\tau)|=1$.
    \end{enumerate}
Moreover, we can express the function $h:\mathbb{C}_{+}^{d}(\tau)\rightarrow\mathbb{C}$ by
$$ h( z) = -\bigg\langle \bigg(\dfrac{\mathds{1}}{\overline{\tau} z}\bigg)_{Y}^{-1}u(\tau),u(\tau)\bigg\rangle_{\mathcal{M}},~~\text{for all}~ z\in\mathbb{C}_{+}^{d}(\tau),$$  
for some desingularized model $(\mathcal{M},u,I)$ of $\varphi$ with $I$ of the form
\begin{equation}\label{Ioperatorgeneralizedmodel-2}
I(\lambda) = 1_{\mathcal{M}} - \bigg( \dfrac{\mathds{1}}{\mathds{1}-\overline{\tau}\lambda}\bigg)_{Y}^{-1} \
\text{for} \ \lambda \in \mathbb{D}^{d},
\end{equation}
for some $d$-tuple of positive contractions $(Y_{1},\hdots,Y_{d})$ on $\mathcal{M}$ such that $\sum_{j=1}^{d}Y_{j}=1_{\mathcal{M}}$. 
\end{theorem}
\begin{proof}
By Theorem \ref{generalizedmodeltheorem}, there exists a desingularized model $(\mathcal{M},u,I)$ of $\varphi$ relative to $\tau$, with
$$I(\lambda)=1_{\mathcal{M}}-\bigg(\dfrac{\mathds{1}}{\mathds{1}-\overline{\tau}\lambda}\bigg)^{-1}_{Y}, \; \text{ for all } \; \lambda \in \D^d,$$
and for some $d$-tuple of positive contractions $(Y_{1},\hdots,Y_{d})$ on $\mathcal{M}$ such that $\sum_{j=1}^{d}Y_{j}=1_{\mathcal{M}}.$
By Theorem \ref{generalizedmodeltheorem}, $\tau$ is a $C$-point of the model $(\mathcal{M},u,I)$, so there exists a vector $u(\tau)\in\mathcal{M}$ such that
\begin{equation}\label{cpointulambda}
\lim_{\lambda{\overset{nt}\rightarrow}\tau}u(\lambda)=u(\tau)
\end{equation}
and, for all $\lambda,\mu\in\mathbb{D}^{d}$, 
\begin{equation}\label{modelequationinthisproof-1}
    1-\overline{\varphi(\mu)}\varphi(\lambda) = \bigg\langle \Big(1_{\mathcal{M}}-I^{*}(\mu)I(\lambda)\Big)u(\lambda),u(\mu)\bigg\rangle_{\mathcal{M}}.
\end{equation}
By Proposition  \ref{phi-nt-omega},   $\ph$ has a non-tangential limit at $\tau$,
$$ \lim_{\lambda{\overset{nt}\rightarrow}\tau}\varphi(\lambda)=\omega,$$
and  $|\omega|=1$. We shall define $\varphi(\tau)$ to be equal to $ \omega$.

Take limits in the last equation as ${\mu{\overset{nt}\rightarrow}\tau}$ to obtain
\begin{equation}\label{generalizedmodelslopes}
    1-\overline{\varphi(\tau)}\varphi(\lambda)=\bigg\langle\Big( 1_{\mathcal{M}}-I(\lambda)\Big)u(\lambda),u(\tau)\bigg\rangle_{\mathcal{M}}.
\end{equation}
On multiplying equation (\ref{generalizedmodelslopes}) by $-\varphi(\tau)$, we have 
\begin{align}
   \varphi(\lambda)-\varphi(\tau)&=-\varphi(\tau)\bigg\langle\Big( 1_{\mathcal{M}}-I(\lambda)\Big)u(\lambda),u(\tau)\bigg\rangle_{\mathcal{M}}\nonumber\\
&=\varphi(\tau)\bigg\langle\Big(I(\lambda)-1_{\mathcal{M}}\Big)u(\lambda),u(\tau)\bigg\rangle_{\mathcal{M}}\nonumber\\
&=\varphi(\tau)\bigg\langle\Big(I(\lambda)-1_{\mathcal{M}}\Big)u(\tau),u(\tau)\bigg\rangle_{\mathcal{M}}\nonumber\\
&\hspace{2cm}+\varphi(\tau)\bigg\langle\Big(I(\lambda)-1_{\mathcal{M}}\Big)\big(u(\lambda)-u(\tau)\big),u(\tau)\bigg\rangle_{\mathcal{M}}.\label{lastequationontheleft}\\\nonumber
\end{align}

Let $ \delta\in\mathbb{C}_{+}^{d}(\tau)$ so that $\lambda_{t}=\tau-t\delta\in\mathbb{D}^{d}$ for small enough $t>0$. By Corollary \ref{invertibleoperatorsforlambdacorollary}, $I(\lambda_{t})$ is well defined since $\bigg(\dfrac{\mathds{1}}{\mathds{1}-\overline{\tau}\lambda_{t}}\bigg)_{Y}^{-1}$ exists. Note that, for small enough $t$, $\mathrm{Re}(1-t\overline{\tau_{j}}\delta_{j})<1$, that is, $\mathrm{Re}(-t\overline{\tau_{j}}\delta_{j})<0$. Therefore
\begin{align}
    I(\lambda_{t})-1_{\mathcal{M}} &=1_{\mathcal{M}}-\bigg(\dfrac{\mathds{1}}{\mathds{1}-\overline{\tau}\lambda_{t}}\bigg)_{Y}^{-1}-1_{\mathcal{M}} \nonumber\\
    &=-\bigg(\sum_{j=1}^{d}\dfrac{1}{1-\overline{\tau_{j}}(\tau_{j}-t \delta_{j})}Y_{j}\bigg)^{-1}\nonumber\\
    &=-\bigg(\sum_{j=1}^{d}\dfrac{1}{t\overline{\tau_{j}} \delta_{j}}Y_{j}\bigg)^{-1}
    =-t\bigg(\dfrac{\mathds{1}}{\overline{\tau} \delta}\bigg)_{Y}^{-1}\label{realparttauz}.
\end{align}
By Theorem \ref{invertibleoperatorsforlambda},  since $\mathrm{Re}(\overline{\tau_{j}} \delta_{j})>0$ for $j=1,\hdots,d$,
equation (\ref{realparttauz}) produces a well defined operator on $\mathcal{B}(\mathcal{M})$. Thus equation (\ref{lastequationontheleft}) in combination with equation (\ref{realparttauz}) gives 
\begin{align*}
\dfrac{\varphi(\lambda_{t})-\varphi(\tau)}{t}&=\varphi(\tau)\bigg\langle \dfrac{I(\lambda_{t})-1_{\mathcal{M}}}{t} u(\tau), u(\tau) \bigg \rangle_{\mathcal{M}}\\ 
&\hspace{2cm}+\varphi(\tau)\bigg\langle \dfrac{I(\lambda_{t})-1_{\mathcal{M}}}{t}\big(u(\lambda_{t})-u(\tau)\big),u(\tau)\bigg\rangle_{\mathcal{M}}\\
&=-\varphi(\tau)\bigg\langle\bigg(\dfrac{\mathds{1}}{\overline{\tau} \delta}\bigg)_{Y}^{-1}u(\tau),u(\tau)\bigg\rangle_{\mathcal{M}}\nonumber\\
&\hspace{2cm}-\varphi(\tau)\bigg\langle \bigg(\dfrac{\mathds{1}}{\overline{\tau} \delta}\bigg)_{Y}^{-1}\big(u(\lambda_{t})-u(\tau)\big),u(\tau)\bigg\rangle_{\mathcal{M}}.
\end{align*}
By equation (\ref{cpointulambda}), $\lim_{\lambda{\overset{nt}{\rightarrow}} \tau}u(\lambda)=u(\tau)$, thus by letting $t\rightarrow 0^{+}$ in the above equation we get that, for all $\delta\in\mathbb{C}_{+}^{d}(\tau)$, 
\begin{equation}\label{directionalderiv123}
    D_{-\delta}\varphi(\tau)=\lim_{t\rightarrow 0^{+}} \dfrac{\varphi(\lambda_{t})-\varphi(\tau)}{t} =-\varphi(\tau)\bigg\langle \bigg(\dfrac{\mathds{1}}{\overline{\tau} \delta}\bigg)_{Y}^{-1}u(\tau),u(\tau)\bigg\rangle_{\mathcal{M}}.
\end{equation}
Thus $\varphi$ is directionally differentiable at $\tau$.

Let us define the function $h:\mathbb{C}_{+}^{d}(\tau)\rightarrow\mathbb{C}$ by
\begin{equation}\label{Def-h}
 h( z) = -\bigg\langle \bigg(\dfrac{\mathds{1}}{\overline{\tau} z}\bigg)_{Y}^{-1}u(\tau),u(\tau)\bigg\rangle_{\mathcal{M}}~~\text{for all}~ z\in\mathbb{C}_{+}^{d}(\tau).
 \end{equation}
Then
$$D_{-\delta}\varphi(\tau) =\varphi(\tau)h(\delta),$$
which proves equation (\ref{directionalderivativeD}).
Let us prove that $h$ is analytic at every point $z\in\mathbb{C}^{d}_{+}(\tau)$. By Proposition \ref{propforIbeinganalytic} statement (3) with $\lambda = \overline{\tau} z$, for all $ z,x\in\mathbb{C}^{d}_{+}(\tau)$,
$$\Bigg(\bigg(\dfrac{\mathds{1}}{\overline{\tau} z}\bigg)_{Y}^{-1}\Bigg)'(x) = \bigg(\dfrac{\mathds{1}}{\overline{\tau} z}\bigg)_{Y}^{-1}\bigg(\dfrac{x}{\overline{\tau} z^{2}}\bigg)_{Y}\bigg(\dfrac{\mathds{1}}{\overline{\tau} z}\bigg)_{Y}^{-1}.$$
Recall that $h$ is defined by equation (\ref{Def-h}), and so, one deduces that
$$h'( z)(x) = -\bigg\langle \Bigg(\bigg(\dfrac{\mathds{1}}{\overline{\tau} z}\bigg)_{Y}^{-1}\Bigg)'(x)u(\tau),u(\tau)\bigg\rangle_{\mathcal{M}}, \; \text{for all} \; x\in\mathbb{C}^{d}_{+}(\tau).$$
Hence $h$ is analytic on $\mathbb{C}_{+}^{d}(\tau)$ and
$$h'( z)(x) = -\bigg\langle \bigg(\dfrac{\mathds{1}}{\overline{\tau} z}\bigg)_{Y}^{-1}\bigg(\dfrac{x}{\overline{\tau} z^{2}}\bigg)_{Y}\bigg(\dfrac{\mathds{1}}{\overline{\tau} z}\bigg)_{Y}^{-1}u(\tau),u(\tau)\bigg\rangle_{\mathcal{M}}.$$
Let us prove that 
\begin{equation*}
            h(\tau) = -\liminf_{\lambda{\rightarrow}\tau} \dfrac{1-\lvert \varphi(\lambda)\rvert}{1-\|\lambda \|_{\infty}}.
\end{equation*}
Note that, for $\tau\in\mathbb{C}_{+}^{d}(\tau)$, 
$$\bigg(\dfrac{\mathds{1}}{\overline{\tau}\tau}\bigg)^{-1}_{Y} = \bigg(\sum_{j=1}^{d}Y_{j} \bigg)^{-1}= 1_{\mathcal{M}}.$$
Thus
\begin{equation}\label{hequationintheproofequalsnormu}
h(\tau)=-\bigg\langle \bigg(\dfrac{\mathds{1}}{\overline{\tau}\tau}\bigg)^{-1}_{Y}u(\tau),u(\tau)\bigg\rangle_{\mathcal{M}}=-\|u(\tau)\|^{2}_{\mathcal{M}}.
\end{equation}
From the model equation (\ref{modelequationinthisproof-1}), for any $\lambda\in\mathbb{D}^{d}$, 
\begin{equation}\label{equationdoublestar}
    1-\lvert \varphi(\lambda) \rvert^{2} = \|u(\lambda)\|^{2}_{\mathcal{M}}-\|I(\lambda)u(\lambda)\|^{2}_{\mathcal{M}}.
\end{equation}
Let $\lambda_{t}=\tau-t\tau$ for $t>0$. Note that
\begin{align}
    I(\lambda_{t}) &= 1_{\mathcal{M}}-\bigg(\dfrac{\mathds{1}}{\mathds{1}-\overline{\tau}\lambda_{t}}\bigg)_{Y}^{-1} \nonumber\\
    &= 1_{\mathcal{M}}-\bigg( \sum_{j=1}^{d} \dfrac{1}{1-\overline{\tau_{j}}(\tau_{j}-t\tau_{j})}Y_{j}\bigg)^{-1}\nonumber\\
    &= 1_{\mathcal{M}}-\bigg( \sum_{j=1}^{d} \dfrac{1}{t}Y_{j}\bigg)^{-1}\nonumber\\
    &= 1_{\mathcal{M}}-t\bigg( \sum_{j=1}^{d} Y_{j}\bigg)^{-1}
    = (1-t)\cdot1_{\mathcal{M}}.\label{1minust}
\end{align}
For small enough $t>0$, by equations (\ref{equationdoublestar}) and (\ref{1minust}),
\begin{align*}
    1-\lvert \varphi(\lambda_{t}) \rvert^{2} &= \|u(\lambda_{t})\|_{\mathcal{M}}^{2} - \|I(\lambda_{t})u(\lambda_{t})\|_{\mathcal{M}}^{2} \\
    &= \|u(\lambda_{t})\|_{\mathcal{M}}^{2}-\|(1-t)u(\lambda_{t})\|^{2}_{\mathcal{M}}\\
    &= (2t-t^{2})\|u(\lambda_{t})\|^{2}_{\mathcal{M}}.
\end{align*}
This implies that
\begin{equation}\label{7starstar}
    \lvert \varphi(\lambda_{t})\rvert^{2}=1-(2t-t^{2})\|u(\lambda_{t})\|^{2}_{\mathcal{M}}.
\end{equation}
For $\lambda_{t}=\tau-t\tau$, we also have
$$\|\lambda_{t}\|_{\infty} = \|\tau-t\tau\|_{\infty} = (1-t)\|\tau\|_{\infty} = 1-t,$$
and so $1-\|\lambda_{t}\|^{2}_{\infty} = 2t-t^{2}>0$ for small $t$. Hence 
$$\lim_{t\rightarrow 0^{+}}\dfrac{1-\lvert \varphi(\lambda_{t})\rvert^{2}}{1-\|\lambda_{t}\|^{2}_{\infty}} = \|u(\tau)\|^{2}_{\mathcal{M}}.$$
Note that equation (\ref{7starstar}) implies that
$$\lvert \varphi(\lambda_{t})\rvert = \big(1-(2t-t^{2})\|u(\lambda_{t})\|^{2}_{\mathcal{M}}\big)^{\frac{1}{2}}.$$
Thus we have
\begin{align*}
    \lim_{t\rightarrow 0^{+}}\dfrac{1-\lvert \varphi(\lambda_{t}) \rvert}{1-\|\lambda_{t}\|_{\infty}}&=\lim_{t\rightarrow 0^{+}}\dfrac{\big(1-\lvert \varphi(\lambda_{t})\rvert^{2}\big)}{\big(1-\|\lambda_{t}\|^{2}_{\infty}\big)} \dfrac{\big(1+\|\lambda_{t}\|_{\infty}\big)}{\big(1+\lvert \varphi(\lambda_{t})\rvert\big)}\\
    &=\lim_{t\rightarrow 0^{+}}\dfrac{\big(1-\lvert \varphi(\lambda_{t})\rvert^{2}\big)}{\big(1-\|\lambda_{t}\|^{2}_{\infty}\big)}\dfrac{\big(2-t\big)}{\big(1+\big(1-(2t-t^{2})\|u(\lambda_{t})\|^{2}_{\mathcal{M}}\big)^{\frac{1}{2}}\big)}\\
    &=\lim_{t\rightarrow 0^{+}}\dfrac{1-\lvert \varphi(\lambda_{t})\rvert^{2}}{1-\|\lambda_{t}\|^{2}_{\infty}}\lim_{t\rightarrow 0^{+}}\dfrac{2-t}{1+\big(1-(2t-t^{2})\|u(\lambda_{t})\|^{2}_{\mathcal{M}}\big)^{\frac{1}{2}}}\\
    &=\lim_{t\rightarrow 0^{+}}\dfrac{1-\lvert \varphi(\lambda_{t})\rvert^{2}}{1-\|\lambda_{t}\|^{2}_{\infty}}\\
    &= \lim_{t\rightarrow 0^{+}} \|u(\lambda_{t})\|^{2}_{\mathcal{M}}
    =\|u(\tau)\|_{\mathcal{M}}^{2}.
\end{align*}
Hence, by equation (\ref{hequationintheproofequalsnormu}), 
$$h(\tau)=-\lim_{t\rightarrow 0^{+}}\dfrac{1-\lvert \varphi(\lambda_{t}) \rvert}{1-\|\lambda_{t}\|_{\infty}}.$$
Moreover, since $\tau$ is a carapoint of $\varphi$, 
$$\lim_{t\rightarrow 0^{+}}\dfrac{1-\lvert \varphi(\lambda_{t}) \rvert}{1-\|\lambda_{t}\|_{\infty}} = \liminf_{\lambda\rightarrow \tau}\dfrac{1-\lvert \varphi(\lambda) \rvert}{1-\|\lambda\|_{\infty}}.$$
Thus equation (\ref{hfunction}) follows.
It remains to show that $\mathrm{Re}(h( z))<0$. For all $ z\in\mathbb{C}^{d}_{+}(\tau)$,
\begin{align}
    \mathrm{Re}(-h( z)) &= \dfrac{1}{2}\bigg[\bigg\langle \bigg(\dfrac{\mathds{1}}{\overline{\tau} z}\bigg)_{Y}^{-1} u(\tau),u(\tau)\bigg\rangle_{\mathcal{M}}+\overline{\bigg\langle \bigg(\dfrac{\mathds{1}}{\overline{\tau} z}\bigg)_{Y}^{-1} u(\tau),u(\tau)\bigg\rangle_{\mathcal{M}}}\bigg]\nonumber\\
    &= \dfrac{1}{2}\bigg[\bigg\langle \bigg(\dfrac{\mathds{1}}{\overline{\tau} z}\bigg)_{Y}^{-1} u(\tau),u(\tau)\bigg\rangle_{\mathcal{M}}+\bigg\langle u(\tau),\bigg(\dfrac{\mathds{1}}{\overline{\tau} z}\bigg)_{Y}^{-1} u(\tau)\bigg\rangle_{\mathcal{M}}\bigg]\nonumber\\
    &= \dfrac{1}{2}\bigg[\bigg\langle \bigg(\dfrac{\mathds{1}}{\overline{\tau} z}\bigg)_{Y}^{-1} u(\tau),u(\tau)\bigg\rangle_{\mathcal{M}}+\bigg\langle \bigg(\dfrac{\mathds{1}}{\tau \overline{z}}\bigg)_{Y}^{-1}u(\tau), u(\tau)\bigg\rangle_{\mathcal{M}}\bigg]\nonumber\\
    &= \bigg\langle \dfrac{1}{2}\bigg[\bigg(\dfrac{\mathds{1}}{\overline{\tau} z}\bigg)_{Y}^{-1}+\bigg(\dfrac{\mathds{1}}{\tau \overline{z}}\bigg)_{Y}^{-1}\bigg] u(\tau),u(\tau)\bigg\rangle_{\mathcal{M}}\nonumber\\
    &= \bigg\langle \mathrm{Re}\bigg(\bigg(\dfrac{\mathds{1}}{\overline{\tau} z}\bigg)_{Y}^{-1}\bigg) u(\tau),u(\tau)\bigg\rangle_{\mathcal{M}}.\label{realpartoftauzstar}
\end{align}
Furthermore, for all $ z\in\mathbb{C}^{d}_{+}(\tau)$,
\begin{align*}
    \mathrm{Re}\bigg(\bigg(\dfrac{\mathds{1}}{\overline{\tau} z}\bigg)_{Y}^{-1}\bigg) &= \dfrac{1}{2}\bigg[\bigg(\dfrac{\mathds{1}}{\overline{\tau} z}\bigg)_{Y}^{-1}+\bigg(\dfrac{\mathds{1}}{\tau \overline{z}}\bigg)_{Y}^{-1}\bigg]\\
    &= \dfrac{1}{2}\bigg[\bigg(\dfrac{\mathds{1}}{\tau \overline{z}}\bigg)_{Y}^{-1}\bigg(\bigg(\dfrac{\mathds{1}}{\tau \overline{z}}\bigg)_{Y}\bigg(\dfrac{\mathds{1}}{\overline{\tau} z}\bigg)_{Y}^{-1}+1_{\mathcal{M}}\bigg)\bigg]\\
    &= \dfrac{1}{2}\bigg[\bigg(\dfrac{\mathds{1}}{\tau \overline{z}}\bigg)_{Y}^{-1}\bigg(\bigg(\dfrac{\mathds{1}}{\tau \overline{z}}\bigg)_{Y}+\bigg(\dfrac{\mathds{1}}{\overline{\tau} z}\bigg)_{Y}\bigg)\bigg(\dfrac{\mathds{1}}{\overline{\tau} z}\bigg)_{Y}^{-1}\bigg]\\
    &= \bigg(\dfrac{\mathds{1}}{\tau \overline{z}}\bigg)_{Y}^{-1}\mathrm{Re}\bigg(\bigg(\dfrac{\mathds{1}}{\overline{\tau} z}\bigg)_{Y}\bigg)\bigg(\dfrac{\mathds{1}}{\overline{\tau} z}\bigg)_{Y}^{-1}.
\end{align*}

It remains to show that $\mathrm{Re}\bigg(\bigg(\dfrac{\mathds{1}}{\overline{\tau} z}\bigg)_{Y}\bigg)$ is positive definite. By Lemma \ref{projectiononcdlemma}, $Y_{j}\in\mathcal{B}(\mathcal{M})$ satisfies $0\leq Y_{j}\leq 1$ for $j=1,\hdots,d$. Thus, for all $z\in\mathbb{C}^{d}_{+}(\tau)$, 
\begin{align*}
    \mathrm{Re}\bigg(\bigg(\dfrac{\mathds{1}}{\overline{\tau}z}\bigg)_{Y}\bigg) &= \dfrac{1}{2}\bigg[\bigg(\dfrac{\mathds{1}}{\overline{\tau}z}\bigg)_{Y} + \bigg(\dfrac{\mathds{1}}{\tau \overline{z}}\bigg)_{Y}\Big]\\
    &= \dfrac{1}{2}\sum_{j=1}^{d}\dfrac{1}{\overline{\tau_{j}}z_{j}}Y_{j}+\sum_{j=1}^{d}\dfrac{1}{\tau_{j}\overline{z_{j}}}Y_{j} \\
    &= \dfrac{1}{2}\sum_{j=1}^{d}\dfrac{\overline{\tau_{j}}z_{j}+\tau_{j}\overline{z_{j}}}{\lvert z_{j}\rvert^{2}}Y_{j} 
    = \sum_{j=1}^{d}\dfrac{\mathrm{Re}(\overline{\tau_{j}}z_{j})}{\lvert z_{j}\rvert^{2}}Y_{j}.
\end{align*}
Since $ z\in\mathbb{C}^{d}_{+}(\tau)$, $\mathrm{Re}(\overline{\tau_{j}} z_{j})>0$. Hence $\mathrm{Re}\bigg(\bigg(\dfrac{\mathds{1}}{\overline{\tau} z}\bigg)_{Y}\bigg)$ is the sum of positive definite operators and in turn $\mathrm{Re}\bigg(\bigg(\dfrac{\mathds{1}}{\overline{\tau} z}\bigg)^{-1}_{Y}\bigg)> 0$. Thus $\mathrm{Re}(-h(z))$ is strictly positive for all $ z\in\mathbb{C}_{+}^{d}(\tau)$.\end{proof} 

We shall call the function $h$ described in Theorem \ref{directionalderivtheorem}
 the {\em slope function} of $\ph$ at the point $\tau$.  Thus $h$ is the slope function of $\ph$ at a carapoint $\tau\in\T^d$  
if, for all $ \delta\in\mathbb{C}_{+}^{d}(\tau)$,
            \begin{equation}\label{directionalderivativeD-2}
        D_{- \delta}\varphi(\tau) = \varphi(\tau)h( \delta).
\end{equation}

\begin{remark} Let  $\varphi\in\mathcal{SA}_{d}$ and let $\tau\in\mathbb{T}^{d}$. It is clear that the existence of a slope function $h$  of $\ph$ at $\tau$ implies that 
\be \label{Carapoint-2}
\liminf_{\lam \to \tau, \lam \in \D^d} \frac{1-|\ph(\lam)|}{1-\|\lam\|_\infty} =-h(\tau) < \infty,
\ee
and so $\tau$ is a carapoint for $\varphi$.
\end{remark}

\begin{question} Let  $d>2$, $\tau\in\mathbb{T}^{d}$ and let $h:\C_+^d(\tau) \to \C$ be  analytic function. Under which conditions on $h$, is $h$ a slope function of some function  $\varphi\in\mathcal{SA}_{d}$ at  a carapoint $\tau$?
\end{question}
 In the case $d=2$, the answer is given in \cite[Theorem 6.2]{ATY2}.

\section{An example of a function $\varphi_3\in\mathcal{SA}_{3}$ with a carapoint}\label{example-in-D3}

Note that Theorem \ref{directionalderivtheorem} can be used to find directional derivatives for Schur-Agler class functions that have a carapoint $\tau\in\mathbb{T}^{d}$ but need not necessarily be continuous at $\tau$. In this section, we shall apply this theorem to the function 
$$
\varphi_3(\lambda) = \dfrac{3\lambda_{1}\lambda_{2}\lambda_{3} - \lambda_{1}\lambda_{2} - \lambda_{1}\lambda_{3}-\lambda_{2}\lambda_{3}}{3-\lambda_{1}-\lambda_{2}-\lambda_{3}}\ \quad \text{for} \ \lam \in \D^3.
$$
We will need the following set.
\begin{definition}
    The symmetrized tridisc $\mathbb{G}_{3}$ is the set \index{$\mathbb{G}_{3}$}
$$\mathbb{G}_{3} =\Big\{(\lambda_{1}+\lambda_{2}+\lambda_{3},\lambda_{1}\lambda_{2}+\lambda_{2}\lambda_{3}+\lambda_{1}\lambda_{3},\lambda_{1}\lambda_{2}\lambda_{3}): \lambda_{1}\in\mathbb{D}, \lambda_{2}\in\mathbb{D}, \lambda_{3}\in\mathbb{D}\Big\}.$$
Alternatively, let $e_{j}$ denote the elementary symmetric function of degree $j$ in three variables for $j=1,2,3$. Thus
\begin{align}
    e_{1}(\lambda_{1},\lambda_{2},\lambda_{3}) &= \lambda_{1} + \lambda_{2}+\lambda_{3},\nonumber\\
    e_{2}(\lambda_{1},\lambda_{2},\lambda_{3}) &= \lambda_{2}\lambda_{3}+\lambda_{3}\lambda_{1}+\lambda_{1}\lambda_{2},\nonumber\\
    e_{3}(\lambda_{1},\lambda_{2},\lambda_{3}) &= \lambda_{1}\lambda_{2}\lambda_{3},\label{symmetrixfunctions}
\end{align}
and 
$$\mathbb{G}_{3}=\bigg\{(e_{1}(\lambda),e_{2}(\lambda),e_{3}(\lambda)):\lambda\in\mathbb{D}^{3}\bigg\}.$$
\end{definition}

\begin{example}\label{ddimensionalexample}
    Let $\varphi_3:\mathbb{D}^{3}\rightarrow \mathbb{C}$ be the function
\begin{equation}\label{kneesefunction}
    \varphi_3(\lambda) = \dfrac{3\lambda_{1}\lambda_{2}\lambda_{3}-\lambda_{1}\lambda_{2}-\lambda_{2}\lambda_{3}-\lambda_{1}\lambda_{3}}{3-\lambda_{1}-\lambda_{2}-\lambda_{3}}.
\end{equation}
The following conditions are satisfied:
    \begin{enumerate}
    \item $\varphi_3\in\mathcal{SA}_{3}$ and  $\varphi_3$ is inner;
\item $\mathds{1}=(1,1,1)\in\mathbb{T}^{3}$ is a carapoint of $\varphi_3$;
\item $\varphi_3$ does not extend continuously to a neighbourhood of the point $\mathds{1}$ and the cluster set of 
$\varphi_3$ at  $\mathds{1}$ is not a singleton;
\item there exists a desingularized model $(\mathcal{M},u, I)$ for $\varphi_3$ at $\mathds{1}$. Furthermore, $\varphi_3$ has a directional derivative $D_{-\delta}\varphi_3(\mathds{1})$ for all directions $\delta\in\mathbb{C}^{3}_{+}(\mathds{1})$ and
$$D_{-\delta}\varphi_3(\mathds{1}) = \varphi_3(\mathds{1})h(\delta)$$
where 
$$ h(\delta) = -\bigg\langle \bigg(\dfrac{\mathds{1}}{\mathds{1} \delta}\bigg)_{Y}^{-1}u(\mathds{1}),u(\mathds{1})\bigg\rangle_{\mathcal{M}}~~\text{for all}~ \delta\in\mathbb{C}_{+}^{3}(\mathds{1}).$$
    \end{enumerate}
\end{example}
\begin{proof} (1)
Recall that in \cite{Knese}, Knese proved that $\varphi_3\in\mathcal{SA}_{3}$, and by Rudin \cite[Theorem 5.2.5]{rudinpoly}, $\varphi_3$ is inner. 

(2) If $0 < r < 1$ and $\lam= r\mathds{1}$ then
$$\ph_3(\lam) =\varphi_3(r\mathds{1}) = \dfrac{3r^{3}-3r^{2}}{3-3r} = -r^{2},
$$
and so
\[
\frac{1-|\ph_3(\lam)|}{1-\|\lam\|_\infty} = \frac{1-r^2}{1-r} = 1+r \leq 2,
\]
thus
\[
\liminf_{\lam \to \mathds{1}} \frac{1-|\ph_3(\lam)|}{1-\|\lam\|_\infty}  < \infty,
\]
which is to say that $\mathds{1}$ is a carapoint for $\ph_3$. 
Observe also that $\varphi_3$ has radial limit $-1$ at $\mathds{1}$,
\begin{equation}\label{radiallimitddimenional}
    \lim_{r\rightarrow 1^{-}}\varphi_3(r\mathds{1})=-1.
\end{equation}

(3) Now suppose $\varphi_3$ does have a continuous extension $\psi$ to $\mathbb{D}^{3}\cup U$ for some open neighborhood $U$ of $\mathds{1}$. Then 
$$
\psi(\mathds{1}) = \lim_{r\rightarrow 1^{-}} \psi(r\mathds{1}) = \lim_{r\rightarrow 1^{-}}\varphi_3(r\mathds{1}) = -1,
$$
and since $\psi$ is continuous on $\mathbb{D}^{3}\cup U$ it follows that, for any path $(\lambda(t))_{(0<t<1)}$ in $\mathbb{D}^{3}$ such that $\lambda(t)\rightarrow \mathds{1}$ as $t\rightarrow 0$, $\varphi_3(\lambda(t)) = \psi(\lambda(t))\rightarrow \psi(\mathds{1}) = -1$. Hence, to show that $\varphi_3$ does not extend continuously to any neighbourhood of $\mathds{1}$, it suffices to construct a path $(\lambda(t))_{(0<t<1)}$ in $\mathbb{D}^{3}$, which need not be continuous, such that $\lambda(t)\rightarrow \mathds{1}$ as $t\rightarrow 0$ and $\varphi_3(\lambda(t))\nrightarrow -1$ as $t\rightarrow 0$.

We can write $\varphi_3$ in the following form
$$\varphi_3(\lambda) = \dfrac{3e_{3}(\lambda)-e_{2}(\lambda)}{3-e_{1}(\lambda)}~~\text{for}~\lambda\in\mathbb{D}^{3},$$
where $e_{j}$ for $j=1,2,3$ are defined by equations (\ref{symmetrixfunctions}).

If we define $f:\mathbb{G}_{3}\rightarrow \mathbb{C}$ by $f(s_{1},s_{2},s_{3}) = \dfrac{3s_{3}-s_{2}}{3-s_{1}}$ then $\varphi_3(\lambda) = f(e_{1}(\lambda),e_{2}(\lambda),e_{3}(\lambda))$ for all $\lambda\in\mathbb{D}^{3}$. Note also that $\lambda$ can be recovered (up to a permutation of its entries) from the triple $(e_{1}(\lambda),e_{2}(\lambda),e_{3}(\lambda))$, since $\lambda_{1},\lambda_{2},\lambda_{3}$ are the zeros of the polynomial
$$z^{3}-e_{1}(\lambda)z^{2}+e_{2}(\lambda)z-e_{3}(\lambda)$$
in $z$. We shall construct a path $s(t)$ such that $0<t<1$ in $\mathbb{G}_{3}$ by using a parametrization of the points of $\mathbb{G}_{3}$ due to Costara, \cite[Section 2.5]{Costara}. Note that, for any $\lambda\in\mathbb{C}^{3}$, $(e_{1}(\lambda),e_{2}(\lambda),e_{3}(\lambda))\in\mathbb{G}_{3}$ if and only if the roots of the equation
$$z^{3}-e_{1}(\lambda)z^{2}+e_{2}(\lambda)z-e_{3}(\lambda)=0$$
all lie in $\mathbb{D}$, that is, if and only if $\lambda\in\mathbb{D}^{3}$.

 From \cite{characterizationspoly}, a point $(s_{1},s_{2},s_{3})\in\mathbb{C}^{3}$ belongs to $\mathbb{G}_{3}$ if and only if $s_{3}\in\mathbb{D}$ and there exists $(z_{1},z_{2})\in\mathbb{G}$ such that 
\begin{equation}\label{eqnforgidentities}
    s_{1}=z_{1}+\overline{z}_{2}s_{3}~~\text{and}~~s_{2}=z_{2}+\overline{z_{1}}s_{3}.
\end{equation}
It is also known that $(z_{1},z_{2})\in\mathbb{G}$ if and only if $z_{2}\in\mathbb{D}$ and there exists $\beta\in\mathbb{D}$ such that $z_{1}=\beta+\overline{\beta}z_{2}$.

 For $t\in(0,1)$, choose a point $z(t)\in\mathbb{G}$ by taking $\beta=z_{2}=1-t$, that is, 
$$z_{1}(t) = 1-t+(\overline{1-t})(1-t)~~\text{and}~~z_{2}(t)=1-t.$$ 
Thus
$$z(t) = \bigg((1-t)(2-t),1-t\bigg).$$
Thence construct a point $s(t)\in\mathbb{G}_{3}$ by taking $(z_{1},z_{2})=z(t)$, $s_{3}=1-t$ in the parametrization (\ref{eqnforgidentities}). We obtain
\begin{align}
    s(t) &= \bigg((1-t)(2-t)+(1-t)^{2},1-t+\overline{(1-t)(2-t)}(1-t),(1-t)\bigg)\nonumber\\
    &=\bigg((1-t)(3-2t),(1-t)(1+(1-t)(2-t)),(1-t)\bigg)\label{identityfordis}
\end{align}
Now lift the path $s(t)$ such that $0<t<1$ in $\mathbb{G}_{3}$ to a path $\lambda(t)$ for $0<t<1$ in $\mathbb{D}^{3}$ as follows.

 Since $s(t)\in\mathbb{G}_{3}$, the three roots of the cubic equation 
\begin{equation}\label{iur983}
    z^{3}-e_{1}(\lambda)z^{2}+e_{2}(\lambda)z-e_{3}(\lambda)=0
\end{equation}
all lie in $\mathbb{D}$. For each $t\in(0,1)$, let $\lambda_{1}(t)$ be the root of the equation (\ref{iur983}) which has minimum distance from $-1$ and let $\lambda_{2}(t), \lambda_{3}(t)$ be the other two roots of equation (\ref{iur983}), numbered so that 
\begin{equation}\label{numberedinequality}
|1-\lambda_{1}(t)|\leq |1-\lambda_{2}(t)|\leq|1-\lambda_{3}(t)|.
\end{equation}
Then, for each $t$ and all $z\in\mathbb{C}$,
\begin{equation}\label{cubictings}
(z-\lambda_{1}(t))(z-\lambda_{2}(t))(z-\lambda_{3}(t))=z^{3}-s_{1}(t)z^{2}+s_{2}(t)z-s_{3}(t).
\end{equation}
Putting $z=1$, we have
$$(1-\lambda_{1}(t))(1-\lambda_{2}(t))(1-\lambda_{3}(t))=1-s_{1}(t)+s_{2}(t)-s_{3}(t).$$
Take absolute values and use the fact that we chose $\lambda_{1}(t)$ so that $|1-\lambda_{1}(t)|\leq |1-\lambda_{j}(t)|$ for $j=2,3$ to deduce that 
$$|1-\lambda_{1}(t)|^{3}\leq |1-\lambda_{1}(t)||1-\lambda_{2}(t)||1-\lambda_{3}(t)|=|1-s_{1}(t)+s_{2}(t)-s_{3}(t)|.$$
From equation (\ref{identityfordis}) we can see that $s(t)\rightarrow (3,3,1)$ as $t\rightarrow 0$ and hence $1-s_{1}(t)+s_{2}(t)-s_{3}(t)\rightarrow 0$. It follows that $\lambda_{1}(t)\rightarrow 1$ as $t\rightarrow 0$. 

In the light of equation (\ref{cubictings}), the function
$$q_{t}(z)=\dfrac{z^{3}-s_{1}(t)z^{2}+s_{2}(t)z-s_{3}(t)}{z-\lambda_{1}(t)}$$
is a monic quadratic polynomial in $z$ for each $t\in(0,1)$: let us write 
$$q_{t}(z)=z^{2}-b_{1}(t)z+b_{0}(t).$$
Then
$$(z-\lambda_{1}(t))(z^{2}-b_{1}(t)z+b_{0}(t)) = z^{3}-s_{1}(t)z^{2}+s_{2}(t)z-s_{3}(t).$$
Thus 
$$\lambda_{1}(t)+b_{1}(t)=s_{1}(t),~~\lambda_{1}(t)b_{1}(t)+b_{0}(t)=s_{2}(t),~~\lambda_{1}(t)b_{0}(t)=s_{2}(t).$$
As $t\rightarrow 0$, $\lambda_{1}(t)\rightarrow 1$, $s(t)\rightarrow (3,3,1)$, and so $b_{1}(t)\rightarrow 2$, $b_{0}(t)\rightarrow 1$. Therefore 
$$q_{t}(z)\rightarrow z^{2}-2z+1~~\text{as}~t\rightarrow 0~~\text{for all}~z\in\mathbb{C}.$$
Note that $\lambda_{2}(t),\lambda_{3}(t)$ are the zeros of the monic quadratic polynomial $q_{t}(z)$. We have numbered them so that $|1-\lambda_{2}(t)|\leq |1-\lambda_{3}(t)|$, see equation (\ref{numberedinequality}). Then, for all $t\in(0,1)$, 
$$|1-\lambda_{2}(t)|^{2}\leq |1-\lambda_{2}(t)||1-\lambda_{3}(t)|=|q_{t}(1)|,$$
and since $q_{t}(1)\rightarrow 0$ as $t\rightarrow 0$, it follows that 
$$\lambda_{2}(t)\rightarrow 1~~\text{as}~t\rightarrow 0.$$
Furthermore, 
$$\lambda_{1}(t)\lambda_{2}(t)\lambda_{3}(t) = s_{3}(t)$$
for all t and $s_{3}(t)\rightarrow 1$ as $t\rightarrow 0$, we have $\lambda_{3}(t)\rightarrow 1$ as $t\rightarrow 0$. Hence
$$\lambda(t)=(\lambda_{1}(t),\lambda_{2}(t),\lambda_{3}(t))\rightarrow (1,1,1)~~\text{as}~t\rightarrow 0.$$
Therefore,
\begin{align*}
    \varphi_3(\lambda(t))&=\dfrac{3s_{3}(t)-s_{2}(t)}{3-s_{1}(t)} \\
    &=\dfrac{3(1-t)-(1-t)(3-3t+t^{2})}{3-(3-5t+2t^{2})}\\
    &=\dfrac{(1-t)(3t-t^{2})}{5t-2t^{2}}=\dfrac{(1-t)(3-t)}{5-2t}\\
    &\rightarrow \dfrac{3}{5}~~\text{as}~t\rightarrow 0.
\end{align*}
Hence $\varphi_3$ has no continuous extension to any neighbourhood of $\mathds{1}$ and the cluster set of 
$\varphi_3$ at $\mathds{1}$ is not a singleton.

(4)  By Theorem \ref{generalizedmodeltheorem}, there exists a desingularized model $(\mathcal{M},u,I)$ for $\varphi_3$ at $\mathds{1}$. By Theorem \ref{directionalderivtheorem}, there exists an analytic function $h:\mathbb{C}^{3}_{+}(\mathds{1})\rightarrow\mathbb{C}$ given by
$$
h(z)=-\bigg\langle \bigg(\dfrac{\mathds{1}}{\mathds{1}z}\bigg)_{Y}u(\mathds{1}),u(\mathds{1})\bigg\rangle_{\mathcal{M}}
$$
for all $z\in\mathbb{C}^{3}_{+}(\mathds{1})$, for some $3$-tuple of positive contractions $Y=(Y_{1}, Y_{2},Y_{3})$ on $\mathcal{M}$ for which $Y_{1}+Y_{2}+Y_{3}=1_{\mathcal{M}}$ such that, for all $\delta\in\mathbb{C}^{3}(\mathds{1})$, $\varphi_3$ is differentiable in the direction $-\delta$ at the point $\mathds{1}$ with 
$$
D_{-\delta}\varphi_3(\mathds{1})=\varphi_3(\mathds{1})h(\delta).
$$
\end{proof}
\begin{remark}
\rm In \cite{Knese} Knese also proved that every Schur-Agler class rational inner function has a finite-dimensional unitary realization and that the function $\varphi_3$ given by equation (\ref{kneesefunction}) has a $10$-dimensional unitary realization. In \cite{Bickel-Knese}  Bickel and Knese, considered finite-dimensional unitary realizations for a class of functions $\varphi \in\mathcal{SA}_{3}$ and, in particular, they proved that the function $\varphi_3$ has a $7$-dimensional unitary realization and cannot be realized with dimension less than $7$, see \cite[Theorem 1.12]{Bickel-Knese}. 
\end{remark}

\section{A model and a unitary realization for $\varphi_3$}\label{model3}
In the dissertation of C. Evans (see \cite[Sections 7.3, 7.4]{Evans}) an explicit $9$-dimensional model $(\mathcal{L},v)$ of $\ph_3$ given by equation (\ref{kneesefunction}) is computed and the corresponding realization is presented. In this model $(\mathcal{L},v)$, $\mathds{1}$ is a $B$-point but not a $C$-point.

Recall the following theorem of  Knese \cite{Knese}.

\begin{theorem}{\normalfont{\cite[Theorem 2.9]{Knese}}}\label{kneseresultunitary}
    Let $f:\mathbb{D}^{d}\rightarrow \mathbb{D}$ be a Schur-Agler class rational inner function and write $f=p/q$, with $p,q\in\mathbb{C}[z_{1},\hdots,z_{d}]$.  Then
    \begin{enumerate}
        \item a sum-of-squares decomposition holds. There exist integers $N_{1},\hdots,N_{d}$ such that
        $$\lvert q(z)\rvert^{2} - \lvert p(z)\rvert^{2} = \sum_{j=1}^{d}(1-\lvert z_{j}\rvert^{2})\sum_{k=1}^{N_{j}}\lvert A_{j,k}(z)\rvert^{2}$$
    where $A_{j,k}\in\mathbb{C}[z_{1},\hdots,z_{d}]$.
    \item $f$ has a finite-dimensional unitary realization.
    \end{enumerate}   
\end{theorem}

We are going to use Knese's results from  \cite[Section 4]{Knese}, to get a model and a unitary realization of $\varphi_3$. Note that the function $\varphi_3$, defined above  by the equation
\eqref{kneesefunction},
can be presented as $\varphi_3(\lambda)=\dfrac{p(\lambda)}{q(\lambda)}$, for $\lambda=(\lambda_{1},\lambda_{2},\lambda_{3})\in\mathbb{D}^{3}$, where
\begin{align*}
    p(\lambda_{1},\lambda_{2},\lambda_{3}) &= 3\lambda_{1}\lambda_{2}\lambda_{3}-\lambda_{1}\lambda_{2}-\lambda_{1}\lambda_{3}-\lambda_{2}\lambda_{3}, \\
    q(\lambda_{1},\lambda_{2},\lambda_{3}) &= 3 - \lambda_{1} - \lambda_{2} -\lambda_{3}.
\end{align*}
Knese showed in \cite{Knese} that, for $\lambda=(\lambda_{1},\lambda_{2},\lambda_{3})\in\mathbb{D}^{3}$, 
\begin{align}
    \lvert &3 -\lambda_{1} - \lambda_{2} -\lambda_{3} \rvert^{2} -\lvert 3\lambda_{1}\lambda_{2}\lambda_{3}-\lambda_{1}\lambda_{2}-\lambda_{1}\lambda_{3}-\lambda_{2}\lambda_{3}\rvert^{2}\nonumber \\
    &= (1-\lvert \lambda_{1}\rvert^{2})S(\lambda_{2},\lambda_{3})+(1-\lvert \lambda_{2}\rvert^{2})S(\lambda_{1},\lambda_{3})+(1-\lvert \lambda_{3}\rvert^{2})S(\lambda_{1},\lambda_{2}).\label{kneselayout}
\end{align}
Here, for $\eta,\zeta\in\mathbb{C}$,
$$
S(\eta,\zeta) = \lvert s_{1}(\eta,\zeta)\rvert^{2}+\lvert s_{2}(\eta,\zeta) \rvert^{2}+\lvert s_{3}(\eta,\zeta)\rvert^{2},
$$
where
\begin{align*}
    s_{1}(\eta,\zeta)&=\sqrt{3}\Big(\eta \zeta-\dfrac{1}{2}\eta-\dfrac{1}{2}\zeta\Big)\\
    s_{2}(\eta,\zeta)&=\sqrt{3}\Big(1-\dfrac{1}{2}\eta-\dfrac{1}{2}\zeta\Big)\\
    s_{3}(\eta,\zeta)&=\dfrac{1}{\sqrt{2}}(\eta-\zeta).
\end{align*}
Define the map $v:\mathbb{D}^{3}\rightarrow\mathbb{C}^{9}$ by
\be\label{defv}
v(\lambda_{1},\lambda_{2},\lambda_{3})=\dfrac{1}{q(\lambda_{1},\lambda_{2},\lambda_{3})}\begin{bmatrix}
    w(\lambda_{2},\lambda_{3})\\
    w(\lambda_{1},\lambda_{3})\\
    w(\lambda_{1},\lambda_{2})
\end{bmatrix},
\ee
where, for $\eta,\zeta \in\D$,
$$
w(\eta,\zeta)=\begin{bmatrix}
    s_{1}(\eta,\zeta)\\
    s_{2}(\eta,\zeta)\\
    s_{3}(\eta,\zeta)
\end{bmatrix}.
$$
Note that, for $\eta,\zeta \in\D$,
\begin{align*}
    S(\eta,\zeta) &= \|w(\eta,\zeta)\|_{\mathbb{C}^{3}}^{2}.
\end{align*}
With respect to the decomposition $\mathbb{C}^{9}=\mathbb{C}^{3}\oplus\mathbb{C}^{3}\oplus\mathbb{C}^{3}$, define $P_{1}, P_{2}$ and $P_{3}$ to be the orthogonal projections onto $\mathbb{C}^{3}\oplus\{0\}\oplus\{0\}$, $\{0\}\oplus\mathbb{C}^{3}\oplus\{0\}$ and $\{0\}\oplus\{0\}\oplus\mathbb{C}^{3}$, respectively. For $\lambda=(\lambda_{1},\lambda_{2},\lambda_{3})\in\mathbb{D}^{3}$, the map $\lambda_{P}:\mathbb{C}^{9}\rightarrow\mathbb{C}^{9}$ is defined by
\beq\label{lambdaP}
\lambda_{P} = \lambda_{1}P_{1}+\lambda_{2}P_{2}+\lambda_{3}P_{3}.
\eeq 

\begin{proposition}{\normalfont{({\bf{Polarization theorem for holomorphic functions}}) \cite[Proposition 1]{DAngelo}}}\label{polarizationtheorem}
Suppose that $\Omega$ is a domain in $\mathbb{C}^{d}$ and let $\Omega^{*}=\{z\in\mathbb{C}^{d} : \overline{z}\in\Omega\}$ be its conjugate domain. Suppose that $H:\Omega \times \Omega^{*}\rightarrow \mathbb{C}$ is a holomorphic function of $2d$ complex variables $(\lambda,\mu)$, and that 
$$H(\lambda,\overline{\lambda})=0$$
for all $\lambda\in\Omega$. Then
$$H(\lambda,\mu)=0$$
for all $(\lambda,\mu)\in\Omega\times\Omega^{*}$.
\end{proposition}

We show that there exists a model $(\mathbb{C}^{9},v)$ for $\varphi_3$. For $\lambda=(\lambda_{1},\lambda_{2},\lambda_{3})\in\mathbb{D}^{3}$, rewrite equation (\ref{kneselayout}) as 
\begin{align*}
    \lvert q(\lambda) \rvert^{2} &-\lvert p(\lambda)\rvert^{2}= (1-\lvert \lambda_{1}\rvert^{2})\|w(\lambda_{2},\lambda_{3})\|_{\mathbb{C}^{3}}^{2}\nonumber \\
    &+(1-\lvert \lambda_{2}\rvert^{2})\|w(\lambda_{1},\lambda_{3})\|_{\mathbb{C}^{3}}^{2}+(1-\lvert \lambda_{3}\rvert^{2})\|w(\lambda_{1},\lambda_{2})\|_{\mathbb{C}^{3}}^{2},
\end{align*}
which is equivalent to
\begin{align}
    1 -\bigg\lvert &\dfrac{p(\lambda)}{ q(\lambda)}\bigg\rvert^{2}= (1-\lvert \lambda_{1}\rvert^{2})\bigg\|\dfrac{w(\lambda_{2},\lambda_{3})}{q(\lambda)}\bigg\|_{\mathbb{C}^{3}}^{2}+(1-\lvert \lambda_{2}\rvert^{2})\bigg\|\dfrac{w(\lambda_{1},\lambda_{3})}{q(\lambda)}\bigg\|_{\mathbb{C}^{3}}^{2}\nonumber \\
    &\hspace{2cm}+(1-\lvert \lambda_{3}\rvert^{2})\bigg\|\dfrac{w(\lambda_{1},\lambda_{2})}{q(\lambda)}\bigg\|_{\mathbb{C}^{3}}^{2}~~\text{for all}~\lambda=(\lambda_{1},\lambda_{2},\lambda_{3})\in\mathbb{D}^{3}.\label{kneselayoutnew}
\end{align}
By the polarization theorem for holomorphic functions, Proposition \ref{polarizationtheorem}, equation (\ref{kneselayoutnew}), we obtain, for $\lambda=(\lambda_{1},\lambda_{2},\lambda_{3})$ and $\mu=(\mu_{1},\mu_{2},\mu_{3})$ in $\mathbb{D}^{3}$,
\begin{align}
    & 1 -\dfrac{\overline{p(\mu)}}{\overline{q(\mu)}}\dfrac{p(\lambda)}{ q(\lambda)}\nonumber = (1-\overline{\mu_{1}}\lambda_{1})\bigg\langle\dfrac{w(\lambda_{2},\lambda_{3})}{q(\lambda)},\dfrac{w(\mu_{2},\mu_{3})}{q(\mu)}\bigg\rangle_{\mathbb{C}^{3}}\\
    &+( 1-\overline{\mu_{2}}\lambda_{2})\bigg\langle\dfrac{w(\lambda_{1},\lambda_{3})}{q(\lambda)},\dfrac{w(\mu_{1},\mu_{3})}{q(\mu)}\bigg\rangle_{\mathbb{C}^{3}}+(1-\overline{\mu_{3}} \lambda_{3})\bigg\langle\dfrac{w(\lambda_{1},\lambda_{2})}{q(\lambda)},\dfrac{w(\mu_{1},\mu_{2})}{q(\mu)}\bigg\rangle_{\mathbb{C}^{3}}.\label{kneselayoutnew2}
\end{align}
To see that the relation \eqref{kneselayoutnew2} holds, note that both sides of equation (\ref{kneselayoutnew2}) are holomorphic functions of $\lambda$ and $\overline{\mu}$ on $\mathbb{D}^{3} \times \mathbb{D}^{3}$, since $p$ and $q$ are polynomials, $q\neq 0$ in $\mathbb{D}^{3}$ and $v$ is a rational vector function without poles in $\mathbb{D}^{3}$. If $H(\lambda,\overline{\mu})$ denotes the difference between the two sides of equation (\ref{kneselayoutnew2}), then, by equation (\ref{kneselayoutnew}), $H(\lambda,\overline{\lambda})=0$ for all $\lambda\in\mathbb{D}^{3}$. Hence, by Proposition \ref{polarizationtheorem}, $H(\lambda,\overline{\mu})=0$ for all $\lambda,\mu\in\mathbb{D}^{3}$.

 We can rewrite the right hand side of equation (\ref{kneselayoutnew2}) as 
$$    \Bigg\langle \left(1_{\mathbb{C}^{9}}-\begin{bmatrix}
    \overline{\mu_{1}}\lambda_{1}P_{1} & & & \\
    &\overline{\mu_{2}}\lambda_{2}P_{2} & & \\
    & &\overline{\mu_{3}}\lambda_{3}P_{3}
\end{bmatrix} \right) \bigg(\dfrac{1}{q(\lambda)}\bigg)\begin{bmatrix}
    w(\lambda_{2},\lambda_{3})\\
    w(\lambda_{1},\lambda_{3})\\
    w(\lambda_{1},\lambda_{2})
\end{bmatrix},\bigg(\dfrac{1}{q(\mu)}\bigg)\begin{bmatrix}
    w(\mu_{2},\mu_{3})\\
    w(\mu_{1},\mu_{3})\\
    w(\mu_{1},\mu_{2})
\end{bmatrix}\Bigg\rangle_{\mathbb{C}^{9}}.$$
Hence, for $\lambda,\mu\in\mathbb{D}^{3}$, the right hand side of equation (\ref{kneselayoutnew2}) can be written in the form
$$   \langle (1_{\mathbb{C}^{9}}-\mu^{*}_{P}\lambda_{P})v(\lambda),v(\mu)\rangle_{\mathbb{C}^{9}},$$
where $v$ is given by equation \eqref{defv}.
Therefore $(\mathbb{C}^{9},v)$ is a model for $\varphi_3$, that is,
$$1-\overline{\varphi_3(\mu)}\varphi_3(\lambda) = \langle (1_{\mathbb{C}^{9}}-\mu^{*}_{P}\lambda_{P})v(\lambda),v(\mu)\rangle_{\mathbb{C}^{9}}~~\text{for all}~\lambda,\mu\in\mathbb{D}^{3}.$$
Note that, by Proposition \ref{prop5.4} and Example \ref{ddimensionalexample} (2), $\mathds{1}=(1,1,1)$ is a $B$-point for $(\mathbb{C}^{9},v)$.

 We also construct a finite-dimensional unitary realization $(a,\beta,\gamma,D)$ for $\varphi_3$, that is, we find  a scalar $a$, vectors $\beta,\gamma\in {\mathbb{C}^{9}}$ and an operator $D$ on ${\mathbb{C}^{9}}$ such that
\begin{equation}\label{unitaryopinknese}L = \begin{bmatrix}
    a & 1\otimes\beta\\
    \gamma\otimes 1 & D
\end{bmatrix}:\mathbb{C}\oplus \mathbb{C}^{9} \rightarrow\mathbb{C}\oplus\mathbb{C}^{9}
\end{equation}
is unitary, and, for all $\la\in\D^3$,
\beq\label{formphi-3}
\ph_3(\la) = a+ \langle \lambda_{P}(1_{\mathbb{C}^{9}}-D\lambda_{P})^{-1}\gamma,\beta\rangle_{\mathbb{C}^{9}}.
\eeq
Choose $a=\ph_3(0)=0$ and find a unitary realization with this choice. 
By the lurking isometry argument applied to the model $(\mathbb{C}^{9},v)$ of $\varphi_3$,
we get
\begin{align}
    \varphi_3(\lambda) &= \langle\lambda_{P}v(\lambda),\beta\rangle_{\mathbb{C}^{9}}\ \text {and} \ \label{varphiinknese}\\
    v(\lambda) &= D\lambda_{P}v(\lambda)+\gamma,\label{secondequationinknese}
\end{align}
for all $\la\in\D^3$.
To find $\beta$, suppose 
$$\beta=\begin{bmatrix}
    \beta_{1}\\
    \vdots\\
    \beta_{9}
\end{bmatrix}.$$
Since $\varphi_3(\lambda)=\dfrac{p(\lambda)}{q(\lambda)}$, for all $\la\in\D^3$, multiply both sides of equation (\ref{varphiinknese}) by $q(\lambda)$ to obtain 
$$p(\lambda) = \Bigg\langle\lambda_{1} w(\lambda_{2},\lambda_{3}),\begin{bmatrix}
    \beta_{1}\\
    \beta_{2}\\
    \beta_{3}
\end{bmatrix}\Bigg\rangle_{\mathbb{C}^{3}}+\Bigg\langle\lambda_{2} w(\lambda_{1},\lambda_{3}),\begin{bmatrix}
    \beta_{4}\\
    \beta_{5}\\
    \beta_{6}
\end{bmatrix}\Bigg\rangle_{\mathbb{C}^{3}}+\Bigg\langle\lambda_{3} w(\lambda_{1},\lambda_{2}),\begin{bmatrix}
    \beta_{7}\\
    \beta_{8}\\
    \beta_{9}
\end{bmatrix}\Bigg\rangle_{\mathbb{C}^{3}}.$$
Expand out each inner product in the equation above and compare each of the terms $\lambda_{1}\lambda_{2}\lambda_{3}$, $\lambda_{1}\lambda_{2}$, $\lambda_{1}\lambda_{3}$, $\lambda_{2}\lambda_{3}$, $\lambda_{1},\lambda_{2},\lambda_{3}$ and any scalars on the left and right hand side of the equation above to find
$$\begin{bmatrix}
    1 & 0 & 0 & 1 & 0 & 0 & 1 & 0 & 0\\
    -\sqrt{3} & 0 & \sqrt{2} & -\sqrt{3} & 0 & \sqrt{2} & 0&0&0\\
    -\sqrt{3} & 0 & -\sqrt{2} & 0&0&0&-\sqrt{3}&0&\sqrt{2}\\
    0&0&0&\sqrt{3} &0&\sqrt{2}&\sqrt{3}&0&\sqrt{2}
\end{bmatrix}\begin{bmatrix}
    \beta_{1}\\
    \vdots\\
    \beta_{9}
\end{bmatrix}=\begin{bmatrix}
    \sqrt{3}\\
    -2\\
    -2\\
    2
\end{bmatrix}.$$
It can be verified, for example by a computer algebra program such as Maple, that the unique solution of this system of equations is
\begin{equation}\label{beta}
\beta = \dfrac{1}{\sqrt{3}}\begin{bmatrix}
    1 \\
    0\\
    0\\
    1\\
    0\\
    0\\
    1\\
    0\\
    0
\end{bmatrix} = \dfrac{1}{\sqrt{3}}\begin{bmatrix}
    e_{1}\\
    e_{1}\\
    e_{1}
\end{bmatrix},
\end{equation}
where $e_{1}=(1,0,0)\in\mathbb{C}^{3}$.

 We now calculate $\gamma$. By Theorem \ref{kneseresultunitary}, the matrix $L$ from equation (\ref{unitaryopinknese}) with $a=0$ is unitary on $\mathbb{C}\oplus\mathbb{C}^{9}$ and satisfies
\begin{equation} 
L \begin{bmatrix}
    1 \\
    \lambda_{P}v(\lambda) 
\end{bmatrix} =
\begin{bmatrix}
    0 & 1\otimes\beta\\
    \gamma\otimes 1 & D
\end{bmatrix}\begin{bmatrix}
    1 \\
    \lambda_{P}v(\lambda) 
\end{bmatrix} = \begin{bmatrix}
    \varphi_3(\lambda)\\
    v(\lambda)
\end{bmatrix},\ \text{ for all} \ \lambda \in \D^3.
\end{equation}
On multiplying the above equation on the left by $L^{*}$, we obtain the equation
$$\begin{bmatrix}
    1 \\
    \lambda_{P}v(\lambda) 
\end{bmatrix} = \begin{bmatrix}
    0 & 1\otimes\gamma\\
    \beta\otimes 1 & D^{*}
\end{bmatrix} \begin{bmatrix}
    \varphi_3(\lambda)\\
    v(\lambda) 
\end{bmatrix}, \ \text{ for all} \ \lambda \in \D^3.$$
This is equivalent to the system of equations
\begin{align}
 1&=\langle v(\lambda), \gamma \rangle,\label{findinggammeeq}\\
 %   1&=(1\otimes\gamma)v(\lambda) =\langle v(\lambda), \gamma \rangle,\label{findinggammeeq}\\
    \lambda_{P}v(\lambda) &= \varphi_3(\lambda)\beta+D^{*}v(\lambda), \ \text{ for all} \ \lambda \in \D^3.\nonumber
\end{align}
Suppose 
$$\gamma=\begin{bmatrix}
    \gamma_{1}\\
    \vdots\\
    \gamma_{9}
\end{bmatrix}.$$
Multiply both sides of equation (\ref{findinggammeeq}) by $q(\lambda)$ to obtain 
\begin{align*}
    q(\lambda) &=  \Bigg\langle\begin{bmatrix}
    w(\lambda_{2},\lambda_{3})\\
    w(\lambda_{1},\lambda_{3})\\
    w(\lambda_{1},\lambda_{2})
\end{bmatrix},\begin{bmatrix}
    \gamma_{1}\\
    \vdots\\
    \gamma_{9}
\end{bmatrix}\Bigg\rangle_{\mathbb{C}^{9}}\\
&= \Bigg\langle
    w(\lambda_{2},\lambda_{3}),\begin{bmatrix}
    \gamma_{1}\\
    \gamma_{2}\\
    \gamma_{3}
\end{bmatrix}\Bigg\rangle_{\mathbb{C}^{9}}+\Bigg\langle
    w(\lambda_{1},\lambda_{3}),\begin{bmatrix}
    \gamma_{4}\\
    \gamma_{5}\\
    \gamma_{6}
\end{bmatrix}\Bigg\rangle_{\mathbb{C}^{9}}+\Bigg\langle
    w(\lambda_{1},\lambda_{2}),\begin{bmatrix}
    \gamma_{7}\\
    \gamma_{8}\\
    \gamma_{9}
\end{bmatrix}\Bigg\rangle_{\mathbb{C}^{9}},
\end{align*}
for all $\lambda \in \D^3$.
As in the method for finding $\beta$, we expand each inner product in the equation above and compare each of the terms $\lambda_{1}\lambda_{2}\lambda_{3}$, $\lambda_{1}\lambda_{2}$, $\lambda_{1}\lambda_{3}$, $\lambda_{2}\lambda_{3}$, $\lambda_{1},\lambda_{2},\lambda_{3}$ and any scalars on the left and right hand side to find
\begin{equation}\label{gamma}
\gamma= \dfrac{1}{\sqrt{3}}\begin{bmatrix}
    0 \\
    1\\
    0\\
    0\\
    1\\
    0\\
    0\\
    1\\
    0
\end{bmatrix}=\dfrac{1}{\sqrt{3}}\begin{bmatrix}
    e_{2}\\
    e_{2}\\
    e_{2}
\end{bmatrix},
\end{equation}
where $e_{2}=(0,1,0)\in\mathbb{C}^{3}$.

 To find $D$, rearrange equation (\ref{secondequationinknese}) in the following way.
\begin{align*}
    v(\lambda) - \gamma &= D\lambda_{P}v(\lambda)~~\text{for}~\lambda\in\mathbb{D}^{3}\\
\iff    \dfrac{1}{q(\lambda)}\begin{bmatrix}
    w(\lambda_{2},\lambda_{3})\\
    w(\lambda_{1},\lambda_{3})\\
    w(\lambda_{1},\lambda_{2})
\end{bmatrix} - \dfrac{1}{\sqrt{3}}\begin{bmatrix}
    e_{2}\\
    e_{2}\\
    e_{2}
\end{bmatrix} &= D\lambda_{P}\bigg(\dfrac{1}{q(\lambda)}\bigg)\begin{bmatrix}
    w(\lambda_{2},\lambda_{3})\\
    w(\lambda_{1},\lambda_{3})\\
    w(\lambda_{1},\lambda_{2})
\end{bmatrix}~~\text{for}~\lambda\in\mathbb{D}^{3}\\
\iff~\sqrt{3}\begin{bmatrix}
    w(\lambda_{2},\lambda_{3})\\
    w(\lambda_{1},\lambda_{3})\\
    w(\lambda_{1},\lambda_{2})
\end{bmatrix} - q(\lambda)\begin{bmatrix}
    e_{2}\\
    e_{2}\\
    e_{2}
\end{bmatrix} &= \sqrt{3}D\lambda_{P}\begin{bmatrix}
    w(\lambda_{2},\lambda_{3})\\
    w(\lambda_{1},\lambda_{3})\\
    w(\lambda_{1},\lambda_{2})
\end{bmatrix}~~\text{for}~\lambda\in\mathbb{D}^{3}.
\end{align*}
Write $D$ as a $9\times9$ matrix with entries $d_{ij}$ for $i,j=1,\hdots,9$. Once more, expand out the equation above and compare each of the terms $\lambda_{1}\lambda_{2}\lambda_{3}$, $\lambda_{1}\lambda_{2}$, $\lambda_{1}\lambda_{3}$, $\lambda_{2}\lambda_{3}$, $\lambda_{1},\lambda_{2},\lambda_{3}$ and any scalars on the left and right hand side to find
\begin{equation}\label{D}
D = \begin{bmatrix}
\dfrac{1}{3} & 0 & 0 &-\dfrac{1}{6} & -\dfrac{1}{2} & -\dfrac{1}{\sqrt{6}} & -\dfrac{1}{6} & -\dfrac{1}{2} & -\dfrac{1}{\sqrt{6}} \\[0.5cm]
0 & \dfrac{1}{3} & \dfrac{1}{3\sqrt{6}} & \dfrac{1}{6} & -\dfrac{1}{6} & \dfrac{2}{3\sqrt{6}} & -\dfrac{1}{6} & -\dfrac{1}{6} & \dfrac{1}{3\sqrt{6}}\\[0.5cm]
0 & 0 & \dfrac{2}{9} & -\dfrac{1}{3\sqrt{6}} & \dfrac{1}{\sqrt{6}} & \dfrac{1}{9} & \dfrac{1}{3\sqrt{6}} & -\dfrac{1}{\sqrt{6}} & -\dfrac{1}{9} \\[0.5cm]
-\dfrac{1}{6} & -\dfrac{1}{2} & -\dfrac{1}{\sqrt{6}} & \dfrac{1}{3} & 0 & 0 & -\dfrac{1}{6} & -\dfrac{1}{2} & \dfrac{1}{\sqrt{6}} \\[0.5cm]
\dfrac{1}{18} & -\dfrac{1}{6} & \dfrac{1}{3\sqrt{6}} & -\dfrac{1}{9} & \dfrac{1}{3} & 0 & \dfrac{1}{18} & -\dfrac{1}{6} & -\dfrac{1}{3\sqrt{6}} \\[0.5cm]
-\dfrac{3}{5\sqrt{6}} & \dfrac{1}{\sqrt{6}} & \dfrac{1}{5} & 0 & 0 & 0 & \dfrac{3}{5\sqrt{6}} & -\dfrac{1}{\sqrt{6}} & \dfrac{1}{5} \\[0.5cm]
-\dfrac{1}{6} & -\dfrac{1}{2} & \dfrac{1}{\sqrt{6}} & -\dfrac{1}{6} & -\dfrac{1}{2} & \dfrac{1}{\sqrt{6}} & \dfrac{1}{3} & 0 & 0\\[0.5cm]
\dfrac{1}{18} & -\dfrac{1}{6} & -\dfrac{1}{3\sqrt{6}} & \dfrac{1}{18} & -\dfrac{1}{6} & -\dfrac{1}{3\sqrt{6}} & -\dfrac{1}{9} & \dfrac{1}{3} & 0 \\[0.5cm]
-\dfrac{1}{3\sqrt{6}} & \dfrac{1}{\sqrt{6}} & -\dfrac{1}{9} & \dfrac{1}{3\sqrt{6}} & -\dfrac{1}{\sqrt{6}} & \dfrac{1}{9} & 0 & 0 & \dfrac{2}{9} 
\end{bmatrix}.
\end{equation} 
Hence $\varphi_3$, given by equation (\ref{kneesefunction}), can be written in the following form
$$\varphi_3(\lambda) = a + \langle \lambda_{P}(1_{\mathbb{C}^{9}}-D\lambda_{P})^{-1}\gamma,\beta \rangle_{\mathbb{C}^{9}},~~\text{for all}~\lambda\in\mathbb{D}^{3},$$
where $a=0$ and   $\beta$, $\gamma$ and $D$ are defined by equations  \eqref{beta}, \eqref{gamma} and \eqref{D} respectively and $\lambda_{P}$ is defined by \eqref{lambdaP}.
This can be verified by the computer program Maple.


\begin{thebibliography}{99}
\bibitem{Abate} M. Abate, The Julia-Wolff-Carath\'{e}odory theorem in polydiscs. {\it{J. Anal. Math.}} {\bf 74} (1998) 275-306.

\bibitem{Ag1} J. Agler, On the representation of certain holomorphic functions defined on a polydisc, In {\em Operator Theory: Advances and Applications, {\rm Vol. 48}}, pages
  47--66. {Birkh\"auser}, Basel, 1990.
  
\bibitem{AMYcara}  J. Agler, J. E. McCarthy and N. J. Young,  A Carath\'eodory Theorem for the bidisk via Hilbert space methods,  {\em Math. Annalen} {\bf 352} (2012) 581-624. Corrected version, arXiv:1002.3727 [math.CV], June 2026. 

\bibitem{AMYOperatorbook} J. Agler, J. E. McCarthy and N.J. Young, {\it Operator Analysis: Hilbert space methods in complex analysis}. Cambridge: Cambridge University Press, 2020.

\bibitem{ATY2} J. Agler, R. Tully-Doyle and N. J. Young,  Boundary behavior of analytic functions of two variables via generalized models, {\em Indagationes Math.} {\bf 23}(4) (2012) 995-1025.

\bibitem{TheoreticalNumericalAnalysis} K. Atkinson and W. Han. {\it Theoretical numerical analysis, a functional analysis framework}. Springer, 2009.

\bibitem{BallandBolot1} J. A. Ball and V. Bolotnikov, A tangential interpolation problem on the distinguished boundary of the polydisc for the Schur-Agler class. {\it J. Math. Anal. Appl.} {\bf 273}(2)
(2002) 328-348.

\bibitem{Bickel-Knese} K. Bickel and G. Knese, Inner functions on the bidisk and associated Hilbert spaces, 
 {\it J. Func. Anal.} {\bf 265} (2013) 2753-2790.
 
\bibitem{completenormedalgebras} F. F. Bonsall and J. Duncan, {\it Complete normed algebras}. Springer-Verlag, 1973.

\bibitem{caratheodory} C. Carath\'{e}odory, \"{U}ber die Winkelderivierten von beschr\"{a}nkten analytischen Funktionen, {\it Sitzunber. Preuss. Akad. Wiss} (1929) 39-52.

\bibitem{Costara} C. Costara, {\it{Le probleme de Nevanlinna-Pick spectral}}. Thesis submitted for the degree of Doctor of Philosophy, Laval University, 2004.

\bibitem{crabbdavie} M. J. Crabb and A. M. Davie, Von Neumann's inequality for Hilbert space operators. {\it{Bull. London Math. Soc.}}  {\bf 7} (1975) 49-50.

\bibitem{Evans} C. Evans, {\em Generalized models in the function
theory of several complex variables}, PhD Thesis, Newcastle University, UK, 2025.

\bibitem{DAngelo} J.P. D’Angelo. {\em{Several Complex Variables and the Geometry of Real Hypersurfaces}}. CRC, Boca Raton, FL, 1993.

\bibitem{jaf93} F.~Jafari,  Angular derivatives in polydisks,
{\em Indian J. Math.} {\bf 35} (1993) 197--212.

\bibitem{FoiasNagy} C. Foias and B. Sz.-Nagy, {\it{Harmonic analysis of operators on Hilbert space}.} North-Holland Publishing Company, Amsterdam - London, 1970.

\bibitem{characterizationspoly} S. Gorai and J. Sarkar, Characterizations of symmetrized polydisc. {\it{Indian J. Pure Appl. Math.}} {\bf 47} (2016) 391–397.
 
\bibitem{Knese} G. Knese, Rational inner functions in the Schur-Agler class of the polydisk, {\em Publicacions Matem\'{a}tiques}, {\bf 55}(2) (2011) 343-257.

\bibitem{rudinpoly} W. Rudin, {\it{Function theory in polydiscs}}. W.A. Benjamin Inc, New York, 1969.

\bibitem{RynneYoungson} B. P. Rynne and M. A. Youngson, {\it Linear Functional Analysis}. 2nd Edition, Springer, 2010.

\bibitem{rudin}
W. Rudin,  {\em Function theory in the unit ball of $\C^n$},  Springer, New York, 1980.

\bibitem{BarrySimon} B. Simon, {\it{Real Analysis, A comprehensive course in analysis, Part 1}}, American Mathematical Society, 2015.

\bibitem{Varopouloscite} N. Varopoulos, On an inequality of von Neumann and an application of the metric theory of tensor products to operator theory, {\it J. Func. Anal.} {\bf 16} (1974) 83-100.

\bibitem{woody}
K. W\l odarczyk,  Julia's Lemma and Wolff's Theorem for $J^*$-algebras,
 {\em Proc. Amer. Math. Soc.} {\bf 99} (1987) 472-476.
\end{thebibliography}
\end{document}